\documentclass[11pt, hidelinks]{article}

\usepackage{etoolbox}

\usepackage{hyperref}
\usepackage{url}
\usepackage[margin=1in]{geometry}
\usepackage{amsmath,amssymb,amsthm, amsfonts}
\usepackage[T1]{fontenc}
\usepackage{mathtools}
\usepackage{xcolor}
\usepackage{enumitem, comment, xifthen}
\usepackage{graphicx}
\usepackage{mathabx} %
\graphicspath{{figs/}}
\DeclarePairedDelimiter{\vecnorm}{\lVert}{\rVert}
\DeclarePairedDelimiter{\matnorm}{\vvvert}{\vvvert}
\newcommand{\twonorm}[2][]{\vecnorm[#1]{#2}_{2}}
\newcommand{\opnorm}[2][]{\matnorm[#1]{#2}_{\rm op}}
\newcommand{\fronorm}[2][]{\matnorm[#1]{#2}_{\rm F}}

\makeatletter
\def\letterdef#1#2#3{\def\letterdef@##1{\expandafter\def\csname #1\endcsname{#2}}%
  \letterdef@@#3{?\@car{}}\@nil}
\def\letterdef@@#1{\@gobble#1\letterdef@{#1}\letterdef@@}
\makeatother

\DeclarePairedDelimiterX{\klx}[2]{(}{)}{%
  #1\;\delimsize\|\;#2%
}
\DeclarePairedDelimiterX{\quantklx}[3]{(}{)}{%
  #1\;\delimsize\|\;#2\;\delimsize\vert\;#3%
}
\DeclarePairedDelimiterX{\inner}[2]{\langle}{\rangle}{%
  #1,#2%
}

\newcommand{\R}{\mathbf R} %
\newcommand{\N}{\mathbf N} %

\letterdef{c#1}{\mathcal{#1}}{ABCDEFGHIJKLMNOPQRSTUVWXYZ}
\letterdef{b#1}{\mathbb{#1}}{ABCDEFGHIJKLMNOPQRSTUVWXYZ}
\let\defn\coloneq
\newcommand{\twomax}[2]{\ensuremath{#1 \lor #2}}
\newcommand{\twomin}[2]{\ensuremath{#1 \land #2}}
\newcommand{\ceil}[1]{\left\lceil #1 \right\rceil}

\newcommand{\ud}[0]{\mathrm{d}}  %
\newcommand{\1}{\mathbf 1} %
\let\ones\1
 
\let\epsilon\varepsilon
\newcommand{\eps}{\varepsilon}

\let\tilde\widetilde

\let\succeq\succcurlyeq
\let\preceq\preccurlyeq

\renewcommand{\leq}{\leqslant}
\renewcommand{\geq}{\geqslant}

\newcommand{\argmin}{\mathop{\rm arg\,min}}

\DeclareMathOperator{\trace}{\bf Tr}
\DeclareMathOperator{\diag}{\bf diag}

\newcommand{\T}{\mathsf{T}}

\newcommand{\iid}{\textnormal{i.i.d.}}

\newcommand{\simiid}{\stackrel{\iid}{\sim}} %

\newcommand{\Cov}{\mathrm{Cov}}

\makeatletter
\newcommand{\E}{\operatorname*{\mathbf{E}}\ilimits@}
\makeatother
\makeatletter
\renewcommand{\P}{\operatorname*{\mathbf{P}}\ilimits@}
\makeatother

\newcommand{\ie}{\textit{i}.\textit{e}., }

\newcounter{algorithmctr}
\renewcommand{\thealgorithmctr}{\arabic{algorithmctr}}
   {\refstepcounter{algorithmctr}\begin{list}{}{%
       \setlength{\rightmargin}{0\linewidth}%
       \setlength{\leftmargin}{0\linewidth}}%
       \rmfamily\small
       \item[]{\setlength{\parskip}{0ex}\hrulefill\par%
        \nopagebreak{\bfseries\textsf{Algorithm \thealgorithmctr~}}}}%
   {{\setlength{\parskip}{-1ex}\nopagebreak\par\hrulefill} \end{list}}

\makeatletter
\long\def\@makecaption#1#2{
        \vskip 0.8ex
        \setbox\@tempboxa\hbox{\small {\bf #1.} #2}
        \parindent 1.5em 
        \dimen0=\hsize
        \advance\dimen0 by -3em
        \ifdim \wd\@tempboxa >\dimen0
                \hbox to \hsize{
                        \parindent 0em
                        \hfil 
                        \parbox{\dimen0}{\def\baselinestretch{0.96}\small
                                {\bf #1.} #2
                                } 
                        \hfil}
        \else \hbox to \hsize{\hfil \box\@tempboxa \hfil}
        \fi
        }
\makeatother

\usepackage[noabbrev]{cleveref}

\theoremstyle{plain}

\newtheorem{theo}{Theorem}

\newtheorem{lem}{Lemma}
\newtheorem{prop}{Proposition}
\newtheorem{cor}{Corollary}

\theoremstyle{definition} 

\newtheorem{nota}{Notation}
\newtheorem{de}{Definition}
\newtheorem{exa}{Example}
\AtBeginEnvironment{exa}{%
  \pushQED{\qed}%
}
\AtEndEnvironment{exa}{\popQED\endexa}

\newtheorem{as}{Assumption}
\newtheorem{alg}{Algorithm}

\Crefname{theo}{Theorem}{Theorems}
\Crefname{lem}{Lemma}{Lemmas}
\Crefname{prop}{Proposition}{Propositions}

\newcommand{\btheo}{\begin{theo}}
\newcommand{\bde}{\begin{de}}
\newcommand{\ble}{\begin{lem}}
\newcommand{\bpr}{\begin{prop}}
\newcommand{\bno}{\begin{nota}}
\newcommand{\bex}{\begin{exa}}
\newcommand{\bcor}{\begin{cor}}
\newcommand{\spro}{\begin{proof}}
\newcommand{\bas}{\begin{as}}
\newcommand{\balg}{\begin{alg}}

\newcommand{\etheo}{\end{theo}}
\newcommand{\ede}{\end{de}}
\newcommand{\ele}{\end{lem}}
\newcommand{\epr}{\end{prop}}
\newcommand{\eno}{\end{nota}}
\newcommand{\eex}{\end{exa}}
\newcommand{\ecor}{\end{cor}}
\newcommand{\fpro}{\end{proof}}
\newcommand{\eas}{\end{as}}
\newcommand{\ealg}{\end{alg}}

\theoremstyle{plain}

\newtheorem{theos}{Theorem}
\newtheorem{props}{Proposition}
\newtheorem{lems}{Lemma}
\newtheorem{cors}{Corollary}

\theoremstyle{definition}
\newtheorem{exas}{Example}
\newtheorem{algs}{Algorithm}
\newtheorem{asss}{Assumption}
\newtheorem{defns}{Definition}

\newcommand{\btheos}{\begin{theos}}
\newcommand{\etheos}{\end{theos}}
\newcommand{\bprops}{\begin{props}}
\newcommand{\eprops}{\end{props}}
\newcommand{\bdes}{\begin{defns}}
\newcommand{\edes}{\end{defns}}
\newcommand{\blems}{\begin{lems}}
\newcommand{\elems}{\end{lems}}
\newcommand{\bcors}{\begin{cors}}
\newcommand{\ecors}{\end{cors}}
\newcommand{\bexs}{\begin{exas}}
\newcommand{\eexs}{\end{exas}}
\newcommand{\balgs}{\begin{algs}}
\newcommand{\ealgs}{\end{algs}}
\newcommand{\bass}{\begin{asss}}
\newcommand{\eass}{\end{asss}}

\usepackage{hyperref}
\usepackage[margin=1in]{geometry}
\usepackage{amsmath,amssymb,amsthm, amsfonts}
\usepackage[T1]{fontenc}
\usepackage{mathtools}
\usepackage{xcolor}
\usepackage{comment, xifthen}
\usepackage{graphicx}
\usepackage{mathscinet}
\graphicspath{{figs/}}
\usepackage{subcaption} 
\makeatletter
\long\def\@makecaption#1#2{
        \vskip 0.8ex
        \setbox\@tempboxa\hbox{\small {\bf #1.} #2}
        \parindent 1.5em 
        \dimen0=\hsize
        \advance\dimen0 by -3em
        \ifdim \wd\@tempboxa >\dimen0
                \hbox to \hsize{
                        \parindent 0em
                        \hfil 
                        \parbox{\dimen0}{\def\baselinestretch{0.96}\small
                                {\bf #1.} #2
                                } 
                        \hfil}
        \else \hbox to \hsize{\hfil \box\@tempboxa \hfil}
        \fi
        }
\makeatother

\usepackage{xparse}
\let\realItem\item %
\makeatletter
\NewDocumentCommand\myItem{ o }{%
   \IfNoValueTF{#1}%
      {\realItem}%
      {\realItem[#1]\def\@currentlabel{#1}}%
}
\makeatother

\usepackage{enumitem}    
\setlist[enumerate]{
    before=\let\item\myItem,       %
    label=\textnormal{(\arabic*)}, %
    widest=(2')                    %
}
\newcommand{\dimension}{d}

\newcommand{\radius}{\varrho}

\newcommand{\numobs}{n}
\newcommand{\thetastar}{\theta^\star}

\newcommand{\CovMat}{\Sigma}
\newcommand{\ConstraintMat}{{K_{c}}}

\newcommand{\EstimationMat}{{K_{e}}}

\newcommand{\EmpCov}{\CovMat_\numobs}

\newcommand{\minimaxrisk}{\mathfrak{M}}
\newcommand{\Normal}[2]{\mathsf{N}\left(#1, #2\right)}

\newcommand{\ball}[1]{\mathsf{B}_\cH(#1)}
\newcommand{\inprod}[2]{\ensuremath{\langle #1 , \, #2 \rangle}}
\newcommand{\Information}[1]{\cI\left(#1\right)}
\DeclareMathOperator{\range}{\bf ran}

\newcommand{\opdim}{n}
\newcommand{\NoiseVariance}{\Sigma_w}
\renewcommand{\iid}{IID}
\let\IID\iid
\renewcommand{\simiid}{\stackrel{\rm \iid}{\sim}}
\let\NoiseCovariance\NoiseVariance

\newcommand{\PlainLinOp}{\ensuremath{T}}
\newcommand{\LinOp}[1]{\ensuremath{\PlainLinOp_{#1}}}

\newcommand{\myhat}[1]{\ensuremath{\widehat{#1}}}
\newcommand{\etahat}{\ensuremath{\myhat{\eta}}}
\newcommand{\thetahat}{\ensuremath{\myhat{\theta}}}

\newcommand{\myparagraph}[1]{\paragraph{#1}}

\newcommand{\fstar}{\ensuremath{f^\star}}

\newcommand{\usedim}{\ensuremath{d}}

\newcommand{\Xspace}{\ensuremath{\mathcal{X}}}

\newcommand{\Exs}{\E}
\newcommand{\real}{\R}
\newcommand{\Nat}{\N}

\newcommand{\tomowindow}{\ensuremath{h}}
\newcommand{\XiSpace}{\ensuremath{\Xi}}
\newcommand{\yhat}{\ensuremath{\widehat{y}}}
\newcommand{\jitter}{\ensuremath{u}}

\theoremstyle{plain}

\newtheorem{theorem}{Theorem}
\newtheorem{proposition}{Proposition}
\newtheorem{lemma}{Lemma}

\theoremstyle{definition}
\newtheorem{example}{Example}

\begin{document}

\begin{center}
  {\bf{\LARGE Noisy recovery from random linear observations: \\ Sharp
      minimax rates under elliptical constraints}} \\
  \vspace*{.2in}

{\large{
 \begin{tabular}{ccc}
  Reese Pathak$^{\diamond}$ & 
  Martin J.\ Wainwright$^{\diamond, \dagger, \star}$  &
  Lin Xiao$^{\ddagger}$  \\
  \texttt{pathakr@berkeley.edu} & 
  \texttt{wainwrigwork@gmail.com} &
  \texttt{linx@meta.com}
 \end{tabular}
}

\vspace*{.2in}

\begin{tabular}{c}
  EECS$^\diamond$ \& 
  Statistics$^\dagger$ \\ 
  UC Berkeley \\
 \end{tabular}
 
\medskip
 
 \begin{tabular}{c}
   EECS \& Mathematics \\
   Laboratory for Information and Decision Systems$^\star$ \\
   Statistics and Data Science Center$^\star$ \\
    Massachusetts Institute of Technology
 \end{tabular}

 \medskip
  \begin{tabular}{c}
Meta AI$^\ddagger$
 \end{tabular}
}
  \vspace*{.2in}
\end{center}

\begin{abstract}
 Estimation problems with constrained parameter spaces arise in
 various settings. In many of these problems, the observations
 available to the statistician can be modelled as arising from the
 noisy realization of the image of a random linear operator; an
 important special case is random design regression. We
 derive sharp rates of estimation for arbitrary compact elliptical
 parameter sets and demonstrate how they depend on the distribution of
 the random linear operator.  Our main result is a
 functional that characterizes the minimax rate of estimation in terms
 of the noise level, the law of the random operator, and elliptical
 norms that define the error metric and the parameter space.  This
 nonasymptotic result is sharp up to an explicit universal constant,
 and it becomes asymptotically exact as the radius of the parameter
 space is allowed to grow. We demonstrate the generality of the result
 by applying it to both parametric and nonparametric regression
 problems, including those involving distribution shift
 or dependent covariates.
\end{abstract}

\section{Introduction}

In this paper, we study the problem of estimating an unknown vector 
$\thetastar$ on the basis of random linear observations corrupted by
noise. More concretely, suppose that we observe a random operator
$\LinOp{\xi}$ and a random vector $y$, which are linked via the
equation
\begin{equation}
\label{eqn:operator-model}
y = \LinOp{\xi}(\thetastar) + w.
\end{equation}
This observation model involves two forms of randomness: the
unobserved vector $w$, which is a form of additive observation
noise, and the observed operator $\LinOp{\xi}$, which is random, as
indicated by its dependence on an underlying random variable $\xi$.

While relatively simple in appearance, the observation
model~\eqref{eqn:operator-model} captures a broad range of statistical
estimation problems.  

\begin{example}[Linear regression]
\label{ex:linear-regression}
We begin with a simple but widely used model: linear regression. The
goal is to estimate the coefficients $\thetastar \in \real^\usedim$
that define the best linear predictor $x \mapsto
\inprod{x}{\thetastar}$ of some real-valued response variable $Y \in
\real$.  In order to do so, we observe a collection of $(x_i, y_i)$
pairs linked via the noisy observation model
\begin{align*}
y_i & = \inprod{x_i}{\thetastar} + w_i \qquad \mbox{for $i = 1,
  \ldots, \numobs$.}
\end{align*}
If we define the concatenated vector $y = (y_1, \ldots, y_\numobs)$,
with an analogous definition for $w$, this is a special case of our
general setup with the random linear operator $\LinOp{\xi}:\R^d
\rightarrow \R^\numobs$ given by
\begin{align}
[\LinOp{\xi}(\theta)]_i & = \inprod{x_i}{\theta} \quad \mbox{for $i =
  1, \ldots, \numobs$.}
\end{align}
Here, the random index corresponds to the covariate vectors
so that $\xi = (x_1, \dots, x_\numobs)$; note that we
have imposed no assumptions on the dependence structure
of these covariate vectors.  In the classical
setting, these covariates are assumed to be drawn in an i.i.d.
manner; however, our general set-up is by no means limited to this
classical setting.  In the sequel, we consider various examples with
interesting dependence structure, and our theory gives some very
precise insights into the effects of such dependence.
\end{example}

\begin{example}[Nonparametric regression]
\label{ex:nonparametric-regression}
In the preceding example, we discussed the problem of predicting a
response variable $Y \in \real$ in a linear manner.  Let us consider
the nonparametric generalization: here our goal is to estimate the
regression function $\fstar(x) \defn \E[Y \mid X = x]$, which need not
be linear as a function of $x$.  Given observations $\{(x_i, y_i)
\}_{i=1}^\numobs$, we can write them in the form
\begin{align*}
y_i & = \fstar(x_i) + w_i, \qquad \mbox{for $i = 1, \ldots, \numobs$,}
\end{align*}
where $w_i = y_i - \Exs[Y \mid X = x_i]$ are zero-mean noise
variables.

Now let us suppose that $\fstar$ belongs to some function class $\cF$
contained with $L^2(\Xspace)$, and show how this observation model can
be understand as a special case of our setup with $\thetastar \in
\ell^2(\Nat)$.  Take some orthonormal basis $\{\phi_j\}_{j \geq 1}$ of
$L^2(\Xspace)$.  Any function in $\cF$ can then be expanded as $f =
\sum_{j \geq 1} \theta_j \phi_j$ for some sequence $\theta \in
\ell^2(\Nat)$.  Letting $\xi = (x_1, \ldots, x_\numobs)$, we can
define the operator $\LinOp{\xi}: \ell^2(\Nat) \rightarrow \R^\numobs$
via
\begin{align*}
\theta \mapsto [\LinOp{\xi}(\theta)]_i & \defn \sum_{j=1}^\infty \theta_j
\phi_j(x_i) \quad \mbox{for $i = 1, \ldots, \numobs$},
\end{align*}
so that this problem can be written in the form of our general
model~\eqref{eqn:operator-model}.  Observe that the randomness in the
observation operator $\LinOp{\xi}$ arises via the randomness in
sampling the covariates $\{x_i\}_{i=1}^\numobs$.
\end{example}

\begin{example}[Tomographic reconstruction]
The problem of tomographic reconstruction refers to the problem of
recovering an image, modeled as a real-valued function $\fstar$ on
some compact domain $\Xspace \subset \real^2$, based on noisy integral
measurements.  Formally, we observe responses of the form
\begin{align*}
y_i & = \int_{\Xspace} \tomowindow(x_i, u) \fstar(u) \, \ud u + w_i \qquad
\mbox{for $i = 1, \ldots, \numobs$,}
\end{align*}
where $\tomowindow: \real^2 \times \real^2 \rightarrow \real$ is a
known window function.  If we again view $\fstar$ as belonging to some
function class $\cF$ within $L^2(\Xspace)$, then we can write this
model in our general form with
\begin{align*}
[\LinOp{\xi}(v)]_i & = \sum_{j \geq 1} v_j \Big[ \int_{\Xspace}
  \tomowindow(x_i, u) \phi_j(u) \, \ud u \Big], \quad \mbox{and $\xi = (x_1,
  \ldots, x_\numobs)$.}
\end{align*}
Here we have followed the same conversion as
in~\Cref{ex:nonparametric-regression}, in particular re-expressing
$\fstar$ in terms of its generalized Fourier coefficients with respect
to an orthonormal family $\{\phi_j \}_{j \geq 1}$.
\end{example}

\begin{example}[Error-in-variables]
Consider the Berkson variant~\cite{Ber50,CarEtal95} of the
error-in-variables problem in nonparametric regression.  In this
problem, an observed covariate $x$---instead of being associated with
a noisy observation of $\fstar(x)$---is associated with a noisy
observation of the ``jittered'' evaluation $\fstar(x + \jitter)$,
where $\jitter \in \R$ is the random jitter.  Formally, we observe
$\numobs$ pairs $(x_i, y_i)$ of the form
\begin{align*}
y_i & = \fstar(x_i + \jitter_i) + \epsilon_i \qquad \mbox{for $i = 1,
  \ldots, \numobs$,}
\end{align*}
where the unobserved random jitter $\jitter_i$ is drawn independently
of the pair $(x_i, \epsilon_i)$.  We can re-write these observations
as a special case of our general model with $\xi = (x_1, \ldots,
x_\numobs)$, and
\begin{align*}
  [\LinOp{\xi}(f)]_i \defn \Exs_\jitter \big[f(x_i + \jitter) \big],
  \quad \mbox{and} \quad w_i \defn \epsilon_i + \Big \{ f(x_i +
  \jitter_i) - \Exs_\jitter \big[f(x_i + \jitter) \big] \Big \} \quad
  \mbox{for $i = 1, \ldots, \numobs$.}
\end{align*}
Note that the new noise variables $w_i$ are again zero-mean, and our
assumption that $\LinOp{\xi}$ is observed means that the distribution
of the jitter $\jitter$ is known.
\end{example}

These examples (and others, as discussed below
in~\Cref{sec:examples-high-level}) motivate our study of the operator
model~\eqref{eqn:operator-model}.  As we discuss in further detail
later, a key advantage of writing the observation model in this form
is that it will allow us to separate three key components of the
difficulty of the problem: (i) the distribution of the random operator
$\LinOp{\xi}$, as expressed via the distribution of $\xi$, (ii) the
distribution of the noise variable $w \defn y - \LinOp{\xi}
\thetastar$, and (iii) the constraints on the unknown parameter
$\thetastar$.

\subsection{Problem formulation, notation, and assumptions}

With these motivating examples in mind, we now turn to a more precise
mathematical formulation of the estimation problem introduced above.

\subsubsection{Assumptions on the random variables $(\xi, w)$}

Let us start by discussing properties of the random operator
$\LinOp{\xi}$.  In the examples previously introduced, the domain of
the observation operator $\LinOp{\xi}$ was either a subset of $\R^d$,
or more generally, a subset of the sequence space $\ell^2(\Nat)$.  The
bulk of our analysis focuses on the finite-dimensional setting
---i.e., with domain $\R^d$---so that $\LinOp{\xi}$ can be identified
with a random matrix $\R^{\numobs \times \dimension}$, for some pair
$(\numobs, \dimension)$ of positive but finite integers.  However, as
we highlight in~\Cref{sec:examples-nonparametric}, simple
approximation arguments can be used to leverage our finite-dimensional
results to determine minimax rates of convergence for estimating an
element $\thetastar$ of the infinite-dimensional sequence space
$\ell^2(\Nat)$.

In terms of the probabilistic structure of $\LinOp{\xi}$, we assume
the random element $\xi$ lies in the measurable space $(\XiSpace,
\cE)$, and is drawn from a probability measure $\bP$ on the same
space. Throughout we take $\cE$ to be large enough such that linear
functionals of $\LinOp{\xi}$ are measurable.

As for the noise vector $w \in \R^\numobs$, we assume it is
drawn---conditionally on $\xi$---from a noise distribution with
conditional mean zero, and bounded conditional covariance. Formally,
we assume that $w \sim \nu(\cdot \mid \xi)$ where $\nu$ is a Borel
regular conditional probability on $\R^\numobs$ that satisfies the
following two conditions:
\begin{enumerate}
\item [(N1)] \label{ass:centering}
For $\bP$-almost every $\xi \in \XiSpace$, we have $\int w \, \nu(\ud w\mid
\xi) = 0$; and
\item [(N2)]
\label{ass:bounded} For $\bP$-almost every $\xi \in \XiSpace$,
  we have
\begin{align*}
  \int (u^\T w)^2 \, \nu(\ud w \mid \xi) \leq u^\T \NoiseVariance u,
  \qquad \mbox{for any fixed}~u \in \R^\opdim.
\end{align*}
\end{enumerate}
We write that the measure $\nu$ lies in the set $\cP(\NoiseVariance)$
when these two conditions are satisfied.

In words, Assumption~\ref{ass:centering} requires that $w$ is
conditionally centered, and Assumption~\ref{ass:bounded} assumes that
the conditional covariance of $w$ is almost surely upper bounded in
the semidefinite ordering by $\NoiseVariance$.  Let $\bP \times \nu$
denote the distribution of the tuple $(\xi, w)$; in explicit terms,
writing $(\xi, w) \sim \bP \times \nu$ means that $\xi \sim \bP$ and
$w \mid \xi \sim \nu(\cdot \mid \xi)$.  Having specified the joint law
of $(\xi, w)$, the random variable $y$ then satisfies the stated
observation model~\eqref{eqn:operator-model}.

\subsubsection{Decision-theoretic formulation}

In this paper, our goal to estimate $\thetastar$ to the best possible
accuracy as measured by a fixed quadratic form. To
make this rigorous, we introduce two symmetric positive definite
matrices $\EstimationMat$ and $\ConstraintMat$, which induce
(respectively) the squared norms
\begin{align*}
\|\theta\|_{\EstimationMat}^2 \defn \langle \theta, \EstimationMat
\theta\rangle \quad \mbox{and} \quad
\|\theta\|_{\ConstraintMat^{-1}}^2 \defn \langle \theta,
\ConstraintMat^{\!-1} \theta \rangle,
\end{align*}
defined for any $\theta \in \R^\dimension$.  We seek estimates
$\thetahat$ of $\thetastar$ that have low squared \emph{estimation
error} $\vecnorm{\thetahat - \thetastar}_{\EstimationMat}^2$, as
defined by the matrix $\EstimationMat$. In parallel, we assume that
underlying parameter is bounded in the \emph{constraint norm}, so that
it lies in the ellipse
\begin{align*}
\Theta(\radius, \ConstraintMat) \defn \Big\{\, \theta \in
\R^\dimension : \vecnorm{\theta}_{\ConstraintMat^{-1}} \leq
\radius\,\Big\}
\end{align*}
with radius $R$, as defined by the matrix $\ConstraintMat$.

With this notation in hand, the central object of study in this paper
is the \emph{minimax risk}
\begin{align}
\label{eqn:minimax-rate}
\mathfrak{M}(T, \bP, \NoiseCovariance, \radius, \EstimationMat,
\ConstraintMat) &\defn \inf_{\thetahat} \sup_{\substack{\thetastar \in
    \Theta(\radius, \ConstraintMat) \\ \nu \in \cP(\NoiseCovariance)}}
\E_{(\xi, w) \sim \bP \times \nu} \Big[\vecnorm{\thetahat -
    \thetastar}^2_{\EstimationMat} \Big],
\end{align}
where the infimum ranges over all measurable functions $\thetahat
\equiv \thetahat(\LinOp{\xi}, y)$ that map the observed pair
$(\LinOp{\xi}, y)$ to $\R^\dimension$.

\subsection{Examples of choices of sampling laws, constraints and error norms}
\label{sec:examples-high-level}

As discussed previously, our general theory accommodates various forms
of the random linear operators $\LinOp{\xi}$.  As might one expect,
the sampling law $\bP$ for $\xi$ changes the statistical structure of
the observations, and so influences the quality of the best possible
estimates.  Moreover, the interaction between $\bP$ and the geometry
of the error norm, as defined by the matrix $\EstimationMat$, plays an
important role.  Finally, both of these factors interact with the
geometry of the constraint set, as determined by the matrix
$\ConstraintMat$.

Below we discuss some examples of these types of interactions.  To be
clear, each of these statistical settings have been considered
separately in the literature previously; one benefit of our approach
is that it provides a unifying framework that includes each of these
problems as special cases.

\begin{example}[Covariate shift in linear regression]
  \label{ex:linear-reg-multiple-pop}
Recall the set-up for linear regression, as introduced in
Example~\ref{ex:linear-regression}.  In practice, the \emph{source
distribution} from which the covariates $x$ are sampled when
constructing an estimate of $\thetastar$ need not be the same as the
\emph{target distribution} of covariates on which the predictor is to
be deployed.  This phenomenon---a discrepancy between the source and
target distributions---is known as \emph{covariate shift}.  It is now
known to arise in a wide variety of applications (e.g., see the
papers~\cite{MolEtAl20,Koh21} and references therein for more
details).

As one concrete example, in healthcare applications, the covariate
vector $x \in \R^d$ might correspond to various diagnostic measures
run on a given patient, and the response $y \in \R$ could correspond
to some outcome variable (e.g., blood pressure).  Clinicians might use
one population of patients to develop a predictive model relating the
diagnostic measures $x$ to the outcome $y$, but then be interested in
making predictions for a related but distinct population of patients.

In our setting, suppose that we use the linear model $\theta \mapsto
\yhat \defn \inprod{\theta}{x}$ to make predictions over a collection
of covariates with distribution $Q$.  A simple computation shows that
the mean-squared prediction error, averaging over both the noise $w$
and random covariates $x$, takes the form
\begin{align*}  
  \Exs\big[ (\yhat - y)^2 \big] & = \underbrace{(\theta -
    \theta^\star)^\T \Sigma_Q (\theta -\theta^\star)}_{\eqcolon \,
    L_Q(\thetahat, \thetastar)} + c, \qquad \mbox{where} \quad\Sigma_Q
  \defn \E_Q[ x\otimes x],
\end{align*}
and $c$ is a constant independent of the pair $(\theta, \thetastar)$.
Thus, the excess prediction error over the
new population $Q$ corresponds to taking $\EstimationMat = \Sigma_Q$
in our general set-up.  Similarly, if one wanted to assess parameter
error, then studying the minimax risk with the choice $\EstimationMat
= I_\dimension$ would be reasonable. Finally, the error in the
original population (denoted $P$) can be assessed with the choice
$\EstimationMat = \Sigma_P \defn \E_P[ x\otimes x]$.

Among the claims in the paper of Mourtada~\cite{Mou22} is the
following elegant result: when no constraints are imposed on
$\thetastar$, the minimax risk in the squared metric $L_Q(\thetahat,
\thetastar) = \|\thetahat - \thetastar\|_{\Sigma_Q}^2$ is equal to
\begin{align}
\label{eqn:mourtada-result}
\inf_{\thetahat} \sup_{\theta^\star \in \R^\dimension}
\E\Big[L_Q(\thetahat, \thetastar)\Big] = \frac{\sigma^2}{\numobs}
\E[\trace (\EmpCov^{-1} \Sigma_Q)],
\end{align}
where $\EmpCov$ denotes the sample covariance matrix $(1/\numobs)
\sum_{i=1}^\numobs x_i \otimes x_i$, and the expectation is over
$x_1,\dots, x_\numobs \simiid P$.  Thus, the fundamental rate of
estimation depends on the distribution of the sample covariance
matrix, the noise level, and the target distribution $Q$.

In this paper, we derive related but more general results that allow
for many other choices of the error metric and, perhaps more
importantly, permit the statistician to incorporate constraints on the
parameter $\thetastar$.  We demonstrate in
Section~\ref{sec:rederive-Mourtada} that these more general results
allow us to recover the known relation~\eqref{eqn:mourtada-result} via
a simple limiting argument where the constraint radius tends to
infinity.
\end{example}
  
\begin{example}[Nonparametric regression with non-uniform sampling]
  \label{ex:nonparametric-reg-non-uniform}
Consider observing covariate-target pairs $\{(x_i,
y_i)\}_{i=1}^\numobs$ where $y_i$ is modeled as being a noisy
realization of a conditional mean function; \ie we have $y_i =
f^\star(x_i) + w_i$ where $f^\star(x) = \E[Y \mid X = x]$, analogously
to Example~\ref{ex:nonparametric-regression}.  When $f^\star$ is
appropriately smooth and the covariates are drawn from a uniform
distribution over some compact domain, this problem has been
intensively studied, and the minimax risks are well-understood.
However, when the sampling of the covariates $x_i$ is non-uniform, the
possible rates of estimation can deteriorate drastically---see for
instance the papers~\cite{Gai05, Gai07, Gai07b, Gai09, GuiKlu11,
  AntPenSap14}.

Using tools from the theory of reproducing kernel Hilbert spaces
(RKHSs), one can formulate this problem as an infinite-dimensional
counterpart to our model~\eqref{eqn:operator-model}, where the
constraint parameters $(\radius, \ConstraintMat)$ are determined by
the Hilbert radius and the eigenvalues of the integral operator
associated with the kernel. Although formally our minimax risk is
defined for finite dimensional problems, via limiting arguments, it is
straightforward to obtain consequences for the infinite-dimensional
problem of the type discussed here, 
which discuss in Section~\ref{sec:examples-nonparametric}.
\end{example} 

\begin{example}[Covariate shift in nonparametric regression]
Combining the scenarios in Examples~\ref{ex:linear-reg-multiple-pop}
and~\ref{ex:nonparametric-reg-non-uniform}, now consider the problem
of covariate shift in a nonparametric setting.  We observe samples
$(x_i, y_i)$ where the covariates have been drawn according to some
law $P$, and our goal is to construct a predictor with low risk in
the squared norm defined by some other covariate law $Q$.

In our study of this setting, the constraint set is determined by the
underlying function class in a manner analogous to
Example~\ref{ex:nonparametric-reg-non-uniform}, and the error metric
is determined by the new distribution of covariates on which the
estimates must be deployed, analogously to
Example~\ref{ex:linear-reg-multiple-pop}.  Some recent work has
studied general conditions on the pair $(P, Q)$ and the corresponding
optimal rates of estimation~\cite{KpoMar21, Gog22, PatMaWai22, MaPatWai22,
  SchZam22, Kai23, SimEtAl23, GogEtAl23}.   Among the consequences of our work are more refined results
that are instance-dependent, in the sense that we characterize
optimality for fixed pairs $(P, Q)$, as opposed to optimality over
broad classes of $(P,Q)$ pairs. See~\Cref{sec:examples-nonparametric-shift} for
a detailed discussion of these refined results.
\end{example}

The examples above share the common feature of being problems where
estimating a conditional mean function is able to be formulated within
the observation model~\eqref{eqn:operator-model}. Additionally, in
these examples, the fundamental hardness of the problem depends on
both the structure of this function (modelled via assumptions on
$\thetastar$) as well as the distribution of the covariates.  The goal
of this paper is to build a general theory for these types of
observation models, which elucidates how both the structure of
$\thetastar$ as well as the covariate law $\bP$ determine the minimax
rate of estimation in finite samples. In Section~\ref{sec:examples},
we give concrete consequences of our general results for these types
of problems.

\subsection{Connections and relations to prior work}

Let us discuss in more detail some connections and relations between
our problem formulation and results, and various branches of the
statistics literature.

\myparagraph{Connections to random design regression} As shown by the
examples discussed so far, our general set-up includes, among other
problems, many variants of \emph{random design regression}.  This is a
classical problem in statistics, with a large literature; see the
sources~\cite{GyoEtAl02, Tsy09, HsuEtAl14} and references therein for
an overview. The recent paper~\cite{Mou22} also studies the analogous
problem studied here when the vector $\thetastar$ is allowed to be
arbitrary; the only assumption made is that $\thetastar \in \R^d$. In
this case, it is possible to use tools from Bayesian decision theory
to exhibit the minimax optimality of the ordinary least squares (OLS)
estimator~\cite[Theorem 1]{Mou22}.  In \Cref{sec:rederive-Mourtada},
we demonstrate how to obtain this result as a corollary of our more
general results. Note that in applications, such as those given by the
preceding examples,
it is important that there is
a constraint on $\thetastar$.  For instance, in a nonparametric
regression problem, the parameter $\thetastar$ denotes the
coefficients of a series expansion corresponding to a conditional mean
function $f^\star(x) = \E [Y \mid X = x]$ in an appropriate
orthonormal family of functions. In this case, one can obtain
consistent estimators of $f^\star$ only if $\thetastar$ lies in a
compact set.

\myparagraph{Random design and Bayesian priors} When the the norm of
the vector $\thetastar$ is constrained, there are relatively few
minimax results in the random design setting. On the other hand, a
related Bayesian setting has been studied.  In this line of work, the
definition of the minimax risk is altered so that the ``worst-case''
supremum over $\thetastar$ in the constraint set is replaced with a
suitable ``average''---namely the expectation over $\thetastar$ drawn
according to a prior distribution over the constraint set.

In addition to the clear differences in the formulation, this
line of work exhibits two main qualitative differences from our paper.
First, these Bayesian results have primarily been established in the
proportional asymptotics framework, in the ratio $d/n$ is assumed to
converge towards some aspect ratio $\gamma > 0$ as both $(d, n)$
diverge to infinity.  Secondly, by selecting ``nice priors'', it is
possible to leverage certain properties---for instance, equivariance to
some group action---that can hold for \emph{both} the prior and
covariate law. On the other hand, our setting is somewhat more
challenging in that we make no \emph{a priori} assumptions about the
covariate law and its relationship to the constraint set.

In more detail, when the covariates are drawn from a multivariate
Gaussian, for certain constraint sets, it is possible to find a prior
such that the minimax and Bayesian risks coincide.  As one example,
Dicker~\cite{Dic16} studies the asymptotic minimax risk when the ratio
$d/n$ is allowed to grow, and by using equivariance arguments, he
obtains asymptotically minimax procedure.  Proposition 3(b) in his
paper gives a prior for which the minimax and Bayesian risks
coincide. The thesis~\cite[Corollary 8.2]{MouThesis20} provides a
matching asymptotic lower bound.  The relation between Bayes and
minimax risks in this line of work cannot be expected in general, as
the arguments repose critically on the rotation invariance of the
standard multivariate Gaussian.  Moreover, this and other classical
work on random design regression using Gaussian covariates typically
hinges on special, closed-form formulae for quantities related to the
distribution of the sample covariance matrix (see, e.g., the
papers~\cite{Ste60, BreFre83, And03}).

\myparagraph{Fixed design results}

Although we focus on minimax estimation of the unknown parameter
$\thetastar$ in the random design setting, we note that the related
fixed design setting is well studied. In fact, in classical work,
Donoho studied a very similar operator-based observation model to
the one considered here; a key difference is that in that work, the
focus is on estimating a (scalar-valued) functional of
$\thetastar$~\cite{Don94}.

By sufficiency arguments, our problem, when instantiated in the
setting of fixed design with Gaussian noise, is equivalent to mean
estimation on an elliptical parameter set. It is therefore related to
classical work on sharp asymptotic minimax estimation in the Gaussian
sequence model~\cite{Pin80, Gol90, DonEtAl90, DonJoh94, BelLev95,
  GolTsy01, GolTsy03}; see also the monograph~\cite{Johnstone2019} for
a pedagogical overview of this topic.  These works extend the
classical line of work on estimating a constrained (possibly
multivariate) Gaussian mean~\cite{CasStr81, Bic81, MelRit87, Ber90,
  Mar93}.  We refer the reader to references~\cite{MarStr04,
  FouEtAl18}, which contain a more thorough overview of prior work on
minimax estimation of a parameter when a notion of `signal to noise
ratio' is fixed.  Of course, applying an optimal fixed design
estimator cannot be expected to yield an optimal random design
estimator in general. This is because in the fixed design formulation,
the worst-case $\thetastar$ could adapt to a single design matrix,
whereas in the random design formulation, the worst-case $\thetastar$
must adapt to the \emph{random ensemble} of design matrices induced by
sampling $\numobs$ samples in an \iid{} fashion from a fixed covariate
law.

\section{Main results}\label{sec:main-result}

We now turn to the presentation of our main results, which are upper
and lower bounds on the minimax rate of estimation as defined in
display~\eqref{eqn:minimax-rate}, matching up to a constant
pre-factor.  These bounds are presented in
Section~\ref{SecMainBounds}.

\subsection{General upper and lower bounds}
\label{SecMainBounds}

Our general upper bounds are stated as the following functional of the
distribution of the operator $\LinOp{\xi}$; the noise covariance
$\NoiseCovariance$; the constraint norm, as determined by the pair
$(\radius, \ConstraintMat)$; and the estimation norm, as defined by
the operator $\EstimationMat$,
\begin{multline}
\label{eqn:defn-functional-general}
\Phi(T, \bP, \NoiseCovariance, \radius, \EstimationMat,
\ConstraintMat) \\\defn \sup_{\Omega} \Big\{\, \E \trace\Big(
\EstimationMat^{\!\!1/2} (\Omega^{-1} + \LinOp{\xi}^\T
\NoiseCovariance^{-1} \LinOp{\xi})^{-1} \EstimationMat^{\!\!1/2} \Big)
: \Omega \succ 0,~\trace(\ConstraintMat^{\!-1/2} \Omega \ConstraintMat^{\!-1/2})
\leq \radius^2 \,\Big\}.
\end{multline}

\noindent Our first main result is a general upper bound.
\begin{theorem}[General minimax upper bound]
  \label{thm:main-upper}
The minimax risk is upper bounded as
\begin{equation}
\label{eqn:gen-minimax-upper-bound}
\mathfrak{M}(T, \bP, \NoiseCovariance, \radius, \EstimationMat,
\ConstraintMat) \leq \Phi(T, \bP, \NoiseCovariance, \radius,
\EstimationMat, \ConstraintMat).
\end{equation}
\end{theorem}
\noindent See~\Cref{sec:proof-upper} for the proof. \\

\medskip

\noindent Our second result is a complementary lower bound.
\begin{theorem}[Lower bound]
  \label{thm:main-lower}
The minimax risk is lower bounded as
\begin{equation}
\label{eqn:gen-minimax-lower-bound}
\mathfrak{M}(T, \bP, \NoiseCovariance, \radius, \EstimationMat,
\ConstraintMat) \geq \, \Phi(T, \bP, \NoiseCovariance,
\tfrac{\radius}{2}, \EstimationMat, \ConstraintMat) \geq \frac{1}{4}
\, \Phi(T, \bP, \NoiseCovariance, \radius, \EstimationMat,
\ConstraintMat).
\end{equation}
\end{theorem}
\noindent See~\Cref{sec:proof-lower} for the proof. \\

\medskip

Note that the functional on the righthand side of the
display~\eqref{eqn:gen-minimax-lower-bound} above matches the quantity
appearing in our minimax upper bound~\eqref{eqn:gen-minimax-upper-bound}.
Thus, in a nonasymptotic fashion, we have determined the minimax
risk for this problem up the prefactor $1/4$.

\paragraph{Sharper lower bound constants}
The constant appearing in the lower bound~\eqref{eqn:gen-minimax-lower-bound}
can typically be substantially sharpened.
To describe how this can be done via our results,
fix a scalar $\tau \in (0, 1]$ and a symmetric positive
definite matrix $\Omega$, and let $Z \in \R^\dimension$ be vector of
\iid{} standard Gaussians. Define the scalar
\begin{equation*}
c \defn \tau^2 (1 - \P\{ \tau^2 \sum_{i=1}^\dimension \lambda_i Z^2_i
> 1\}),
\end{equation*}
where $\{\lambda_i\}_{i=1}^\dimension$ are the the eigenvalues of the matrix
$(1/\radius^2)\EstimationMat^{1/2} \Omega \EstimationMat^{1/2}$.
Then, we are able to establish the following minimax lower bound,
\begin{equation}
\label{eqn:gen-sharp-lower}
\mathfrak{M}(T, \bP, \NoiseCovariance, \radius, \EstimationMat,
\ConstraintMat) \geq \E \trace\Big( \EstimationMat^{\!\!1/2}
(\frac{1}{c} \Omega^{-1} + \LinOp{\xi}^\T \NoiseCovariance^{-1}
\LinOp{\xi})^{-1} \EstimationMat^{\!\!1/2} \Big),
\end{equation}
provided that the parameter $\tau \in (0,1]$ and the symmetric positive definite
matrix $\Omega$ is
  such that $\trace(\ConstraintMat^{\!-1/2} \Omega \ConstraintMat^{\!-1/2})= \radius^2$.

With appropriate choices of the pair $(\tau, \Omega)$, the lower
bound~\eqref{eqn:gen-sharp-lower} can lead to pre-factors that are
much closer to $1$, and in some cases, converge to one under various
scalings. In Section~\ref{sec:gaussian-linear-regression}, we give one
illustration of how the family of bounds~\eqref{eqn:gen-sharp-lower}
can be exploited to obtain an improvement of this type.

\paragraph{Form of an optimal procedure}
Inspecting the proof of~\Cref{thm:main-upper}---specifically, as a
consequence of Proposition~\ref{prop:minimizer-of-linear-risk}---if
the supremum on the righthand side of
\eqref{eqn:defn-functional-general}
is attained at the matrix
$\Omega_{\star}$, then the following estimator, in view of the lower
bound~\eqref{eqn:gen-minimax-lower-bound}, is near minimax-optimal,
\begin{equation}
\label{eqn:opt-estimator}
\thetahat(\LinOp{\xi}, y) \defn \big(\Omega_\star^{-1} +
\LinOp{\xi}^\T \NoiseCovariance^{-1} \LinOp{\xi}\big)^{-1}
\LinOp{\xi}^\T \NoiseCovariance^{-1} y.
\end{equation}
It is perhaps instructive to write this estimator in its ``ridge''
formulation
\begin{align*}
\thetahat(\LinOp{\xi}, y) = \argmin_{\vartheta \in \R^\dimension}
\Big\{\, \|y - \LinOp{\xi} \vartheta\|_{\NoiseCovariance^{-1}}^2 +
\|\vartheta\|_{\Omega_\star^{-1}}^2 \,\Big\}.
\end{align*}
In the language of Bayesian statistics, our order-optimal procedure is
a maximum \emph{a posteriori} (MAP) estimate for $\thetastar$ when $y
\sim \Normal{\LinOp{\xi} \thetastar}{\NoiseCovariance}$ and the
parameter follows the prior distribution $\thetastar \sim
\Normal{0}{\Omega_\star}$.  The optimal prior is identified via the
choice of $\Omega_\star$ which is determined by the functional
appearing in Theorems~\ref{thm:main-upper} and~\ref{thm:main-lower}.
If the supremum in~\eqref{eqn:defn-functional-general}
is not attained, then by selecting a sequence of
matrices $\Omega_k$ that approach the maximal value of the functional,
one can similarly argue there exists a sequence of estimators that
approach the order-optimal minimax risk.

\subsection{Independent and identically distributed regression models}
\label{sec:IID-regression-models}

An important application of our general result is for independent and
identically distributed (\iid) regression models of the form
\begin{equation}
\label{eqn:iid-regression}
y_i = \langle \thetastar, \psi(x_i) \rangle + \sigma z_i, \quad
\mbox{for}~i = 1, \ldots, \numobs.
\end{equation}
Above, we assume that $x_i$ are independent and identical draws from a
fixed covariate distribution $P$, on some measurable space $\cX$, and
that $\psi \colon \cX \to \R^\dimension$. The covariates $\{x_i
\}_{i=1}^\numobs$ are independent and the conditional distribution of
$z \mid x$ is an element of $\cP(I_\numobs)$. The parameter $\sigma >
0$ indicates the noise level; it is an upper bound on the conditional
standard deviation of $y_i - \langle \thetastar, \psi(x_i)\rangle$.

For the model described above, the following minimax risk of
estimation provides the best achievable performance of any estimator,
when $\thetastar$ lies in a compact ellipse and the error is measured
in the quadratic norm
\begin{equation}
\label{defn:minimax-rate-iid}
\mathfrak{M}^{\rm \iid}_\numobs\Big(\psi, P, \radius, \sigma^2,
\ConstraintMat, \EstimationMat\Big) \defn \inf_{\thetahat}
\sup_{\substack{\thetastar \in \Theta(\radius, \ConstraintMat) \\ \nu
    \in \cP(\sigma^2 I_\numobs)}}
\E\Big[\vecnorm[\big]{\thetahat(y_1^\numobs, x_1^\numobs) -
    \thetastar}_{\EstimationMat}^2\Big].
\end{equation}
Note that this problem can be formulated as an instance of our general
operator formulation~\eqref{eqn:operator-model} where we take $y =
(y_1,\dots, y_\numobs)$, $w = \sigma (z_1, \dots, z_\numobs)$, and
$\xi=(x_1, \dots, x_\numobs)$, so that $\bP = P^\numobs$. The operator
$\LinOp{\xi}$ is given by the $\numobs \times \dimension$-matrix with
rows $\psi(x_i)^\T$.  In this context the following random matrix,
which is a rescaling of the operator $\LinOp{\xi}^\T \LinOp{\xi}$,
plays an important role:
\begin{equation}
\label{eqn:empirical-covariance}
\EmpCov \defn \frac{1}{\numobs}\sum_{i=1}^\numobs \psi(x_i) \otimes
\psi(x_i).
\end{equation}

In order to state the consequence of our more general results for this
problem, let us introduce a functional. We denote it by $d_\numobs$ to
indicate that it is essentially an ``effective statistical dimension''
for this problem,
\begin{equation}\label{eqn:d-functional-iid}
\resizebox{0.91\hsize}{!}{ $d_\numobs(\psi, P, \radius, \sigma^2,
  \EstimationMat, \ConstraintMat) \defn \sup_{\Omega} \! \Big\{\trace
  \E_{P^\numobs}\!\big[ \EstimationMat^{\!\!1/2} (\EmpCov +
    \Omega^{-1} )^{-1}\EstimationMat^{\!\!1/2}\big] : \Omega \succ 0,
  \trace(\ConstraintMat^{\!\!-1/2} \Omega \ConstraintMat^{\!\!-1/2})
  \leq \frac{\numobs \radius^2}{\sigma^2}\Big\}$ }.
\end{equation}
Then an immediate corollary to Theorems~\ref{thm:main-upper}
and~\ref{thm:main-lower} is the following pair of inequalities for the
\iid{} minimax risk.\footnote{ Strictly speaking, this result follows immediately if we
had defined the minimax risk over estimators which are measurable
functions of the variables $\{(y_i, \psi(x_i))\}$.  Nonetheless, since
our lower bounds use Gaussian noise, the stated inequalities hold even
when defining the minimax risk for estimators which operate on
$\{(y_i, x_i)\}$, by a standard sufficiency argument.  }

\bcor
\label{cor:iid-result}
Under the IID regression model~\eqref{eqn:iid-regression}, the minimax
rate of estimation as defined in
equation~\eqref{defn:minimax-rate-iid} satisfies the following
inequalities,
\begin{multline}\label{eqn:rate-for-iid-models}
\frac{1}{4} \, \frac{\sigma^2}{\numobs} d_\numobs(\psi, P, \radius,
\sigma^2, \EstimationMat,  \ConstraintMat) \leq
\frac{\sigma^2}{\numobs} d_\numobs(\psi, P, \tfrac{\radius}{2},
\sigma^2, \EstimationMat, \ConstraintMat) \\ \leq \mathfrak{M}^{\rm
  \iid}_\numobs\Big(\psi, P, \radius, \sigma^2,
\EstimationMat, \ConstraintMat\Big)
\leq \frac{\sigma^2}{\numobs} d_\numobs(\psi, P,
\radius, \sigma^2, \EstimationMat, \ConstraintMat).
\end{multline}
\ecor
So as to lighten notation, in the sequel, 
when the feature map $\psi$ is the identity mapping $\psi(x) = x$, 
we drop the parameter $\psi$ from the functional $d_\numobs$ and the minimax 
rate $\mathfrak{M}^{\rm \iid}_\numobs$.

\subsection{Some properties of the functional appearing in Theorems~\ref{thm:main-upper} and~\ref{thm:main-lower}}

As indicated by Theorem~\ref{thm:main-upper} and the subsequent
discussion, the extremal quantity
\begin{equation}
\label{eqn:general-functional}
\sup_{\Omega} \Big\{\, \E \trace\Big( \EstimationMat^{\!\!1/2}
(\Omega^{-1} + \LinOp{\xi}^\T \NoiseCovariance^{-1} \LinOp{\xi})^{-1}
\EstimationMat^{\!\!1/2} \Big) : \Omega \succ
0,~\trace(\ConstraintMat^{\!-1/2} \Omega \ConstraintMat^{\!-1/2})\leq
\radius^2 \,\Big\}
\end{equation}
is fundamental in that it determines our minimax risk; moreover when
the supremum is attained, the maximizer defines an order-optimal
estimation procedure (see equation~\eqref{eqn:opt-estimator}).
Conveniently, it turns out that the maximization problem implied by
the display~\eqref{eqn:general-functional} is concave.

\begin{proposition}[Concavity of functional]
  \label{prop:concavity-functional}
The optimization problem
\begin{equation}
\label{eqn:gen-functional-program}
\begin{aligned}
& \text{maximize} & & f(\Omega) \defn \trace
  \E\big[\EstimationMat^{\!1/2} \big(\Omega^{-1} + \LinOp{\xi}^\T
    \NoiseCovariance^{-1} \LinOp{\xi} \big)^{-1}
    \EstimationMat^{\!1/2}\big] \\ & \text{subject to} & & \Omega
  \succ 0, \quad \trace(\ConstraintMat^{\!-1/2} \Omega
  \ConstraintMat^{\!-1/2})\leq \radius^2,
\end{aligned}
\end{equation}
is equivalent to a convex program, with variable $\Omega$. Formally,
the constraint set above is convex, and function $f$ is concave over this set.
\end{proposition}
\noindent See Appendix~\ref{app:proof-concavity-functional} for the
proof.\\

\medskip

Note that this claim implies that, provided oracle access to the
objective function $f$ appearing above, one can in principle obtain a
maximizer in a computationally tractable manner, by leveraging
algorithms for convex optimization~\cite{BV04}.

The functional~\eqref{eqn:general-functional} depends on the
distribution of $\LinOp{\xi}^\T \NoiseCovariance^{-1} \LinOp{\xi}$.
In general, Jensen's inequality along with the convexity of the trace
of the inverse of positive matrices~\cite[Exercise 1.5.1]{Bha07}
implies that it is always lower bounded by
\begin{equation}
\label{eqn:general-functional-lower}
\sup_{\Omega} \Big\{\, \trace\Big( \EstimationMat^{\!\!1/2}
(\Omega^{-1} + \E \LinOp{\xi}^\T \NoiseCovariance^{-1}
\LinOp{\xi})^{-1} \EstimationMat^{\!\!1/2} \Big) : \Omega \succ
0,~\trace(\ConstraintMat^{\!-1/2} \Omega \ConstraintMat^{\!-1/2})\leq
\radius^2 \,\Big\}
\end{equation}
Comparing displays~\eqref{eqn:general-functional}
and~\eqref{eqn:general-functional-lower}, we have simply moved the
expectation over $\xi$ into the inverse. For certain \IID{} regression
models, as described in Section~\ref{sec:IID-regression-models}, we
can give a complementary upper bound. To state our result, we define
\begin{equation}\label{eqn:d-functional-iid-pop}
  \resizebox{0.91\hsize}{!}{%
$\overline{d}_\numobs(P, \radius, \sigma^2, 
    \EstimationMat, \ConstraintMat) \defn \sup_{\Omega} \! \Big\{\trace
    \big(\EstimationMat^{\!\!1/2} (\E_{P^\numobs} \EmpCov +
    \Omega^{-1} )^{-1}\EstimationMat^{\!\!1/2}\big) : \Omega \succ 0,
    \trace(\ConstraintMat^{\!\!-1/2} \Omega \ConstraintMat^{\!\!-1/2})
    \leq \frac{\numobs \radius^2}{\sigma^2}\Big\}$ }.
\end{equation}
Note that this quantity only depends on the distribution $P^\numobs$
through the matrix $\E_{P^\numobs} \EmpCov$.

\begin{proposition}[Comparison of $d_\numobs$ to $\overline{d}_\numobs$]
\label{prop:to-population}
Suppose that $\Sigma_P \defn \E_P[\psi(x) \otimes \psi(x)]$ is
nonsingular.  Define $\kappa$ to be the $P$-essential supremum of $x
\mapsto \|\ConstraintMat^{1/2} \psi(x)\|_2$. If $\kappa < \infty$,
then for any $\radius > 0, \sigma > 0$, we have
\begin{align*}
\overline{d}_\numobs(\psi, P, \radius, \sigma^2, 
\Sigma_P, \ConstraintMat)
\leq d_\numobs(\psi, P, \radius, \sigma^2, \Sigma_P, \ConstraintMat
) \leq \Big(1 + \frac{\radius^2 \kappa^2}{\sigma^2}\Big)
\overline{d}_\numobs(\psi, P, \radius, \sigma^2, 
\Sigma_P, \ConstraintMat).
\end{align*}
\end{proposition}
\noindent Unpacking this result, when $K_c^{1/2} \psi(x)$ is
essentially bounded, for problems where the error is measured in the
norm induced by the covariance $\Sigma_P$, we see that the functionals
$\overline{d}_\numobs$ and $d_\numobs$ are of the same order when the
signal-to-noise ratio satisfies the relation
$\tfrac{\radius^2}{\sigma^2} \lesssim \tfrac{1}{\kappa^2}$.  As
mentioned in the discussion above, the first inequality above is a
consequence of a generic lower bound.  See
Appendix~\ref{app:proof-to-population} for the proof of the upper
bound in the claim.

\subsection{Asymptotics for a diverging radius}

In this section, we develop an asymptotic limit relation for the
minimax risk~\eqref{eqn:minimax-rate} as the radius $\radius$ of the
constraint set $\Theta(\radius, \ConstraintMat)$ tends to
infinity. The relation reveals that the lower bound constant $1/4$
appearing in the lower bound Theorem~\ref{thm:main-lower} can actually
be made quite close to $1$ for large radii.

\bcor \label{cor:limit-relation}
Suppose that $\LinOp{\xi}^\T \NoiseCovariance^{-1} \LinOp{\xi}$ is
$\bP$-almost surely nonsingular.  Then the minimax
risk~\eqref{eqn:minimax-rate} satisfies
\begin{align*} 
\mathfrak{M}(T, \bP, \NoiseCovariance, \radius, \EstimationMat,
\ConstraintMat) = \big(1 - o(1)\big) \, \Phi(T, \bP, \NoiseCovariance,
\radius, \EstimationMat, \ConstraintMat) , \quad \mbox{as}~\radius \to
\infty.
\end{align*}
\ecor
\noindent See Appendix~\ref{app:proof-limit-relation} for a proof of
this claim.

An immediate consequence is that for \iid{} regression settings as in
Section~\ref{sec:IID-regression-models}, we have the following limit
relation.

\bcor\label{cor:limit-relation-iid} Suppose that that the empirical
covariance matrix $\EmpCov$ from
equation~\eqref{eqn:empirical-covariance} is $P^\numobs$-almost surely
invertible.  Then, the minimax risk for an \iid{} observation
model~\eqref{eqn:iid-regression} satisfies the relation
\begin{align*}
\mathfrak{M}^{\rm \iid}_\numobs\Big(\psi, P, \radius, \sigma^2,
\EstimationMat, \ConstraintMat\Big) = \big(1 - o(1)\big) \,
\frac{\sigma^2}{\numobs} d_\numobs\big(\psi, P, \radius, \sigma^2,
\EstimationMat, \ConstraintMat\big), \quad \mbox{as}~\radius \to
\infty.
\end{align*}
\ecor

\section{Consequences of main results}
\label{sec:examples}

In this section, we demonstrate consequences of our main results for a
variety of estimation problems. In
Section~\ref{sec:examples-parametric}, we develop consequences of our
main results for problems where the underlying parameter to be
estimated is finite-dimensional. In
Section~\ref{sec:examples-nonparametric}, we develop consequences of
our main results for problems where the underlying parameter is
infinite-dimensional. In both cases, we are able to derive minimax
rates of estimation, which to the best of our knowledge, are not yet
in the literature.  Additionally, we are also able to re-derive
classical as well as recent results in a unified fashion via our main
theorems.

\subsection{Applications to parametric models}\label{sec:examples-parametric}

We begin by developing the consequences of our main results for regression problems
where the statistician is aiming to estimate a finite-dimensional parameter. 
Sections~\ref{sec:gaussian-linear-regression},~\ref{sec:underdetermined-regression}, 
and~\ref{sec:rederive-Mourtada} concern \iid{} regression settings of the form 
described in Section~\ref{sec:IID-regression-models}. In Section~\ref{sec:markov-regression}, 
we consider a non-\iid{} regression setting. 

\subsubsection{Linear regression with Gaussian covariates}\label{sec:gaussian-linear-regression}

As in the prior work~\cite{Dic16}, consider a random design \iid{}
regression setting of the form presented in the
display~\eqref{eqn:iid-regression}, but with Gaussian data.  Formally,
we assume Gaussian noise, so that $z_i \simiid \Normal{0}{1}$, and
Gaussian covariates, so that $x_i \simiid \Normal{0}{I_\dimension}$
and $\psi(x) = x$.  Here $x$ and $z$ are assumed independent.  Then we
define
\begin{align*}
r(\numobs, \dimension, \radius, \sigma) \defn \inf_{\thetahat}
\sup_{\twonorm{\theta} \leq \radius} \E \Big[\twonorm{\thetahat -
    \theta}^2\Big], \quad \mbox{and} \quad d_{\rm Dicker}(\numobs,
\dimension, \radius, \sigma)\defn \trace \E\Big[ (\EmpCov +
  \tfrac{\sigma^2}{\numobs} \tfrac{\dimension}{\radius^2}
  I_\dimension)^{-1} \Big],
\end{align*}
where the expectations are over the Gaussian covariates and noise
pairs $\{(x_i, z_i)\}_{i=1}^\numobs$.  These quantities correspond,
respectively, to the minimax risk and the worst-case risk (rescaled by
$\numobs/\sigma^2$), of a certain ridge estimator~\cite[Corollary
  1]{Dic16} on the sphere $\{\|\theta\|_2 = \radius\}$.

Dicker~\cite[Corollary 3]{Dic16} proves the following limiting
result. Under the proportional asymptotics $d/\numobs \to \gamma$,
where the limiting ratio $\gamma$ lies in $(0, \infty)$, the minimax risk satisfies
\begin{align}
\label{eqn:dicker}
\lim_{\dimension/\numobs \to \gamma} \Big|r(\numobs, \dimension,
\radius, \sigma) - \frac{\sigma^2}{\numobs} d_{\rm Dicker}(\numobs,
\dimension, \radius, \sigma) \Big| = 0,
\end{align}
for any radius $\radius > 0$ and noise level $\sigma > 0$.

Let us now demonstrate that our general theory yields a nonasymptotic
counterpart of this claim, and taking limits recovers the asymptotic
relation~\eqref{eqn:dicker}.

\bcor
\label{cor:our-version-dicker} For linear regression over the
$\radius$-radius Euclidean sphere with Gaussian covariates, the
minimax risk satisfies the sandwich relation
\begin{subequations}
\begin{align}
\label{eqn:our-dicker-result}
c_\dimension \, \frac{\sigma^2}{\numobs} d_{\rm Dicker}(\numobs,
\dimension, \radius, \sigma) \leq \frac{\sigma^2}{\numobs} d_{\rm
  Dicker}(\numobs, \dimension, \sqrt{c_d} \radius, \sigma) \leq
r(\numobs, \dimension, \radius, \sigma) \leq \frac{\sigma^2}{\numobs}
d_{\rm Dicker}(\numobs, \dimension, \radius, \sigma),
\end{align}
where
\begin{align}
\label{eqn:sharp-lower-bound-constant-dicker}
c_d \defn \begin{cases} (1 - \tfrac{1}{2d - 1}) (1 -
  \exp(-\tfrac{d^{3/2}}{4})) & d \geq 2 \\ 1/4 & d = 1 \end{cases}.
\end{align}
\end{subequations}
\ecor 
\noindent Note that since $c_d = (1 - o(1))$ as $d \to \infty$, the
inequalities~\eqref{eqn:our-dicker-result} allow us to immediately
recover Dicker's result. It should be emphasized, however, that
Corollary~\ref{cor:our-version-dicker}, holds for
\emph{any} quadruple $(\numobs, \dimension, \radius,
\sigma)$. In particular, it is valid in a completely nonasymptotic
fashion and with explicit constants.

We now sketch how this result follows from our main results.  As
calculated in Appendix~\ref{app:proof-gauss-calc}, our functional for
this problem satisfies
\begin{subequations}
\begin{equation}
  \label{eqn:gaussian-regression-functional}
d_\numobs(\Normal{0}{I_\dimension}, \radius, \sigma^2, I_\dimension,
I_\dimension) = d_{\rm Dicker}(\numobs, \dimension, \radius, \sigma).
\end{equation}
Hence, our Corollary~\ref{cor:iid-result} implies the following
characterization of the minimax risk,\footnote{Although
Corollary~\ref{cor:iid-result} takes the supremum over a larger family
of noise distributions, note that our lower bounds are obtained with
Gaussian noise, so that the result applies even if we restrict to
Gaussian noise.}
\begin{align}
\label{eqn:loose-dicker-lower}
\frac{1}{4} \, \frac{\sigma^2}{\numobs} d_{\rm Dicker}(\numobs,
\dimension, \radius, \sigma) \leq r(\numobs, \dimension, \radius,
\sigma) \leq \frac{\sigma^2}{\numobs} d_{\rm Dicker}(\numobs,
\dimension, \radius, \sigma^2).
\end{align}
\end{subequations}
To establish our sharper result~\eqref{eqn:our-dicker-result}, we
leverage the stronger lower bound~\eqref{eqn:gen-sharp-lower}.  The
details of this calculation are presented in
Appendix~\ref{app:lower-bound-tail-bound-calc}.  Note that in
Section~\ref{sec:random-design-simulation}, we simulate this problem
and find that as suggested by Corollary~\ref{cor:our-version-dicker},
that, indeed, the gap between our upper and lower bounds is tiny, even
for problems with small dimension (see
Figure~\ref{fig:random-reg-sim}).

\subsubsection{Underdetermined linear regression}\label{sec:underdetermined-regression}

Consider observing samples from a standard linear regression model; that is, we observe 
pairs $\{(x_i, y_i)\}$ according to the model~\eqref{eqn:iid-regression}, with 
$\psi(x) = x$. 
A practical scenario in which some assumption regarding the norm of the underlying parameter 
is necessary is when the sample covariance matrix $\EmpCov$, defined in 
display~\eqref{eqn:empirical-covariance} is singular with positive $P^\numobs$-probability. 
This occurs if $\numobs < d$, or if there is a hyperplane $H \subset \R^\dimension$ such that 
$x \sim P$ lies in $H$ with positive probability. 

In this setting, the correct dependence of the minimax risk on the geometry 
of the constraint set and the distribution of sample covariance matrix 
is relatively poorly understood. For simplicity---although our results are more general 
than this---let us assume that error is measured in the Euclidean norm and that 
it is assumed that the underlying parameter $\thetastar$ has Euclidean norm bounded by 
$\radius > 0$, and that the noise is independent Gaussian with variance $\sigma^2$. 
Then Corollary~\ref{cor:iid-result} demonstrates that 
\begin{multline*}
\inf_{\thetahat} \sup_{\|\theta\|_2 \leq \radius}  
\E[\|\thetahat - \theta\|_2^2] \asymp 
\frac{\sigma^2}{\numobs} d_\numobs( P, \radius, \sigma^2, I_\dimension, I_\dimension) 
=
\frac{\sigma^2}{\numobs} 
\sup_{\Omega \succ 0 } \! \Big\{\trace \E_{P^\numobs}\!\big[ (\EmpCov + \Omega^{-1} )^{-1} 
\big]  : \trace(\Omega) \leq \frac{\numobs \radius^2}{\sigma^2}\Big\}.
\end{multline*} 
Taking $\Omega = \tfrac{\numobs}{\dimension} \tfrac{\radius^2}{\sigma^2} I_\dimension$, 
we obtain the following lower bound on the minimax risk for any covariate law $P$, 
\begin{equation}\label{eqn:lower-bound-singular}
\resizebox{0.93\hsize}{!}{$\frac{\sigma^2}{\numobs} \trace
  \E_{P^\numobs}\!\big[ (\EmpCov + \tfrac{\sigma^2}{\radius^2}
    \tfrac{\dimension}{\numobs} I_\dimension )^{-1}\big] \asymp
  \underbrace{\E \Big[ \sum_{i=1}^{\dimension}
      \tfrac{\sigma^2}{\numobs} \tfrac{1}{\lambda_i(\EmpCov)}
      \1\{\lambda_i(\EmpCov) \geq \tfrac{\sigma^2}{\numobs} \tfrac
                 {\dimension}{\radius^2}\}
                 \Big]}_{\substack{\text{Estimation error from}
      \\ \text{large eigenvalues of $\EmpCov$}}} + \underbrace{\E
    \Big[ \sum_{i=1}^{\dimension} \tfrac{\radius^2}{\dimension}
      \1\{\lambda_i(\EmpCov) < \tfrac{\sigma^2}{\numobs} \tfrac
                 {\dimension}{\radius^2}\}
                 \Big]}_{\substack{\text{Approximation error due}
      \\ \text{to small eigenvalues of $\EmpCov$}}}.  $}
\end{equation}
The lower bound~\eqref{eqn:lower-bound-singular} is sharp in certain cases. 
For instance, when $x_i \simiid \Normal{0}{I_\dimension}$ but there are fewer samples than the dimension, so that 
$\numobs < d$, it is equal to the minimax risk up to universal constants, following the same 
argument as in Section~\ref{sec:gaussian-linear-regression}. 

Note that above, $\lambda_i$ denotes the $i$th largest (nonnegative)
eigenvalue of a symmetric positive semidefinite matrix.  One possible
interpretation of this lower bound is as follows: the first term
indicates the estimation error incurred in directions where the
effective signal-to-noise ratio is high; on the other hand, the second
term indicates the bias or approximation error that must be incurred
in directions where the effective signal-to-noise ratio is low. In
fact, the message of this lower bound is that in these directions, no
procedure can do much better than estimating $0$ there. One concrete
and interesting takeaway is that if $\EmpCov$ has an eigenvalue equal
to zero, it increases the minimax risk by essentially the same amount
as if the eigenvalue were positive and in the interval $(0,
\tfrac{\sigma^2}{\numobs}\tfrac{\dimension}{\radius^2})$.

\subsubsection{Linear regression with an unrestricted parameter space}\label{sec:rederive-Mourtada}

In recent work, Mourtada~\cite{Mou22} characterizes the minimax risk
for random design linear regression problem for an \emph{unrestricted}
parameter space. Consider observing samples $\{(x_i,
y_i)\}_{i=1}^\numobs$ following the \iid{} 
model~\eqref{eqn:iid-regression} with $\psi(x) = x$, where the
covariates are drawn from some distribution $P$ on $\R^\dimension$.
As argued by Mourtada (see his Proposition 1), or as can be seen by
taking $\radius \to \infty$ in our singular lower
bound~\eqref{eqn:lower-bound-singular} from
Section~\ref{sec:underdetermined-regression}, if we impose no
constraint on the underlying parameter $\thetastar$, then it is
necessary to assume that the sample covariance matrix $\EmpCov$ is
invertible with probability $1$ in order to obtain finite minimax
risks. Theorem 1 in Mourtada's paper then asserts that under this
condition, we have
\begin{equation}
\label{eqn:mourtada-result}
\inf_{\thetahat} \sup_{\substack{\thetastar \in \R^\dimension \\ \nu
    \in \cP(\sigma^2 I_\numobs)}} \E \Big[\vecnorm[\big]{\thetahat -
    \thetastar}^2_{\Sigma_P}\Big] = \frac{\sigma^2}{\numobs} \E
\big[\trace(\EmpCov^{-1} \Sigma_P)\big],
\end{equation}
where the expectation is over the data $\{(x_i, y_i)\}_{i=1}^\numobs$,
and $\Sigma_P \defn \E_P[x \otimes x]$ is the population covariance
matrix under $P$.

We now show that this result, with the exact constants, is a
consequence of our more general results. We focus on establishing the
lower bound, because it is well-known (and easy to show) that the
upper bound is achieved by the ordinary least squares
estimator.\footnote{Alternatively, note that if we define
$\thetahat_\radius$ to be the order-optimal estimator we derive for
the constraint set $\{\|\thetastar\|_2^2 \leq \radius^2\}$ (see
equation~\eqref{eqn:opt-estimator}, with $\ConstraintMat =
I_\dimension$, $\NoiseCovariance = \sigma^2 I_\dimension$, and
$\LinOp{\xi} = X$, where $X$ is the design matrix.), then it
converges compactly to the ordinary least squares estimate as $\radius
\to \infty$.} Thus for the lower bound, our results imply that
\begin{subequations}
\begin{align}
\label{eqn:functional-radius-mourtada}   
\inf_{\thetahat} \sup_{\substack{\thetastar \in \R^\dimension \\ \nu
    \in \cP(\sigma^2 I_\numobs)}} \E\Big[\vecnorm[\big]{\thetahat -
    \thetastar}^2_{\Sigma_P}\Big] &\geq \sup_{\radius > 0} \,
\bigg\{\inf_{\thetahat} \sup_{\substack{\|\thetastar\|_2 \leq \radius
    \\ \nu \in \cP(\sigma^2 I_\numobs)}}
\E\Big[\vecnorm[\big]{\thetahat - \thetastar}^2_{\Sigma_P} \Big]
\bigg\} \\
\label{eqn:limit-lb-mourtada}
& = \frac{\sigma^2}{\numobs} \lim_{\radius \to \infty} d_\numobs(P,
\radius, \sigma^2,  \Sigma_P, I_\dimension).
\end{align}
\end{subequations}
In order to obtain the relation~\eqref{eqn:limit-lb-mourtada}, we have
used the fact that the constrained minimax risk over the set
$\{\|\thetastar\|_2 \leq \radius\}$ is nondecreasing in $\radius > 0$,
and have applied our limit relation in
Corollary~\ref{cor:limit-relation-iid}.  A short calculation, which we
defer to Appendix~\ref{app:mourtada-limit-proof}, demonstrates that
\begin{equation}
\label{eqn:mourtada-limit-relation}
\lim_{\radius \to \infty} d_\numobs(P, \radius, \sigma^2,
\Sigma_P, I_\dimension) = \E \big[\trace(\EmpCov^{-1} \Sigma_P)\big].
\end{equation}
Thus, after combining displays~\eqref{eqn:limit-lb-mourtada}
and~\eqref{eqn:mourtada-limit-relation}, we have obtained the lower
bound in Mourtada's result~\eqref{eqn:mourtada-result}. One
consequence of this argument is that the
inequality~\eqref{eqn:functional-radius-mourtada} is, as may be
expected, an equality. That is, we have
\begin{align*}
\inf_{\thetahat} \sup_{\substack{\thetastar \in \R^\dimension \\ \nu
    \in \cP(\sigma^2 I_\numobs)}} \E\Big[\vecnorm[\big]{\thetahat -
    \thetastar}^2_{\Sigma_P}\Big] = \sup_{\radius > 0} \,
\bigg\{\inf_{\thetahat} \sup_{\substack{\|\thetastar\|_2 \leq \radius
    \\ \nu \in \cP(\sigma^2 I_\numobs)}}
\E\Big[\vecnorm[\big]{\thetahat - \thetastar}^2_{\Sigma_P} \Big]
\bigg\}.
\end{align*}
Note that establishing this equality directly is somewhat cumbersome,
as it requires essentially applying a form of a min-max theorem, which
in turn requires compactness and continuity arguments.

\subsubsection{Regression with Markovian covariates}
\label{sec:markov-regression}

We consider a dataset $\{(x_t,y_t)\}_{t=1}^T$ comprising of
covariate-response pairs. The covariates are initialized with $x_0 =
0$, and then proceed via the recursion
\begin{equation}
\label{eqn:markov-covariates}
x_{t} = \sqrt{r_t} \; x_{t-1} + \sqrt{1 - r_t} \; z_t \quad
\mbox{for}~t = 1, \ldots, T,
\end{equation}
for some collection of parameters $\{r_t\}_{t=1}^T \subset [0, 1]$, and
family of independent standard Gaussian variates $\{z_t\}_{t=1}^T$.
By construction, the samples $\{x_t\}_{t=1}^T$ form
  a Markov chain---a time-varying $\mathrm{AR}(1)$ process with
  stationary distribution being the standard Gaussian law.
  At the extreme $r_t \equiv 0$, the sequence $\{x_i\}_{i=1}^\numobs$ is
  \iid{}, whereas for $r_t \in (0,1)$, is a dependent sequence,
  and its mixing becomes slower as the parameters $\{r_t\}$ get closer to $1$.
In addition to these random covariates, suppose that we also observe responses $\{y_t\}_{t=1}^T$ from the model
\begin{equation}
    \label{eqn:markov-responses}
y_t = x_t  \theta^\star + \sigma w_t,
\qquad\mbox{for}~t=1,\ldots, T,
\end{equation}
where $\sigma > 0$ is a noise standard deviation, and the noise
sequence $\{w_t\}_{t=1}^T$ consists of \iid{} standard Gaussian
variates.  We assume that $z_t$ and $x_t$ are independent for all $t = 1, \ldots, T$.

We now describe how our main results apply to this setting.  Let us
define a matrix $M \in \R^{T \times T}$ which is associated to the dynamical system~\eqref{eqn:markov-covariates}.
It has entries
\begin{equation}\label{eqn:gen-markov-matrix}
M_{s s'}  = \sum_{t = \twomax{s}{s'}}^T \sqrt{c_{st} c_{s't}},
\quad \mbox{where} \quad c_{st} \defn (1 - r_s) \prod_{\tau = s+1}^t r_\tau. 
\end{equation}
To give one example, in the special case that $r_t \equiv \alpha \in (0, 1)$ for all $t$, then
the matrix $M$ is similar under permutation to the matrix with entries
\begin{equation*}
M_{st} = \sqrt{\alpha}^{|s-t|} -
\sqrt{\alpha}^{s + t}.
\end{equation*}
Evidently, this matrix is a rank-one update to the covariance matrix for the underlying $\mathrm{AR}(1)$ process
(\ie the Kac–Murdock–Szegö matrix~\cite{KMS53}); it is easily checked to be symmetric positive definite. 

We now state the consequences of our main results for this problem.
\bcor \label{cor:markov}
The minimax risk for the Markovian observation model described above satisfies 
\begin{equation}\label{eqn:minimax-relation-markovian}
\inf_{\thetahat} \sup_{|\theta^\star|\leq \radius} 
\E \big[(\thetahat - \thetastar)^2 \big]
\asymp \Phi_T(\radius, \sigma) \defn \E\bigg[\Big(\frac{1}{\radius^2} + \frac{z^\T M z}{\sigma^2}\Big)^{-1}\bigg].
\end{equation}
\ecor 
\noindent See Appendix~\ref{app:proof-minimax-relation-markovian} for details of this calculation. 

Note that in the result above, the expectation on the lefthand side is
over the dataset $\{(x_i, y_i)\}_{i=1}^T$, under the Markovian model~\eqref{eqn:markov-covariates} for the covariates, and the expectation on the
righthand size is over the Gaussian vector $z = (z_1, \dots, z_T) \sim
\Normal{0}{I_{T}}$. Corollary~\ref{cor:markov} gives one example of how our
general results can even establish sharp rates for regression problems
of the form described in Section~\ref{sec:IID-regression-models}, but
with additional dependence among the covariates.

\subsection{Applications to infinite-dimensional and nonparametric models}
\label{sec:examples-nonparametric}

In this section, we derive some of the consequences of our main results for 
infinite-dimensional models, such as those arising in nonparametric regression. 
The basic idea will be to identify an infinite dimensional parameter space 
$\Theta$, typically lying in the Hilbert space $\ell^2(\N)$. We then find a nested sequence of subsets 
\begin{align*}
\Theta_1 \subset \Theta_2 \subset \cdots \subset \Theta_k \subset \cdots \subset \Theta, 
\end{align*}
where $\Theta_k$ are finite-dimensional truncations of $\Theta$. 
Under regularity conditions, we can show that the minimax risk for the $k$-dimensional 
problems converge to the minimax risk for the infinite dimensional problem 
as $k \to \infty$. Thus, since we have determined the minimax risk for each subset 
$\Theta_k$ up to universal constants (importantly, constants independent of the underlying 
dimension), we take the limit of our functional in the limit $k \to \infty$ to obtain 
a tight characterization of the minimax risk for the infinite-dimensional set $\Theta$.

In the next few sections, we carry this program out in a few examples.
We begin with a study of the canonical Gaussian sequence model in
Section~\ref{sec:Gaussian-sequence-model}.  We then turn, in
Sections~\ref{sec:examples-RKHS}
and~\ref{sec:examples-nonparametric-shift}, to nonparametric
regression models arising from reproducing kernel Hilbert spaces.  In
this setting, we are able to derive some classical results for Sobolev
spaces, derive new and sharper forms of bounds on nonparametric
regression with covariate shift, and obtain new results for random
design nonparametric models with non-uniform covariate laws.

\subsubsection{Gaussian sequence model}\label{sec:Gaussian-sequence-model}

In the canonical Gaussian sequence model, we make a countably infinite
sequence of observations of the form
\begin{equation}
\label{eqn:gauss-sequence-model}
y_i = \theta^\star_i + \eps_i z_i, \qquad \mbox{for}~i = 1, 2, 3,
\ldots
\end{equation}
Here the variables $\{z_i \}$ are a sequence of \iid{}
standard Gaussian variates, and $\eps \defn \{\eps_i\}$
indicate the noise level (\ie the standard deviation)
of the entries of the observation $y$.
It is typically assumed that there is a nondecreasing sequence
of divergent, nonnegative numbers $a \defn \{a_i\}$ and radius $C > 0$ such
that
\begin{align*}
\thetastar \in \Theta(a, C) \defn \Big\{\, \theta \in \R^\N : \sum_{j \geq 1} a_j^2 \theta_j^2 
\leq C^2 \,\Big\}.
\end{align*}
The minimax risk for this problem is then defined by 
\begin{align*}
\mathfrak{M}\Big(\eps, a, C\Big) 
\defn \inf_{\thetahat} \sup_{\thetastar \in \Theta(a, C)} 
\E\Big[ \sum_{j=1}^\infty (\thetahat_j(y) - \thetastar_j)^2 \Big],
\end{align*}
where the expectation is over $y$ according to the observation model~\eqref{eqn:gauss-sequence-model}. 

Let us define a $k$-dimensional truncation, 
\begin{align*}
\Theta_k(a, C) \defn \Big\{\, \theta \in \Theta(a, C) : \theta_{j} = 0, ~\text{for all}~
j > k \,\Big\}. 
\end{align*}
Evidently $\Theta_k(a, C)$ may be regarded as a subset of $\R^k$.
Note that the class $\{\Theta_k(a, C)\}_{k \geq 1}$ forms a nested sequence of 
subsets within $\Theta$. Moreover, we can define the minimax risk for the $k$-dimensional 
problem 
\begin{align*}
\mathfrak{M}_k\Big(\eps, a, C\Big) 
\defn \inf_{\thetahat} \sup_{\thetastar \in \Theta_k(a, C)} 
\E\Big[ \sum_{j=1}^k (\thetahat_j(y) - \thetastar_j)^2 \Big].
\end{align*}
Slightly abusing notation, above we regard $y, \thetastar \in \R^k$, where 
$y$ is distributed as the first $k$ components of the observation model~\eqref{eqn:gauss-sequence-model}. 
Then, this sequence of minimax risks satisfies the limit relation 
\begin{equation}\label{eqn:limit-relation-gauss-seq-model}
\lim_{k \to \infty} \mathfrak{M}_k\Big(\eps, a, C\Big) 
= 
\mathfrak{M}\Big(\eps, a, C\Big).
\end{equation}
See Appendix~\ref{app:proof-limit-relation-gauss-seq-model} for a proof of this relation. The $k$-dimensional problem 
can be seen as a special case of our operator model~\eqref{eqn:operator-model},
with parameters $T^{(k)}, \NoiseCovariance^{(k)}, \EstimationMat^{(k)}, \radius^{(k)}, 
\ConstraintMat^{(k)}$ defined as,
\begin{equation}\label{eqn:parameters-k-trunc-gauss-sequence-model}
\begin{gathered}
T^{(k)}(\xi) \equiv I_k, 
\qquad  
\NoiseCovariance^{(k)} = \diag(\eps_1^2, \dots, \eps_k^2), 
\qquad \EstimationMat^{\!\!(k)} = I_k, \\ 
\qquad \ConstraintMat^{\!\!(k)} = \diag\Big(\frac{1}{a_1^2}, 
\dots, \frac{1}{a_k^2}\Big), 
\quad \mbox{and},\quad 
\radius^{(k)} = C. 
\end{gathered}
\end{equation}
Computing the functional~\eqref{eqn:general-functional} for the $k$-dimensional problem,
we find it is equal to 
\begin{equation}\label{eqn:functional-k-trunc-gauss-seq-model}
R_k^\star\Big(\eps, a, C\Big) \defn 
\sup_{\tau_1, \dots, \tau_k} \Big\{\, \sum_{j=1}^k \frac{\tau_j^2 \eps_j^2}{
  \tau_j^2 + \eps_j^2
} : \sum_{j=1}^k \tau_j^2 a_j^2 \leq C^2 \,\Big\}.
\end{equation}
Hence, define the following functional 
of $\eps \defn \{\eps_j\}_{j \geq 1}, a \defn \{a_j\}_{j\geq 1}$, and $C >0$,
\begin{equation}\label{eqn:functional-full-gauss-seq-model}
R^\star(\eps, a, C) 
\defn 
\sup_{\tau = \{\tau_j\}_{j=1}^\infty} 
\Big\{\, \sum_{j=1}^\infty \frac{\tau_j^2 \eps_j^2}{
  \tau_j^2 + \eps_j^2
} : \sum_{j=1}^\infty \tau_j^2 a_j^2 \leq C^2 \,\Big\}. 
\end{equation}
Then our main results, Theorems~\ref{thm:main-upper} and~\ref{thm:main-lower} 
imply the sandwich relation
\begin{equation}\label{eqn:sandwich-relation-gauss-sequence-model}
\frac{1}{4} \, R^\star(\eps, a, C) 
\leq
\mathfrak{M}\Big(\eps, a, C\Big) 
\leq 
R^\star(\eps, a, C). 
\end{equation}
See Appendix~\ref{app:proof-sandwich-relation-gauss-seq-model} for
verification of this relation as a consequence of our results.  Note
that this recovers a well-known result for the Gaussian sequence
model~\cite{Tsy09,Johnstone2019}.  Some previous
work~\cite{DonEtAl90} has shown that the lower bound constant can be
slightly improved to $\frac{1}{1.25}$ by arguments specific to the
Gaussian sequence model. Importantly, the Gaussian sequence model is a
``deterministic'' operator model in the sense that the operator
$\LinOp{\xi}$ has no dependence on $\xi$ for this problem. The next
few examples show some consequences of our theory for
infinite-dimensional problems where the corresponding operator
$\LinOp{\xi}$ is truly random.

\subsubsection{Nonparametric regression over reproducing kernel Hilbert spaces (RKHSs)} 
\label{sec:examples-RKHS}

In this section, we consider a nonparametric regression model of the
form
\begin{align}
\label{eqn:RKHS-observations}
y_i & = f^\star(x_i) + w_i, \quad \mbox{for}~i = 1, \ldots, \numobs.
\end{align}
We assume that $\{x_i\}_{i=1}^\numobs$ are \iid{} samples covariate
law $P$ and $w_i$ being conditionally centered with conditional
variance bounded above by $\sigma^2$.  Equivalently, the noise
variables are drawn from a conditional distribution satisfying the
noise conditions~\ref{ass:centering} and~\ref{ass:bounded} with
$\NoiseCovariance = \sigma^2 I_\numobs$.\footnote{The discussion below
is unaffected by imposing additional structure on the noise, so long
as the family of possible noise distributions includes $w \sim
\Normal{0}{\sigma^2 I_\numobs}$.}  We will assume that $f^\star$ lies
in a reproducing kernel Hilbert space $\cH$, and has bounded Hilbert
norm $\vecnorm{f^\star}_{\cH} \leq \radius$. The goal is to estimate
$f^\star$.

\myparagraph{Relating the RKHS observation
  model~\eqref{eqn:RKHS-observations} with the
  model~\eqref{eqn:iid-regression}} We now show that the observation
model when $f^\star \in \cH$ is an infinite-dimensional version of the
observation model~\eqref{eqn:iid-regression}, as can be made precise
with RKHS theory.  Indeed, fix a measure space $(\cX, \cA, \nu)$, and
a measurable positive definite kernel $k \colon \cX \times \cX \to \R$
and let $\cH$ denote its reproducing kernel Hilbert
space~\cite{Aro50}.  Under mild regularity assumptions\footnote{The
elliptical representation~\eqref{eqn:elliptical-rep} is available in
great generality.  Indeed, a sufficient condition is for the map $x
\mapsto \sqrt{k(x, x)}$ to lie in $L^2(\nu)$.  It can be
shown~\cite[see Lemma 2.3]{SteSco12} that in this case, $\cH$
compactly embeds into $L^2(\nu)$ and that there is a series expansion
\begin{align}
\label{eqn:mercer-representation}
k(x, x') = \sum_{j=1}^\infty \mu_j \phi_j(x) \phi_j(x'), \quad
\mbox{for any}~x, x'\in \cX.
\end{align}
Here $\{\mu_j\}_{j=1}^\infty $ denotes a summable sequence of
non-negative eigenvalues, whereas the sequence
$\{\phi_j\}_{j=1}^\infty$ is an orthonormal family of functions $\cX
\to \R$ that lie in $L^2(\nu)$. Finally, the series converges
absolutely, for each $x, x' \in \cX$.  Note that the
infinite-dimensional series representation~\eqref{eqn:elliptical-rep}
of $\cH$ follows from the series expansion of the underlying
kernel~\eqref{eqn:mercer-representation}; see Cucker and
Smale~\cite{CucSma02} for details.}, the RKHS $\cH$ can be put into
one-to-one correspondence with a mapping of $\ell^2(\N)$. Formally, we
have
\begin{align}
\label{eqn:elliptical-rep}
\cH = \Big\{\,f \defn \sum_{j=1}^\infty \theta_j \sqrt{\mu_j} \phi_j \mid
\sum_{j=1}^\infty \theta_j^2 < \infty \,\Big\}.
\end{align}
for a nonincreasing sequence $\mu_j \to 0$ as $j \to \infty$, 
and for an orthonormal sequence $\{\phi_j\}$ in $L^2(\nu)$.  
This allows us to equivalently write the observations~\eqref{eqn:RKHS-observations} 
in the form
\begin{align}\label{eqn:equivalent-representation-RKHS}
y_i  = \langle \theta^\star, \Phi(x_i) \rangle 
+ w_i, \quad \mbox{for}~~i=1,\ldots, \numobs.
\end{align}
Above, we have defined the sequence $\theta^\star \defn (\thetastar_j)_{j=1}^\infty$
and ``feature map'' $\Phi(x) \in \ell^2(\N)$, by the formulas
\begin{align*}
\theta^\star_j \defn 
\frac{ 
\int_{\cX} f^\star(x) \phi_j(x) \, \ud \nu(x)}{\sqrt{\mu_j}}, 
\quad \mbox{and} \quad 
\big(\Phi(x)\big)_j \defn 
\sqrt{\mu_j} \phi_j(x), 
\qquad \mbox{for all}~j \geq 1. 
\end{align*}
With these definitions, note that the inner product in equation~\eqref{eqn:equivalent-representation-RKHS} is taken in the sequence space $\ell^2(\N)$. 
From the display~\eqref{eqn:equivalent-representation-RKHS}, 
we see that the RKHS observation model~\eqref{eqn:RKHS-observations} 
is in fact an infinite-dimensional version 
of the observation model~\eqref{eqn:iid-regression}. 
The remainder of this section is devoted to deriving consequences of our 
results for this model by various truncation and limiting arguments. 

\myparagraph{Truncation argument for RKHS minimax risks} Given the
RKHS ball $\ball{\radius} \defn \big \{ \, g \in \cH : \vecnorm{g}_\cH
\leq \radius \, \big\}$, our goal is to characterize the minimax risk
\begin{align}
\label{eqn:minimax-rate-RKHS}
\minimaxrisk_\numobs(\radius, \sigma^2, P) & \defn \inf_{\hat
  f} \sup_{\substack{f^\star \in \ball{\radius} \\ 
  \nu \in \cP(\sigma^2 I_\numobs)}}  \E \Big[\vecnorm[\big]{\hat f -
    f^\star}_{L^2(\nu)}^2\Big].
\end{align}
It should be noted here that the covariates are drawn from $P$ and 
the error is measured in $L^2(\nu)$. In classical work on estimation over RKHSs, it is 
typical to assume that $P = \nu$. However, we develop in this section and in 
Section~\ref{sec:examples-nonparametric-shift} some interesting consequences of 
our theory when $P \neq \nu$, and so this generality is important for our discussion. 

To apply our results to this setting, we need to define certain
finite-dimensional truncations. We start by defining
\begin{align*}
\cH_k \defn \Big\{\, f \defn \sum_{j=1}^\infty \theta_j \sqrt{\mu_j} \phi_j \mid
\theta_j = 0,~~\mbox{for all}~j > k \, \Big\}.
\end{align*}
We then define the minimax risk over the the ball $\ball{\radius}$ restricted to 
$\cH_k$, 
\begin{align}
\label{eqn:minimax-rate-k-trunc-RKHS}
\minimaxrisk^{(k)}_\numobs(\radius, \sigma^2, P) & \defn \inf_{\hat
  f} \sup_{\substack{f^\star \in \ball{\radius} \cap \cH_k \\ 
  \nu \in \cP(\sigma^2 I_\numobs)}}  \E \Big[\vecnorm[\big]{\hat f -
    f^\star}_{L^2(\nu)}^2\Big].
\end{align}
In analogy to the limit relation~\eqref{eqn:limit-relation-gauss-seq-model}
for the Gaussian sequence model, we can show that 
\begin{equation}\label{eqn:limit-relation-rkhs-model}
\lim_{k \to \infty} \minimaxrisk^{(k)}_\numobs(\radius, \sigma^2, P)
= \minimaxrisk_\numobs(\radius, \sigma^2, P).
\end{equation}
See Appendix~\ref{app:proof-limit-relation-rkhs-model} for a proof of
this relation.  The $k$-dimensional problem associated with the
risk~\eqref{eqn:minimax-rate-k-trunc-RKHS} can be seen, using the
representation~\eqref{eqn:equivalent-representation-RKHS}, as a
special case of our \iid{} observation
model~\eqref{eqn:iid-regression}, with parameters, $P, \radius,
\sigma$ and
\begin{equation}
\label{eqn:parameters-k-trunc-rkhs-model}
\psi(x) = \Phi_k(x) \defn \big(\sqrt{\mu_j}\phi_j(x)\big)_{j=1}^k,
\quad \EstimationMat = M_k \defn \diag(\mu_1, \dots, \mu_k), \quad
\mbox{and} \quad \ConstraintMat = I_k.
\end{equation}
Let us define the $k \times k$ empirical covariance matrix 
\begin{align*}
\EmpCov^{(k)} \defn \frac{1}{\numobs} \sum_{i=1}^\numobs \Phi_k(x_i) \otimes 
\Phi_k(x_i).
\end{align*}
Then the using~\eqref{eqn:parameters-k-trunc-rkhs-model}, we see that the 
functional~\eqref{eqn:d-functional-iid} for the $k$-dimensional problem is equal
to 
\begin{equation}\label{eqn:functional-k-trunc-rkhs-model}
d_\numobs^{(k)} \defn 
\sup_{\Omega \succ 0} \, \Big\{\trace \E_{P^\numobs}\!\big[ M_k^{1/2} (\EmpCov^{(k)} + \Omega^{-1} )^{-1}M_k^{1/2} \big]  : \trace(\Omega) \leq \frac{\numobs \radius^2}{\sigma^2}\Big\}
\end{equation}

\myparagraph{Characterizations of RKHS minimax risks of estimation} 
We now state the consequence of our results for the rate of estimation~\eqref{eqn:minimax-rate-RKHS}. 

\bcor \label{cor:RKHS-minimax-rate}
Define $d_{\numobs}^\star = \limsup_{k \to \infty} d_\numobs^{(k)}$, 
where the sequence $\{d_\numobs^{(k)}\}_{k \geq 1}$ is 
defined in display~\eqref{eqn:functional-k-trunc-rkhs-model}. 
Then the RKHS minimax risk satisfies 
satisfies the inequalities, 
\begin{equation}\label{eqn:sandwich-relation-rkhs-model}
\frac{1}{4}\, \frac{\sigma^2}{\numobs} d_\numobs^\star
\leq \minimaxrisk_\numobs(\radius, \sigma^2, P) \leq 
\frac{\sigma^2}{\numobs} d_\numobs^\star.
\end{equation}
\ecor
\noindent Note that this result is an immediate consequence of
Theorems~\ref{thm:main-upper} and~\ref{thm:main-lower}, together with
the limit relation~\eqref{eqn:limit-relation-rkhs-model}.

Let us now further simplify the
characterization~\eqref{eqn:sandwich-relation-rkhs-model} in the
classical situation where $P = \nu$. We can give an explicit
calculation of the minimax risk as a function of the kernel
eigenvalues $\{\mu_j\}$, using Proposition~\ref{prop:to-population},
under the additional assumption that the map $x \mapsto k(x, x)$ is
essentially bounded by a finite number $\kappa$ under $P$.  Let us
define two parameters $\lambda_\numobs^\star,
\overline{d}_\numobs^\star$ by
\begin{subequations}
\label{eqn:loosened-quantities}
\begin{align}
\sum_{k = 1}^\infty \frac{1}{\sqrt{\mu_k}}\Big(\lambda_\numobs^\star -
\frac{1}{\sqrt{\mu_k}}\Big)_+ &= \frac{\numobs \radius^2}{\sigma^2},
\quad
\mbox{and,} \label{eqn:lambda-eqn-pop}\\ \overline{d}_\numobs^\star
&\defn \sum_{k=1}^\infty \frac{1}{\lambda_\numobs^\star}
\Big(\lambda_\numobs^\star - \frac{1}{\sqrt{\mu_k}}\Big)_+.
\end{align}
\end{subequations}
When $\kappa < \infty$, the
characterization~\eqref{eqn:sandwich-relation-rkhs-model} can be
further simplified as
\begin{align}\label{eqn:loosened-characterization-rkhs}
\frac{1}{4} \, \frac{\sigma^2 }{\numobs} \overline{d}_\numobs^\star
\leq \minimaxrisk_\numobs(\radius, \sigma^2, P) \leq \Big(1 +
\frac{\kappa^2 \radius^2}{\sigma^2}\Big) \frac{\sigma^2 }{\numobs}
\overline{d}_\numobs^\star.
\end{align}
It should be noted that
relations~\eqref{eqn:sandwich-relation-rkhs-model}
and~\eqref{eqn:loosened-characterization-rkhs} establish the
nonasymptotic minimax risk of estimation for the RKHS ball of radius
$\rho$, apart from universal constants, in fairly general fashion.
The loosened inequalities~\eqref{eqn:loosened-characterization-rkhs}
permit easier calculation, but require $P = \nu$, $P$-essential
boundedness of the diagonal of the kernel, and the signal-to-noise
ratio $\tfrac{\rho}{\sigma} \lesssim \tfrac{1}{\kappa}$.  Indeed,
compared with~\eqref{eqn:sandwich-relation-rkhs-model}, the key
quantity $\overline{d}_\numobs^\star$
in~\eqref{eqn:loosened-characterization-rkhs} can be easier to
compute. The cost is that we require additional assumptions and gain
the additional prefactor $(1 + \frac{\kappa^2 \radius^2}{\sigma^2})$,
which can be large when the signal-to-noise ratio is large. Although
we have suppressed the dependence of $\lambda_\numobs^\star,
\overline{d}_\numobs^\star$ on the parameters $\sigma, \radius$ in the
notation, it should be noted that they do vary with $\sigma, \radius$
in general; see display~\eqref{eqn:loosened-quantities}.  Leveraging
our main results, we present the proofs of the
characterizations~\eqref{eqn:sandwich-relation-rkhs-model}
and~\eqref{eqn:loosened-characterization-rkhs}, respectively, in
Appendices~\ref{app:proof-sandwich-relation-rkhs-model}
and~\ref{app:proof-loosened-characterization-rkhs-model}.

Interestingly, we note that our characterizations---even the loosened
characterization~\eqref{eqn:loosened-characterization-rkhs}---does not
need the kernel to satisfy an additional eigenvalue decay condition.
Indeed, our results hold even if the kernel eigenvalues do not satisfy
the requirement of a \emph{regular kernel} as proposed in prior
work~\cite{YanEtAl17}.

Finally, we mention that---as a sanity check, classical results can be
easily derived from~\eqref{eqn:loosened-characterization-rkhs}. To
provide one concrete example, when $P = \nu$ is the uniform
distribution on $[0, 1]^d$, and $\cH$ is the Sobolev space of order
$\beta > d/2$, it can be shown that $\tfrac{\sigma^2
  \overline{d}_\numobs^\star}{\numobs} \asymp \radius^2
(\tfrac{\sigma^2}{\numobs \radius^2})^{\tfrac{2\beta}{2\beta + d}}$.
This recovers the classical minimax risk of estimation over this
function class~\cite{IbrKha80, Sto82}.  We defer this calculation to
Appendix~\ref{app:proof-sobolev-rate-calculation}, making use
of~\eqref{eqn:loosened-characterization-rkhs}.

\subsubsection{Kernel regression under covariate shift} 
\label{sec:examples-nonparametric-shift}

We now discuss one important case in which we have $P \neq \nu$ in the
RKHS model~\eqref{eqn:RKHS-observations}.  In the setting of covariate
shift, the model~\eqref{eqn:RKHS-observations} comprises of covariates
$x_i$ drawn from a \emph{source} distribution $P$ that is different
from the \emph{target} distribution $Q$ of covariates on which
estimates of the regression function are to be deployed. In this
setting, then we take $\nu = Q$ and $P \neq Q$.

For any such pair, following the argument given previously in
Section~\ref{sec:examples-nonparametric}, we find that
\begin{equation}
\label{eqn:minimax-rate-sharp-covshift}
\inf_{\hat f} \sup_{f^\star \in \ball{\radius}}
\E\Big[\vecnorm[\big]{\hat f - f^\star}_{L^2(Q)}^2\Big] \asymp
\frac{\sigma^2}{\numobs} \limsup_{k \to \infty} d_\numobs^{(k)},
\end{equation}
where the quantity $d_\numobs^{(k)}$ is defined as in
display~\eqref{eqn:functional-k-trunc-rkhs-model}.  Above, the
expectation on the lefthand side is over the noise and the covariates
drawn from $P$ as described by the
model~\eqref{eqn:RKHS-observations}.  Note that the eigenvalues
$\{\mu_j\}_{j \geq 1}$ here correspond to the diagonalization of the
integral kernel operator under the target distribution $Q$.

Let us now compare to past work due to Ma et al.~\cite{MaPatWai22},
who studied the covariate shift problem in RKHSs.  In contrast to this
work, our result is \emph{source-target distribution-dependent}: it
characterizes, apart from universal constants, the minimax risk for
any kernel, any radius, any noise level, and any covariate shift pair
$(P, Q)$.  By contrast, the results in the paper~\cite{MaPatWai22}
consider a more restrictive setup in which pair $(P, Q)$ satisfy an
absolute continuity condition ($Q \ll P$), and moreover, the
likelihood ratio is $P$-essentially bounded, meaning that there exists
some $B \in [1, \infty)$ such that
\begin{equation}
\label{eqn:def-b-bounded-pair}
\frac{\ud Q}{\ud P}(x) \leq B, \quad \mbox{for $P$-almost every $x$}.
\end{equation}
Let $d_\infty(P, Q)$ denote the $P$-essential supremum of the
likelihood ratio $\ud Q/\ud P$ when $Q \ll P$ and $d_\infty(P, Q) =
+\infty$ otherwise.  ``Uniform'' results, where minimax risks of
estimation are studied over families of covariate shifts $P$ relative
to $Q$ where $d_\infty(P, Q) \leq B$ for some parameter $B$ can be
derived as a corollary to the sharper rate
description~\eqref{eqn:minimax-rate-sharp-covshift}.

To give one simple and concrete illustration of this, we will show how
one can derive Theorem 2 in the paper~\cite{MaPatWai22}. By Jensen's
inequality, we have
\begin{align}\label{eqn:k-trunc-lower-covshift}
d_\numobs^{(k)} \geq \sup_{\Omega \succ 0} \, \Big\{\trace
(\E_{P^\numobs} M_k^{-1/2} \EmpCov^{(k)} M_k^{-1/2} + \Omega^{-1}
)^{-1} : \trace(M_k^{-1} \Omega) \leq \frac{\numobs
  \radius^2}{\sigma^2}\Big\}.
\end{align}
If $P$ satisfies $d_\infty(P, Q) \leq B$, then it follows that we have
the ordering
\begin{equation}\label{eqn:cov-shift-lower-bound-ordering}
\E_{P^\numobs} M_k^{-1/2} \EmpCov^{(k)} M_k^{-1/2} \succeq \frac{1}{B}
I_k.
\end{equation}
Moreover, this lower bound can be achieved by a shift $P$ whenever the
zero sets of the eigenfunctions $\phi_j$ in $L^2(Q)$ of the integral
operator associated with the kernel $k$ have nontrivial
intersection. Equivalently, when there exists
\begin{equation}\label{eqn:zero-set-condition}
x_0 \in \bigcap_{j \geq 1} \phi_j^{-1}(\{0\}),
\end{equation}
then the bound~\eqref{eqn:cov-shift-lower-bound-ordering} is achieved
by the distribution $P_{x_0} \defn \frac{1}{B} Q + \Big(1 -
\frac{1}{B}\Big) \delta_{x_0}$.  This choice is evidently a
$B$-bounded shift relative to $Q$.  To give an example where the zero
set condition~\eqref{eqn:zero-set-condition} holds, note that in the
case of where the kernel $k$ is associated with the periodic
$\beta$-order Sobolev class on $[0, 1]$ and $Q$ is the uniform law on
$[0, 1]$, one can take $x_0 = 0$ as the eigenfunctions are sinusoids.

Now, combining relations~\eqref{eqn:minimax-rate-sharp-covshift}
and~\eqref{eqn:k-trunc-lower-covshift} with the choice of $P =
P_{x_0}$ given above, we have
\begin{align}
\sup_{P : d_\infty(P, Q) \leq B} \inf_{\hat f} \sup_{f^\star \in
  \ball{\radius}} \E\Big[\vecnorm[\big]{\hat f -
    f^\star}_{L^2(Q)}^2\Big] &\gtrsim \frac{\sigma^2}{\numobs}
\sup_{\omega \succ 0} \Big\{ \sum_{j=1}^\infty \frac{B
  \omega_j}{\omega_j + B} : \sum_{j=1}^\infty
\frac{\omega_j}{\lambda_j} = \frac{\numobs \radius^2}{\sigma^2} \Big\}
\nonumber \\ &\asymp \radius^2 \, \sup_\lambda \Big\{
\sum_{j=1}^\infty \twomin{\frac{\sigma^2 B}{\numobs \radius^2}}{
\lambda_j \mu_j} : \lambda_j \geq 0,~\sum_{j=1}^\infty \lambda_j
= 1\, \Big\}.
\label{eqn:loose-covshift-lower}
\end{align}
Suppose, following the paper~\cite{MaPatWai22}, we additionally impose
a regularity condition on the decay of the eigenvalues $\mu_j$ of
kernel integral operator in $L^2(Q)$.  Namely, that there exists a
constant $c \in (0, \infty)$ such that
\begin{equation}\label{eqn:regular-kernel}
\sup_{\delta > 0} \frac{\sum_{j > d(\delta)} \mu_j}{\delta^2
  d(\delta)} \leq c, \quad \mbox{where} \quad d(\delta) \defn \inf\{j
\geq 1: \mu_j \leq \delta^2\}.
\end{equation}
Under this condition, we can further lower
bound~\eqref{eqn:loose-covshift-lower}, up to universal constants, by
\begin{equation}\label{eqn:very-loose-covshift-lower}
\radius^2 \, \inf_{\delta > 0} \Big\{\delta^2 + \frac{\sigma^2
  B}{\radius^2 \numobs} d(\delta) \Big\}.
\end{equation}
The details of this calculation can be found in
Appendix~\ref{app:proof-of-covshift-relation}.  Note that by
establishing the lower bound~\eqref{eqn:very-loose-covshift-lower}, we
have recovered Theorem 2 from the paper~\cite{MaPatWai22}.  We remark
that---as seen from the steps taken to arrive at this lower
bound---our more general determination of the minimax
rate~\eqref{eqn:minimax-rate-sharp-covshift} is sharper in that it
holds for a fixed pair $(P, Q)$ rather than uniformly over the larger
class $\{P : d_\infty(P, Q) \leq B\}$. Moreover, our result, as
compared to the work~\cite{MaPatWai22}, requires fewer regularity
assumptions on the underlying kernel and its diagonalization in the
target Hilbert space $L^2(Q)$.  In fact, as demonstrated in
Appendix~\ref{app:proof-of-covshift-relation}, the regularity
condition~\eqref{eqn:regular-kernel} is \emph{not} necessary for us to
establish the lower bound~\eqref{eqn:very-loose-covshift-lower}.

\section{Proofs of Theorems~\ref{thm:main-upper} and~\ref{thm:main-lower}}
\label{sec:proof-main}

In this section, we present the proofs of our main results.
In~\Cref{sec:proof-upper}, we provide the proof of our minimax upper
bound (cf.~\Cref{thm:main-upper}). In~\Cref{sec:proof-lower}, we
provide the proof of our minimax lower bound. Some calculations and
routine verifications are deferred to~\Cref{app:proofs}.

\subsection{Proof of~\Cref{thm:main-upper}}
\label{sec:proof-upper}

In this section, we develop an upper bound on the minimax risk.  In
order to do so, so, we define the risk function
\begin{align*}
r(\thetahat, \thetastar) \defn \sup_{\nu \in \cP(\NoiseCovariance)}
\E_{(\xi, w) \sim \bP \times \nu}
\E\Big[\vecnorm[\big]{\thetahat(\LinOp{\xi}, \LinOp{\xi} \theta^\star
    + w) - \thetastar}_\EstimationMat^2\Big].
\end{align*}
defined for any measurable estimator $\thetahat$ of $(\LinOp{\xi},
y)$, and any $\thetastar \in \Theta(\radius, \ConstraintMat)$.
Evidently, the minimax risk we are bounding is then expressible as
\begin{equation}
\label{eqn:minimax-risk-relation}
\mathfrak{M}(T, \bP, \NoiseCovariance, \radius, \EstimationMat,
\ConstraintMat) = \inf_{\thetahat} \sup_{\thetastar \in
  \Theta(\radius, \ConstraintMat)} r(\thetahat, \thetastar).
\end{equation}
In order to derive an upper bound, we restrict our focus to estimators
that are \emph{conditionally linear}.  Formally, we consider the class
of procedures
\begin{equation}
\label{eqn:linear-estimator}
\thetahat_C(\LinOp{\xi}, y) \defn C(\LinOp{\xi}) \LinOp{\xi}^\T
\NoiseCovariance^{-1} y,
\end{equation}
where $C$ is a $\R^{\dimension \times \dimension}$-valued measurable
function of $\LinOp{\xi}$.  Our strategy involves the following three
steps:
\begin{itemize}
  \item [(i)] First, we compute the supremum risk over the parameter
    set $\Theta(\radius, \ConstraintMat)$ and all $\nu \in
    \cP(\NoiseCovariance)$.
  \item [(ii)] Second, compute the minimizer of the supremum risk in
    the choice of $C$ in~\eqref{eqn:linear-estimator}.
  \item [(iii)] Finally, by using the curvature of the supremum risk
    and appealing to a min-max theorem, we put the pieces together to
    determine the final minimax risk.
\end{itemize}
The following subsections are devoted to the details associated with
each of these three steps. In all cases, we defer routine calculations
and verification to~\Cref{app:proofs-upper}.

\subsubsection{Supremum risk of estimator $\thetahat_C$}

Starting with the definition~\eqref{eqn:linear-estimator}, for any
matrix $C$, we have
\begin{align*}
\thetahat_C - \thetastar = (C(\LinOp{\xi}) \LinOp{\xi}^\T
\NoiseCovariance^{-1} \LinOp{\xi} - I_\dimension) \thetastar +
C(\LinOp{\xi}) \LinOp{\xi}^\T \NoiseCovariance^{-1} w.
\end{align*}
Therefore, the risk $r(\thetahat_C, \theta^\star)$ associated with
$\thetahat_C$ can be bounded as
\begin{align}
r(\thetahat_C, \thetastar) &\defn \sup_{\nu \in \cP(\NoiseCovariance)}
\E\Big[\vecnorm{\thetahat_C(X, y) - \thetastar}_{\EstimationMat}^2
  \Big] \nonumber \\
& = \trace \bigg\{ \EstimationMat^{\!1/2} \E_\xi \Big[(C(\LinOp{\xi})
  \LinOp{\xi}^\T \NoiseCovariance^{-1} \LinOp{\xi} - I_\dimension)
  \thetastar \otimes \thetastar (C(\LinOp{\xi}) \LinOp{\xi}^\T
  \NoiseCovariance^{-1} \LinOp{\xi} - I_\dimension)^\T \nonumber \\
& \qquad\qquad\quad + C(\LinOp{\xi}) \LinOp{\xi}^\T
  \NoiseCovariance^{-1} \LinOp{\xi} C(\LinOp{\xi})^\T \Big]
\EstimationMat^{\!1/2} \bigg\}. \label{eqn:risk-upper}
\end{align}
The equality above uses the property~\ref{ass:bounded} of
distributions $\nu \in \cP(\NoiseCovariance)$; note that it is
achieved by the Gaussian distribution $\nu =
\Normal{0}{\NoiseCovariance}$.

\subsubsection{Curvature and minimizers of the functional
  $r(\thetahat_C, \theta^\star)$}

We begin by observing that the function $r(\thetahat_C, \cdot) \colon
\Theta(\radius, \ConstraintMat) \to \R_+$ can be replaced by an
equivalent mapping---which, with a slight abuse of notation we denote
by the same symbol $r$--- on the space of symmetric positive definite
matrices of the form
\begin{align*}
\cK(\radius, \ConstraintMat) \defn \Big\{\, \Omega \succeq 0 \mid
\trace(\ConstraintMat^{\!-1/2} \Omega \ConstraintMat^{\!-1/2}) \leq
\radius^2\, \Big\}.
\end{align*}
We define (in a sense, this is can be regarded as an extension to the
set $\cK(\radius, \ConstraintMat)$)
\begin{multline}
\label{eqn:definition-linear-risk-matrix-version}
r(\thetahat_C, \Omega) \defn \trace \Big\{ \EstimationMat^{\!1/2}
\E_\xi \Big[(C(\LinOp{\xi}) \LinOp{\xi}^\T \NoiseCovariance^{-1}
  \LinOp{\xi} - I_\dimension) \Omega (C(\LinOp{\xi}) \LinOp{\xi}^\T
  \NoiseCovariance^{-1} \LinOp{\xi} - I_\dimension)^\T \\ +
  C(\LinOp{\xi}) \LinOp{\xi}^\T \NoiseCovariance^{-1} \LinOp{\xi}
  C(\LinOp{\xi})^\T \Big]\EstimationMat^{\!1/2}\Big\}.
\end{multline}
Note that $r(\thetahat_C, \thetastar) = r(\thetahat_C, \thetastar
\otimes \thetastar)$ for $\thetastar \in \Theta(\radius,
\ConstraintMat)$.  We claim that the suprema over $\Theta(\radius,
\ConstraintMat)$ and $\cK(\radius, \ConstraintMat)$ are the same.

\begin{lemma} \label{lem:suprema-equality} The suprema of the risk functional
$r$ taken over either the set $\Theta(\radius, \ConstraintMat)$ or the
  set $\cK(\radius, \ConstraintMat)$ are equal---that is, we have
\begin{align*}
\sup_{\thetastar \in \Theta(\radius, \ConstraintMat)} r(\thetahat_C,
\thetastar) = \sup_{\Omega \in \cK(\radius, \ConstraintMat)}
r(\thetahat_C, \Omega),
\end{align*}
for every conditionally linear estimator $\thetahat_C$ of the
form~\eqref{eqn:linear-estimator}.
\end{lemma} 
\noindent See Appendix~\ref{app:proof-suprema-equality} for the proof
of this claim. \\

\medskip
\noindent Our next result characterizes some properties of the mapping
$(C, K) \mapsto r(\thetahat_C, K)$.

\begin{lemma}
\label{lem:curvature-of-linear-risk}
Over the set of measurable functions $C$ and matrices $\Omega \in
\cK(\radius, \ConstraintMat)$, the mapping $(C, \Omega) \mapsto
r(\thetahat_C, \Omega)$ is affine in $\Omega$ and convex in $C$.
\end{lemma} 
\noindent See Appendix~\ref{app:proof-curvature-risk} for the proof of
this claim. \\

\medskip

\noindent Our next claim determines the minimizer of $r(\cdot,
\Omega)$ over estimators $\thetahat_C$ of the
form~\eqref{eqn:linear-estimator}, provided that $\Omega$ is strictly
positive definite.

\begin{proposition}\label{prop:minimizer-of-linear-risk}
Let $\Omega$ be a symmetric positive definite matrix. Then
\begin{equation}\label{eqn:minimizer-linear-risk}
\inf_C r(\thetahat_C, \Omega) = \trace \Big\{ \EstimationMat^{\!1/2}
\E_\xi (\Omega^{-1} + \LinOp{\xi}^\T \NoiseCovariance^{-1}
\LinOp{\xi})^{-1} \EstimationMat^{\!1/2} \Big\}
\end{equation}
Moreover, the infimum is attained with the choice $C(\LinOp{\xi}) =
(\Omega^{-1} + \LinOp{\xi}^\T \NoiseCovariance^{-1}
\LinOp{\xi})^{-1}$.
\end{proposition}
\noindent See~\Cref{app:proof-minimizer-linear-risk} for the proof.

\subsubsection{Proof of~\Cref{thm:main-upper}}

We now piece together the previous lemmas to establish our main upper
bound, as claimed in~\Cref{thm:main-upper}.  In view of the
relation~\eqref{eqn:minimax-risk-relation} and the
bound~\eqref{eqn:risk-upper}, we find that
\begin{subequations}
\begin{align}
\mathfrak{M}(T, \bP, \NoiseCovariance, \radius, \EstimationMat,
\ConstraintMat) &\leq \inf_{C} \sup_{\thetastar \in \Theta(\radius,
  \ConstraintMat)} r(\thetahat_C,
\thetastar) \label{eqn:min-max-justify-C}\\ &= \inf_{C} \sup_{\Omega
  \in \cK(\radius, \ConstraintMat)} r(\thetahat_C,
\Omega) \label{eqn:min-max-replace-matrix}\\ &= \sup_{\Omega \in
  \cK(\radius, \ConstraintMat)} \inf_{C} r(\thetahat_C,
\Omega) \label{eqn:min-max-apply-ky-fan} \\ &= \sup_{\Omega} \Big\{\,
\E \trace\Big( \EstimationMat^{\!\!1/2} (\Omega^{-1} + \LinOp{\xi}^\T
\NoiseCovariance^{-1} \LinOp{\xi})^{-1} \EstimationMat^{\!\!1/2} \Big)
: \nonumber \\ &\qquad\qquad\qquad\qquad\qquad\qquad\qquad \Omega
\succ 0,\trace(\ConstraintMat^{\!-1/2} \Omega
\ConstraintMat^{\!-1/2})\leq \radius^2
\,\Big\}. \label{eqn:upper-bound-conclusion}
\end{align}
\end{subequations}
To clarify, in the first display~\eqref{eqn:min-max-justify-C} and
below, the infimum over $C$ denotes an infimum over all
$\R^{\dimension \times \dimension}$-valued measurable functions of
$\LinOp{\xi}$. In display~\eqref{eqn:min-max-replace-matrix}, we have
applied
Lemma~\ref{lem:suprema-equality}. Relation~\eqref{eqn:min-max-apply-ky-fan}
follows from the Ky Fan min-max theorem~\cite{Fan53,BorZhu86} together
with~\Cref{lem:curvature-of-linear-risk}. Note that the set
$\cK(\radius, \ConstraintMat)$ is evidently a compact convex subset of
$\R^{\dimension \times \dimension}$.  The final
equality~\eqref{eqn:upper-bound-conclusion} is essentially an
application of~\Cref{prop:minimizer-of-linear-risk};
see~\Cref{app:proof-upper-conclusion} for the details of this
verification.

\subsection{Proof of lower bound, Theorem~\ref{thm:main-lower}}
\label{sec:proof-lower}
In this section, we prove our lower bound on the minimax risk.  In
order to do so, we focus on lower bounding the Gaussian minimax risk
\begin{align*}
\mathfrak{M}^{\rm G}(\PlainLinOp, \bP, \NoiseCovariance, \radius,
\EstimationMat, \ConstraintMat) \defn \inf_{\thetahat}
\sup_{\thetastar \in \Theta(\radius, \ConstraintMat)} \E_{(\xi, w)
  \sim \bP \times \Normal{0}{\NoiseCovariance}}
\Big[\vecnorm{\thetahat(\LinOp{\xi}, \LinOp{\xi}\thetastar + w) -
    \thetastar}^2_{\EstimationMat} \Big].
\end{align*}
Evidently, the Gaussian minimax risk lower bounds the general minimax
risk, so that we have $\mathfrak{M}^{\rm G} \leq \mathfrak{M}$.  In
Section~\ref{sec:lower-proof-reduction}, we reduce this Gaussian
minimax risk to yet another Gaussian observation model. A minimax
lower bound for this auxiliary problem is then presented as
Proposition~\ref{prop:lower-reduced-risk} in
Section~\ref{sec:reduction-and-proof}. This result is the bulk of the
proof of the lower bound, and it quickly allows us to establish our
main result, Theorem~\ref{thm:main-lower}. In
Section~\ref{sec:proof-prop-lower}, we then complete the proof of
Proposition~\ref{prop:lower-reduced-risk}.

\subsubsection{Reduction to an alternate observation model}
\label{sec:lower-proof-reduction}

To establish the lower bound, we first show that the minimax risk
associated with our estimation problem is equivalent to another,
perhaps simpler, minimax risk.

\myparagraph{An auxiliary observation model} This observation model is
defined by a random quadruple $(r, V, \Lambda, \Upsilon)$.  The triple
$(r, V, \Lambda)$ comprises a random integer $r$, a random orthogonal
matrix $V \in \R^{\dimension \times r}$ satisfying $V^\T V = I_r$, and
a random, $r \times r$ diagonal positive definite matrix $\Lambda$.
Conditional on $(r, V, \Lambda)$, the observation $\Upsilon$ is a
Gaussian random variable, satisfying the equation
\begin{equation}\label{eqn:aux-model}
\Upsilon = VV^\T \eta^\star  + V\Lambda^{-1/2} z, \quad \mbox{where}\quad z\sim \Normal{0}{I_r}. 
\end{equation}
Above, the random vector $z$ is drawn from the multivariate Gaussian
with identity covariance in $\R^r$; it is independent of $(r, V,
\Lambda)$. If $\omega \defn (r, V, \Lambda)$ is distributed according
to $\bQ$, we denote the minimax risk for this observation model as
\begin{align*}
\mathfrak{M}^{\rm G}_{\rm red}(\bQ, K) 
\defn \inf_{\etahat} \sup_{\eta \in \Theta(K)} 
\E_{(\omega, \Upsilon)} 
\Big[ \vecnorm{\etahat(\omega, \Upsilon) - \eta}_2^2\Big].
\end{align*}
Above, the expectation indexed by $(\omega, \Upsilon)$ is over $\omega \sim \bQ$ and 
$\Upsilon$ as in~\eqref{eqn:aux-model}. The infimum is over measurable functions of 
$(\omega, \Upsilon)$. The set $\Theta(K)$ is a shorthand for the set 
$\Theta(1, K) = \{\|\theta\|_K \leq 1\}$. 

\myparagraph{Reduction to the new observation model}

We formally reduce the minimax risk $\mathfrak{M}^{\rm G}$ to the reduction $\mathfrak{M}^{\rm G}_{\rm red}$, as 
follows. 
\begin{lemma} \label{lem:minimax-reduction}
Let $\tilde \bP$ denote the distribution of the triple $(r(\xi), V_{\xi}, \Lambda_\xi)$ under $\bP$, where
$r(\xi)$ is the (finite) rank of $Q_\xi = \EstimationMat^{\!-1/2} \LinOp{\xi}^\T \NoiseCovariance^{-1} \LinOp{\xi} \EstimationMat^{\!-1/2}$, and 
$Q_\xi = V_{\xi} \Lambda_\xi V_{\xi}^\T$ denotes the diagonalization of this positive definite matrix. 
Then, for any $(T, \bP, \NoiseCovariance, \radius, \ConstraintMat, \EstimationMat)$, we have 
\begin{align*}
\minimaxrisk^{\rm G}(T, \bP, \NoiseCovariance, \radius, \ConstraintMat, \EstimationMat) 
= \minimaxrisk^{\rm G}_{\rm red}(\tilde \bP, \radius^2 \EstimationMat^{\!1/2} \ConstraintMat \EstimationMat^{\!1/2} ).
\end{align*}
\end{lemma} 
\noindent See Appendix~\ref{app:proof-minimax-reduction} for a proof of this claim.

\subsubsection{Lower bounding the minimax risk}\label{sec:reduction-and-proof}

We now focus on lower bounding $\mathfrak{M}^{\rm G}_{\rm red}$. The following result is a formal 
statement of the lower bound for the ``reduced'' minimax risk. 

\begin{proposition}\label{prop:lower-reduced-risk}
  For any $\tau \in (0, 1]$ and any $\Pi \succ 0$ such that $\trace(K^{-1/2} \Pi K^{-1/2}) \leq 1$, we have
\begin{equation}\label{eqn:strong-reduced-lower-bound}
\mathfrak{M}^{\rm G}_{\rm red}(\bQ, K) 
\geq
\E \trace\Big( 
(\tfrac{1}{c(\tau, \Pi)} \Pi^{-1} + V \Lambda V^\T )^{-1}\Big),
\end{equation}
where the constant $c(\tau, \Pi)$ is defined in Lemma~\ref{lem:fisher-info-lower}. 
Moreover, we have the lower bounds
\begin{subequations}
\begin{align}
\mathfrak{M}^{\rm G}_{\rm red}(\bQ, K) 
&\geq 
\sup_{\Pi}  \Big\{\, \E \trace\Big( 
(\Pi^{-1} + V \Lambda V^\T )^{-1}\Big) 
: \Pi \succ 0,~\trace(K^{-1/2} \Pi K^{-1/2})\leq 1/4
\,\Big\} \label{eqn:strong-reduced-lb}\\ 
&\geq 
\frac{1}{4} \, \sup_{\Pi}  \Big\{\, \E \trace\Big( 
(\Pi^{-1} + V \Lambda V^\T )^{-1}\Big) 
: \Pi \succ 0,~\trace(K^{-1/2} \Pi K^{-1/2})\leq 1 
\,\Big\}. \label{eqn:weak-reduced-lb}
\end{align}
\end{subequations}
\end{proposition}

\myparagraph{Proof of Theorem~\ref{thm:main-lower}}

We take the claim of~\Cref{prop:lower-reduced-risk} as given for the
moment, and use it to derive our minimax lower bound. As mentioned, we
may restrict to Gaussian noise to establish the lower bound; formally,
we have $\mathfrak{M} \geq \mathfrak{M}^{\rm G}$. Additionally, the
reduction given in Lemma~\ref{lem:minimax-reduction} combined with the
stronger lower bound~\eqref{eqn:strong-reduced-lb} in
Proposition~\ref{prop:lower-reduced-risk} gives us
\begin{multline*}
\mathfrak{M}(T, \bP, \NoiseCovariance, \radius, \EstimationMat, \ConstraintMat) 
\\\geq \resizebox{.915\hsize}{!}{$
\sup_{\Pi}  \Big\{\, \E \trace\Big( 
(\Pi^{-1} + \EstimationMat^{\!-1/2} \LinOp{\xi}^\T \NoiseCovariance^{-1} \LinOp{\xi} \EstimationMat^{\!-1/2})^{-1}\Big) 
: \Pi \succ 0, \trace(\EstimationMat^{\!-1/2} \Pi \EstimationMat^{\!-1/2} \ConstraintMat^{-1} )\leq \tfrac{\radius^2}{4}\,\Big\}$}.
\end{multline*}
Now define the matrix $\Omega = \EstimationMat^{\!-1/2} \Pi \EstimationMat^{\!-1/2}$. Then, the quantity on the righthand side is equal to 
\begin{align*}
\sup_{\Omega}  \Big\{\, \E \trace\Big( \EstimationMat^{\!1/2} 
(\Omega^{-1} + \LinOp{\xi}^\T \NoiseCovariance^{-1} \LinOp{\xi})^{-1} \EstimationMat^{\!1/2} \Big) 
: \Omega \succ 0,~\trace(\ConstraintMat^{\!-1/2} \Omega \ConstraintMat^{\!-1/2})\leq \tfrac{\radius^2}{4} \,\Big\},
\end{align*}
which furnishes the first inequality in Theorem~\ref{thm:main-lower}. With similar manipulations
to the weaker lower bound~\eqref{eqn:weak-reduced-lb}
in Proposition~\eqref{prop:lower-reduced-risk}, or by arguing directly from the display above,
the second inequality in Theorem~\ref{thm:main-lower} follows.
In order to establish the more detailed lower bound~\eqref{eqn:gen-sharp-lower}, 
we repeat the argument above but use~\eqref{eqn:strong-reduced-lower-bound}.

\subsubsection{Proof of Proposition~\ref{prop:lower-reduced-risk}}\label{sec:proof-prop-lower}
The lower bound proceeds in five steps: 
\begin{enumerate}[label=(\roman*)]
  \item We first lower bound the minimax risk in terms of the expected conditional Bayesian risk over any prior on the 
  parameter set $\Theta(K)$.
  \item We then demonstrate that, conditionally, 
  there is a family of auxiliary Bayesian estimation problems, indexed by a parameter $\lambda > 0$, which are all no harder than the Bayesian estimation problem implied by the conditional Bayesian risk. 
  \item We compute, in closed form, the Bayesian risk for any prior and any parameter $\lambda > 0$. We are able to show that the Bayesian risk is a functional of 
  the Fisher information of the marginal distribution of the observed data under the prior and sampling model. 
  \item For each $\lambda > 0$, we then calculate a lower bound on the Fisher information for a prior obtained by conditioning a Gaussian distribution with mean zero and covariance 
  $\Pi$ to the parameter space.
  \item We put the pieces together: optimizing over all covariance operators $\Pi$, and the family of ``easier'' problems (\ie optimizing over $\lambda > 0$), we obtain our claimed lower bound. 
\end{enumerate}

Next, we present the details of the steps outlined above. Extended calculations and routine 
verification are deferred to Appendix~\ref{app:proofs-lower}.

\myparagraph{Step 1: Reduction to conditional Bayesian risk}
We begin by lower bounding the minimax risk via the Bayes risk. 
Owing to the standard relation between minimax and Bayesian risks, we have for any 
prior $\pi$ on $\Theta(K)$ that
\begin{multline}\label{ineq:minimax-to-Bayes}
\mathfrak{M}^{\rm G}_{\rm red}(\bQ, K) 
= \inf_{\etahat} \sup_{\eta \in \Theta(K)} 
\E_{(\omega, \Upsilon)} 
\Big[ \vecnorm{\etahat(\omega, \Upsilon) - \eta}_2^2\Big] 
\geq 
\inf_{\etahat} 
\E_{\eta \sim \pi} 
\E_{(\omega, \Upsilon)} 
\big[\vecnorm{\etahat- \eta}_2^2\big] \eqcolon B(\pi).
\end{multline}
The quantity $B(\pi)$ appearing above is the Bayesian risk when the parameter 
$\eta$ is drawn from the prior $\pi$. 
The following observation is key for the lower bound. 
After moving to Bayesian risks, we can condition on the ``design'', 
denoted by the random tuple $\omega = (r,V, \Lambda)$, and consider 
the conditional Bayesian risk. 
Formally, we have 
\begin{align}
B(\pi) &= \inf_{\etahat} \E_{\eta \sim \pi} 
\E_{(\omega, \Upsilon) \sim \cD_\eta} 
\Big[\vecnorm[\big]{\etahat- \eta}_2^2\Big] 
\geq \E_{\omega \sim \bQ} \bigg[
\inf_{\etahat_\omega} \E_{\eta \sim \pi} \E_{\Upsilon} \vecnorm[\big]{\etahat_\omega(\Upsilon) - \eta}_2^2 \bigg]. 
\label{eqn:conditional-Bayes-risk}
\end{align}
Above, the inequality follows by observing that if the function 
$\etahat\colon (\omega, \Upsilon)\mapsto \etahat\in \R^\dimension$ is measurable, 
then $\etahat_\omega(\Upsilon) \defn \etahat(\omega, \Upsilon)$ is a 
measurable of $\Upsilon$. Note that the infimum on the righthand side is restricted to those maps which are measurable function of $\omega$; note that 
they may depend on $\omega$, and therefore we have included a subscript depending on $\omega$ to indicate this.\footnote{In some cases, 
this inequality may hold with equality. However, to be clear, in general the inequality arises since if $\{\etahat_\omega\}_{\omega}$ 
is a family of measurable functions (of $\Upsilon$) for each $\omega$ in the support of $\bQ$, 
it is not necessarily the case that $\etahat(\omega, \Upsilon) \defn \etahat_{\omega}(\Upsilon)$ is measurable.} 
To lighten notation in the subsequent discussion, we define the \emph{conditional Bayesian risk} under $\pi$ and for 
a realization of the random variable $\omega = \omega_0$,
\begin{align*}
B(\pi \mid \omega_0) \defn \inf_{\etahat} \E_{\eta \sim \pi} \E_{z \sim \Normal{0}{I_{r_0}}} \Big[\vecnorm[\big]{\etahat(V_0 V_0^\T \eta +  V_0 \Lambda_0^{-1/2} z) - \eta}_2^2\Big], 
\quad \mbox{where}~~\omega_0 = (r_0, V_0, \Lambda_0). 
\end{align*}
Using this definition, along with the two inequalities~\eqref{ineq:minimax-to-Bayes} and~\eqref{eqn:conditional-Bayes-risk}, we have 
demonstrated 
\begin{equation}\label{eqn:final-bayesian-risk-lower}
\mathfrak{M}^{\rm G}_{\rm red}(\bQ, K) 
\geq 
\E_{\omega \sim \bQ} \big[B(\pi \mid \omega)\big], \qquad \mbox{for any prior}~\pi~\mbox{on}~\Theta(K).
\end{equation}
Therefore, it suffices for us to lower bound $B(\pi \mid \omega)$.

\myparagraph{Step 2: Reduction to a family of easier problems}

In this step, we fix a parameter $\lambda > 0$, which will index yet another auxiliary Bayesian estimation problem.
The intuition will be that as $\lambda \to 0^+$, we are ``approaching'' the difficulty of the original Bayesian estimation problem. 

Formally, fix $\omega = (r, V, \Lambda)$. Throughout we will let $V_\perp \colon \R^\dimension \to \range(V)^\perp$ denote the projection of 
an element $\eta \in \R^\dimension$ to the orthogonal complement of the closed subspace $\range(V)$.
We now consider the observation, where for an independent random Gaussian variable $z \sim \Normal{0}{I_\dimension}$
\begin{equation}\label{eqn:lambda-obs-model}
\Upsilon_\lambda = \underbrace{(VV^\T + \lambda V_\perp)}_{\eqcolon X_\lambda} \eta +  V \Lambda^{-1/2} w  + \sqrt{\lambda} V_\perp z=
X_\lambda \eta +  (V \Lambda^{-1} V^\T + \lambda V_\perp)^{1/2} w',
\end{equation}
where the last equality holds in distribution. Define $\Sigma_\lambda \defn V\Lambda^{-1} V^\T + \lambda V_\perp$; 
evidently $\Sigma_\lambda$ is a symmetric positive definite matrix for any $\lambda > 0$. Then, 
$\Upsilon_\lambda$ has distribution $\Normal{X_\lambda \eta}{\Sigma_\lambda}$. We remark that the observation $\Upsilon_\lambda$
is more convenient than $\Upsilon$ as its covariance is nonsingular and moreover its mean is a nonsingular linear transformation of 
$\eta$---note that neither of these properties hold for $\Upsilon$.

Our goal is to show that the observation $\Upsilon_\lambda$ is more ``informative'' than $\Upsilon$. 
To do this, we now define the (conditional) Bayesian risk for $\Upsilon_\lambda$, 
\begin{align*}
B_\lambda(\pi \mid \omega) \defn 
\inf_{\etahat} \Big\{ B_\lambda(\etahat, \pi \mid \omega) \defn \E\big[\vecnorm{\etahat(\Upsilon_\lambda) - \eta}_2^2\big]\Big\}. 
\end{align*}
The main claim is that this provides a lower bound on our original conditional Bayesian risk. 
\begin{lemma} \label{lem:lambda-lower}
For any $\omega$ and $\lambda > 0$, we have 
\begin{align*}
B(\pi \mid \omega) \geq B_\lambda(\pi \mid \omega). 
\end{align*}
\end{lemma} 
\noindent See Appendix~\ref{app:proof-lambda-lower} for a proof of this claim.

\myparagraph{Step 3: Calculation of Bayesian risk $B_\lambda(\pi \mid \omega)$, for a fixed prior $\pi$ and parameter $\lambda > 0$}

To compute the Bayesian risk for a fixed prior $\pi$ and parameter $\lambda > 0$, we develop a variant of 
Tweedie's formula (also sometimes referred to as Brown's identity, when applied to Bayesian risks)~\cite{Twe47,Rob56,Bro71}.

To state the result, we need to introduce some notation. 
We define the marginal and conditional densities of $\Upsilon_\lambda$---disregarding normalization constants---as,
\begin{align*}
p(y) \defn \int p(y \mid \eta) \, \pi(\ud \eta) \qquad\mbox{where} \quad 
p(y \mid \eta) \defn \exp\Big(-\frac{1}{2} \|y- X_\lambda \eta\|_{\Sigma_\lambda^{-1}}^2\Big).
\end{align*}
Finally we define the Fisher information of the marginal distribution of $\Upsilon_\lambda$, which is given by 
\begin{align*}
\Information{\Upsilon_\lambda} \defn \E[\nabla \log p(\Upsilon_\lambda) \otimes \nabla \log p(\Upsilon_\lambda)].
\end{align*}
With this notation in hand, we can now state our formula for the Bayesian risk under the prior $\pi$ and for parameter $\lambda >0$. 
\begin{lemma} \label{lem:lambda-bayes-risk}
Fix $\omega = (r, V, \Lambda)$. Define $X_\lambda \defn VV^\T + \lambda V_\perp$ and 
$\Sigma_\lambda \defn V \Lambda^{-1} V^\T + \lambda V_\perp$. 
Fix prior $\pi$, and parameter $\lambda > 0$. Then the conditional Bayesian risk is given by 
\begin{align*}
B_\lambda(\pi \mid \omega) = \trace\Big(X_\lambda^{-1} \Sigma_\lambda \big[\Sigma_\lambda^{-1} - \Information{\Upsilon_\lambda}\big]\Sigma_\lambda X_{\lambda}^{-1}\Big).
\end{align*}
\end{lemma} 
\noindent See Appendix~\ref{app:proof-lambda-bayes-risk} for a proof of this claim.

\myparagraph{Step 4: Lower bound on Fisher information for conditioned
  Gaussian prior} Consider a prior $\pi$ which is absolutely
continuous with respect to Lebesgue measure on $\R^\dimension$.
Furthermore, suppose that its Lebesgue density $f_\pi \defn \frac{\ud
  \pi}{\ud \eta}$ has logarithmic gradient almost everywhere.  Define
\begin{align*}
\Information{\pi} \defn \int \nabla \log f_\pi(\eta) \otimes \nabla
\log f_\pi(\eta) \, \ud \pi(\eta).
\end{align*}
Recall also that the Fisher information associated with a Gaussian distribution $\Normal{\mu}{\Pi}$ for nonsingular $\Pi$ 
is given by $\Pi^{-1}$~\cite[Example 6.3]{LehCas98}. 
Therefore, applying well-known results for the Fisher information~\cite[eqn.~(8) and Corollary 1]{Zam98}
\begin{equation}\label{eqn:fisher-info-ineq}
\Information{\Upsilon_\lambda} \preceq (X_\lambda \Information{\pi}^{-1} X_\lambda + \Sigma_\lambda)^{-1}.
\end{equation}

Next, we select a prior distribution and calculate the Fisher information $\Information{\Upsilon_\lambda}$ for the 
marginal density under this prior. For a parameter $\tau \in (0,
1]$ and symmetric positive definite covariance matrix $\Pi$, we define the probability measures
\begin{equation}  
  \label{eqn:display-distributions}
  \pi^{\rm G}_{\tau, \Pi} = \Normal{0}{\tau^2 \Pi} \quad \mbox{and}
  \quad \pi_{\tau, \Pi} = \pi_{\tau, \Pi}^{\rm G}\big(\cdot \mid
  \Theta(K)\big).
\end{equation}
In other words, $\pi_{\tau, \Pi}$ denotes the probability measure
$\Normal{0}{\tau^2 \Pi}$ conditioned on the constraint set. Formally, it is 
defined by the relation, 
\begin{align*}
\pi_{\tau, \Pi}(A) \defn \frac{\pi^{\rm G}_{\tau, \Pi}\big(A \cap \Theta(K)\big)}{\pi^{\rm G}_{\tau, \Pi}\big(\Theta(K)\big)},
\end{align*}
for any event $A$. 
For these priors, we have the following claim. 
\begin{lemma}\label{lem:fisher-info-lower}
Let $\tau \in (0, 1]$ and $\Pi$ be a symmetric positive definite matrix satisfying the relation 
$\trace(\Pi^{1/2} K^{-1} \Pi^{1/2}) \leq 1$. Then the Fisher information of the conditioned 
prior $\pi_{\tau, \Pi}$ satisfies the inequality
\begin{align*}
\Information{\pi_{\tau, \Pi}}^{-1} \succeq c(\tau, \Pi) \Pi, 
\end{align*}
where $c(\tau, \Pi) = \tau^2 (1 - \pi_{\tau, \Pi}^{\rm G}(\Theta(K)^c)) > 0$.
\end{lemma} 
\noindent See Appendix~\ref{app:proof-fisher-info-lower} for the proof of this claim.

\myparagraph{Step 5: Putting the pieces together}

Combining Lemmas~\ref{lem:lambda-lower} and~\ref{lem:lambda-bayes-risk} along with 
the inequality~\eqref{eqn:fisher-info-ineq} and Lemma~\ref{lem:fisher-info-lower}, we 
find that for any $\tau \in (0, 1]$ and symmetric positive definite matrix $\Pi$
satisfying $\trace(\Pi^{1/2} K^{-1} \Pi^{1/2}) \leq 1$, that 
\begin{align*}
B(\pi \mid \omega) &\geq \sup_{\lambda > 0} 
\trace\Big(X_\lambda^{-1} \Sigma_\lambda \big[\Sigma_\lambda^{-1} - 
(c(\tau, \Pi) X_\lambda \Pi X_\lambda + \Sigma_\lambda)^{-1}\big]\Sigma_\lambda X_{\lambda}^{-1}\Big) \\ 
&= \sup_{\lambda > 0} 
\trace\Big( 
(\tfrac{1}{c(\tau, \Pi)} \Pi^{-1} + X_\lambda \Sigma_\lambda^{-1} X_\lambda)^{-1}\Big). 
\end{align*}
Above, we used the relation $A(A^{-1} - (B + A)^{-1})A = (A^{-1} + B^{-1})^{-1}$, valid 
for any pair $(A, B)$ of symmetric positive definite matrices. Our particular choice 
of matrices was $A = \Sigma_\lambda$ and $B = X_\lambda$. 
Note that 
\begin{align*}
X_\lambda \Sigma_\lambda^{-1} X_\lambda = V \Lambda V^\T + \lambda V_\perp. 
\end{align*}
Therefore, by continuity, we have 
\begin{equation}\label{eqn:conditional-Bayes-lower}
B(\pi \mid \omega) \geq \lim_{\lambda \to 0^+}
\trace\Big( 
(\tfrac{1}{c(\tau, \Pi)} \Pi^{-1} + V \Lambda V^\T + \lambda V_\perp)^{-1}\Big)
= \trace\Big( 
(\tfrac{1}{c(\tau, \Pi)} \Pi^{-1} + V \Lambda V^\T )^{-1}\Big).
\end{equation}
Taking the expectation over $\omega$, and
applying our minimax lower bound~\eqref{eqn:final-bayesian-risk-lower}, we
have established lower bound~\eqref{eqn:strong-reduced-lower-bound}.
Note that since $c(\tau, \Pi) \in (0, 1]$,  we evidently have from the above display that 
\begin{align*}
B(\pi \mid \omega) \geq c(\tau, \Pi) \trace\Big( 
(\Pi^{-1} + V \Lambda V^\T )^{-1}\Big).
\end{align*}
Let us define the constant 
\begin{align*}
c_\ell(K) \defn \inf_{\substack{\Pi \succ 0 \\ \trace(\Pi K^{-1}) \leq 1}} 
\sup_{\tau \in (0, 1]} c(\tau, \Pi). 
\end{align*}
Then combining the conditional lower bound~\eqref{eqn:conditional-Bayes-lower} 
with our minimax lower bound~\eqref{eqn:final-bayesian-risk-lower}, we obtain 
\begin{align*}
\mathfrak{M}^{\rm G}_{\rm red}(\bQ, K) 
&\geq \sup_{\Pi}  
\Big\{\, \E \trace\Big( 
(\Pi^{-1} + V \Lambda V^\T )^{-1}\Big) 
: \Pi \succ 0,~\trace(\Pi^{1/2} K^{-1} \Pi^{1/2})\leq c_\ell(K) 
\,\Big\} \\ 
&=  \sup_{\Pi}  \Big\{\, \E \trace\Big( 
(\tfrac{1}{c_\ell(K)} \Pi^{-1} + V \Lambda V^\T )^{-1}\Big) 
: \Pi \succ 0,~\trace(\Pi^{1/2} K^{-1} \Pi^{1/2}) \leq 1 
\,\Big\} \\ 
&\geq c_\ell(K) \, 
 \sup_{\Pi}  \Big\{\, \E \trace\Big( 
( \Pi^{-1} + V \Lambda V^\T )^{-1}\Big) 
: \Pi \succ 0,~\trace(\Pi^{1/2} K^{-1} \Pi^{1/2})\leq 1 
\,\Big\}.
\end{align*}
To complete the proof, we simply need to lower bound the constant $c_\ell(K)$ universally. 
\begin{lemma} \label{lem:constant-lower}
The constant $c_\ell(K)$ is lower bounded, for any symmetric positive definite $K$, as 
\begin{align*}
c_\ell(K) \geq \frac{1}{4}.
\end{align*} 
\end{lemma}
\noindent See Appendix~\ref{app:proof-constant-lower} for a proof of this claim.

\section{Discussion}

In this work, we determined the minimax risk of estimation for
observation models of the form~\eqref{eqn:operator-model}, where one
observes the image of a unknown parameter under a random linear
operator with additive noise. Our results reveal the dependence of the
rate of convergence on the covariate law, the parameter space, the
error metric, and the noise level. We conclude our paper by presenting
some simulation results; see Section~\ref{sec:simulation-results}

Finally, we note that in this work we studied minimax risks of
convergence in expectation.  This is convenient, as it requires
relatively minor assumptions of the distribution of $\LinOp{\xi}$.  On
the other hand, for the setting of random design regression,
high-probability results, such as those obtained in the
papers~\cite{AudEtAl11, Men15, HsuSab16, LecMen16, Oli16}, typically
require stronger assumptions such as the sub-Gaussianity of the
covariate distribution.  Nonetheless, high-probability guarantees
provide a complementary perspective on the problem we
consider. Indeed, when the covariate law can be considered
``heavy-tailed,'' it may be more relevant to develop robust estimators
that have low risk with high probability.  We refer to the survey
article~\cite{LugMen19} for a overview of work in this direction.

\subsection{Some illustrative simulations}
\label{sec:simulation-results}

We conclude our paper by presenting the results of some simulations
reveal how changes in the distribution of the random operator
$\LinOp{\xi}$ can lead to dramatic changes in the overall minimax
risk.

In this section, we present simulation results to illustrate the
behavior of the functionals appearing in our main results for two
versions of random design linear regression. In
Section~\ref{sec:random-design-simulation}, we present simulation
results for a multivariate, random design linear regression setting
with \iid{} covariates. Concretely, we provide two different covariate
laws, where the minimax error for the same parameter space differs by
at least two orders of magnitude. We emphasize this difference in
\emph{entirely} due to the covariate law; the noise, observation
model, error metric, and parameter space are fixed in this comparison.

Additionally, in Section~\ref{sec:markov-simulation}, we present
simulation results for a univariate regression setting where the
covariates are sampled from a Markov chain.  In both cases, the
functional is able to capture the dependence of the minimax rate of
estimation on the underlying covariate distribution.

\subsubsection{Higher-order effects in \iid{} random design linear regression}
\label{sec:random-design-simulation}

For random design linear regression, higher order properties of the
covariate distribution over the covariates can have striking effects
on the minimax risk.  In order to illustrate this phenomenon, we
consider the regression model~\eqref{eqn:iid-regression} with feature
map $\psi(x) = x$, and parameter vector $\thetastar$ constrained to a
ball in the Euclidean norm.  We then construct a family of
distributions over the covariates that are all zero-mean with identity
covariance, but differ in interesting ways in terms of their
higher-order moment properties.  More precisely, we let $\delta_0$
denote the Dirac measure with unit mass at $0$, and for a mixture
weight $\lambda \in [0,1]$, we consider covariates generated from the
probability distribution
\begin{equation}
\label{eqn:covariate-ensemble}
P_{\lambda} \defn \lambda \delta_0 + (1-\lambda)
\Normal{0}{\frac{1}{1-\lambda} I_\dimension}.
\end{equation}
By construction, all members of the ensemble have the same behavior
with respect to their first and second moments,
\begin{equation}
\label{eqn:sim-moment-restrictions}
\E_{P_\lambda}[x] = 0 \quad \mbox{and} \quad \Cov_{P_\lambda}(x) =
\E_{P_\lambda}[x\otimes x] = I_\dimension, \quad \mbox{for
  all}~\lambda \in [0, 1].
\end{equation}
In the special case $\lambda = 0$, the distribution $P_\lambda$
corresponds to the standard Gaussian law on $\R^\dimension$, whereas
it becomes an increasingly ill-behaved Gaussian mixture distribution
as $\lambda \rightarrow 1^{-}$.

Following the argument in
Section~\ref{sec:gaussian-linear-regression}, in this case, the
minimax risk is upper and lower bounded as
\begin{equation}
\label{eqn:gen-expression-simulation-minimax}
\frac{\sigma^2}{\numobs} \E_{P_\lambda^\numobs}[\trace((\EmpCov +
  \tfrac{c_d \sigma^2 \dimension}{\numobs \radius^2} I_\dimension
  )^{-1})] \\
\leq \mathfrak{M}^{\rm \iid}_\numobs\Big(P_\lambda, \radius, \sigma^2,
I_\dimension, I_\dimension\Big) \leq \frac{\sigma^2}{\numobs}
\E_{P_\lambda^\numobs}[\trace((\EmpCov + \tfrac{\sigma^2
    \dimension}{\numobs \radius^2} I_\dimension )^{-1})].
\end{equation}
Above, the lower bound constant $c_d$ is defined in
display~\eqref{eqn:sharp-lower-bound-constant-dicker}.

To understand the effect of the covariate law, we fix the
signal-to-noise ratio such that $\tfrac{\radius}{\sigma} = \tau$, for
$\tau \in \{1, 10\}$.  Note that after renormalizing the minimax risk
by $\radius^2$, it only depends on $\tau$ (and not on the particular
choices of $(\radius, \sigma)$).  Similarly, this invariance relation
holds for the functionals appearing on the left- and righthand sides
of the display~\eqref{eqn:gen-expression-simulation-minimax}---after
normalization by $1/\radius^2$, they no longer depend on $(\radius,
\sigma)$ except via the ratio $\tau = \tfrac{\radius}{\sigma}$.
Additionally, we fix the aspect ratio $\gamma =
\tfrac{\dimension}{\numobs}$.\footnote{ Specifically, we take $d = \ceil{\gamma
  n}$.}%
By varying $\gamma \in [0.05, 4]$ we are able to illustrate the
behavior of the minimax risk, as characterized by our functional, for
problems which are both under- and overdetermined.

Having fixed the SNR at $\tau$ and aspect ratio at $\gamma$, we can
somewhat simplify the
display~\eqref{eqn:gen-expression-simulation-minimax}, by introducing
the following quantities which only depend on the parameters $\tau,
\gamma$ and the sample size $\numobs$ and the mixture parameter
$\lambda$,
\begin{subequations}
\label{eqn:rescaled-bounds}
\begin{align}
\mathfrak{m}_{\numobs}(\lambda, \tau, \gamma) &\defn
\frac{\mathfrak{M}^{\rm \iid}_\numobs\Big(P_\lambda, \tau \sigma ,
  \sigma^2, I_{\ceil{\gamma n}}, I_{\ceil{\gamma n}}\Big)}{\tau^2
  \sigma^2}, \\ u_\numobs(\lambda, \tau, \gamma) &\defn
\frac{1}{\tau^2 \numobs} \E_{P_\lambda^\numobs}[\trace((\EmpCov +
  \tfrac{ \ceil{\gamma \numobs}}{\numobs \tau^2} I_{\ceil{\gamma
      \numobs}} )^{-1})], \\ \ell_\numobs(\lambda, \tau, \gamma)
&\defn \frac{1}{\tau^2 \numobs} \E_{P_\lambda^\numobs}[\trace((\EmpCov
  + \tfrac{c_\dimension \ceil{\gamma \numobs}}{\numobs \tau^2}
  I_{\ceil{\gamma \numobs}} )^{-1})].
\end{align}
\end{subequations}
Then, the relations~\eqref{eqn:gen-expression-simulation-minimax}, can
be equivalently expressed as
\begin{align*}
\ell_\numobs(\lambda, \tau, \gamma) \leq
\mathfrak{m}_{\numobs}(\lambda, \tau, \gamma) \leq u_\numobs(\lambda,
\tau, \gamma),
\end{align*}
and moreover this holds for all $\lambda \in [0, 1], \tau > 0, \gamma
> 0$.  In our simulation, we use Monte Carlo simulation with 50 trials
to estimate the upper and lower bound functionals $\ell_\numobs$ and
$u_\numobs$.

\begin{figure}
\centering
\begin{subfigure}{0.48\textwidth}
\includegraphics[width=\textwidth]{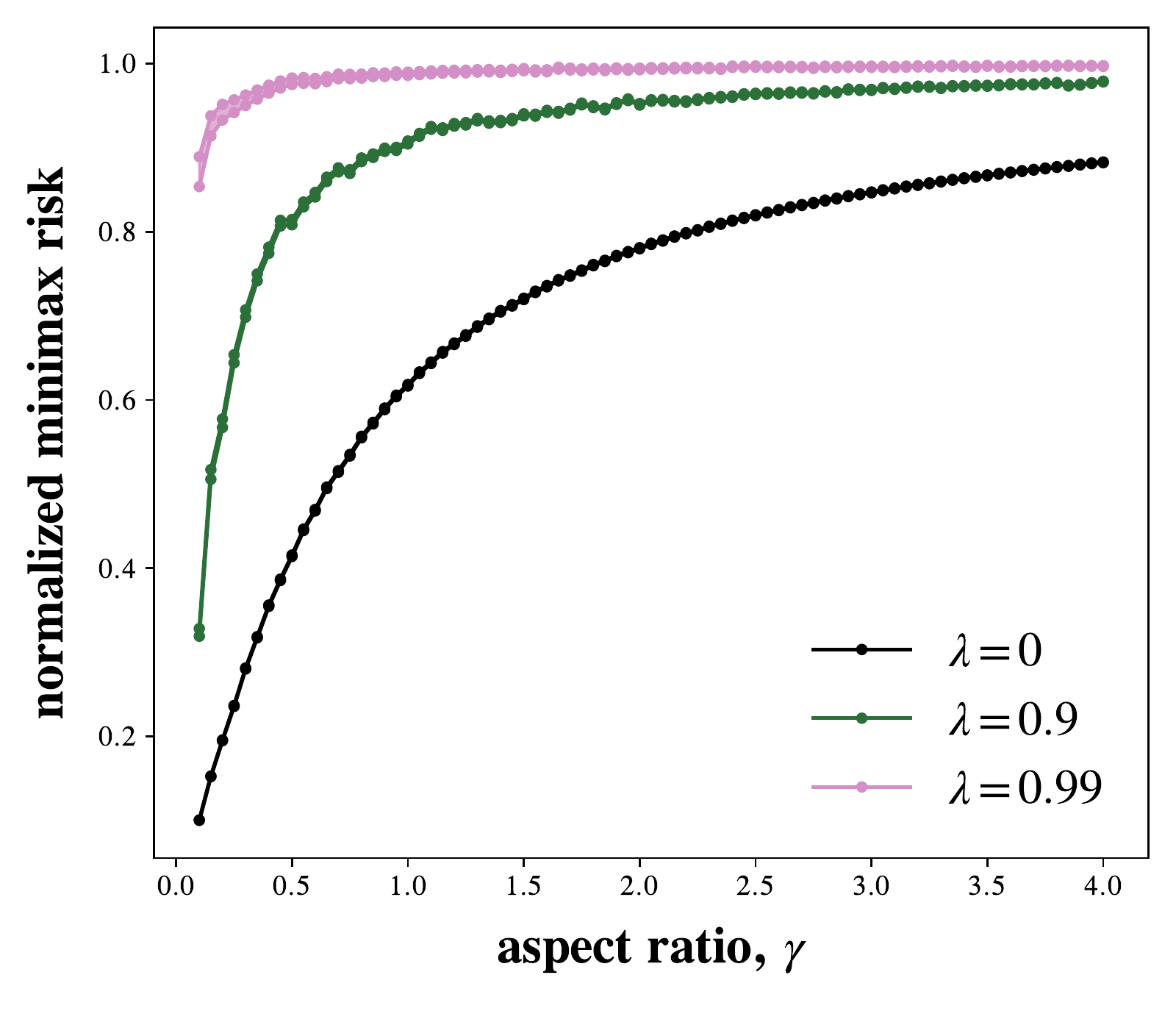}
\caption{$n = 128,~\tau = 1$}
\vspace*{1em}
\label{fig:n-128-SNR-1}
\end{subfigure}
\hfill
\begin{subfigure}{0.48\textwidth}
\includegraphics[width=\textwidth]{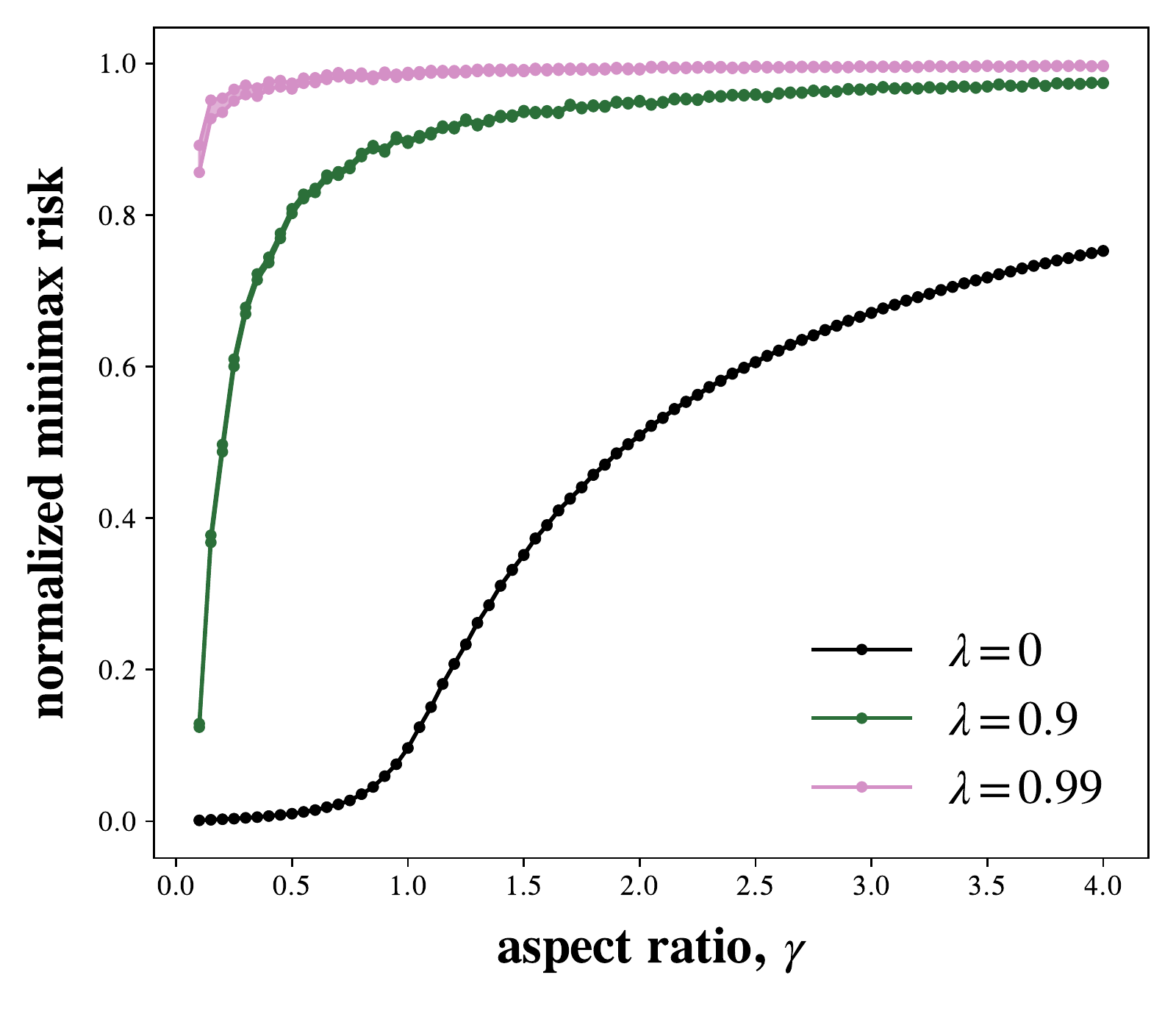}
\caption{$n = 128,~\tau= 10$}
\vspace*{1em}
\label{fig:n-128-SNR-100}
\end{subfigure}
\begin{subfigure}{0.48\textwidth}
\includegraphics[width=\textwidth]{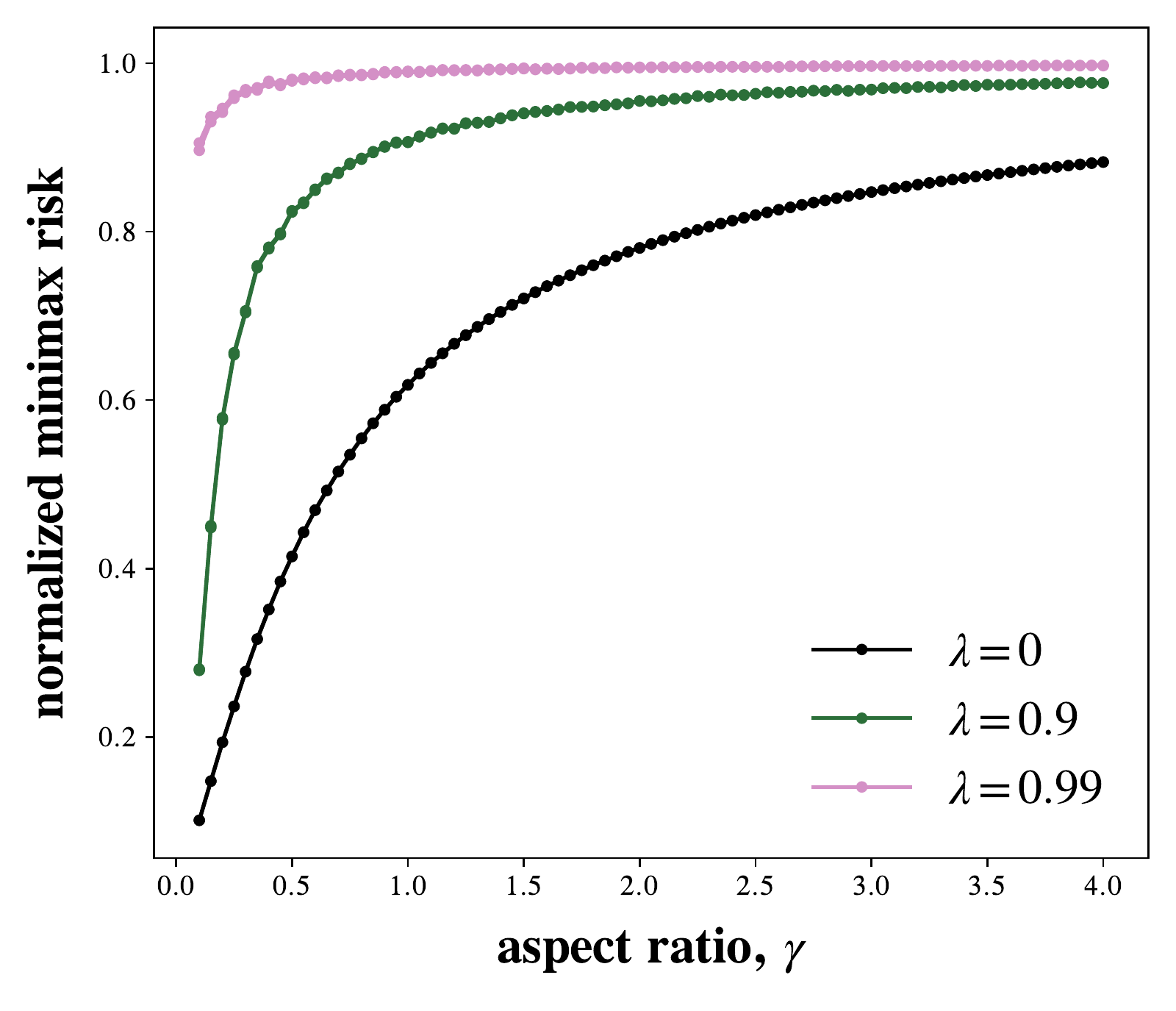}
\caption{$n = 512,~\tau = 1$}
\label{fig:n-512-SNR-1}
\end{subfigure}
\hfill
\begin{subfigure}{0.48\textwidth}
\includegraphics[width=\textwidth]{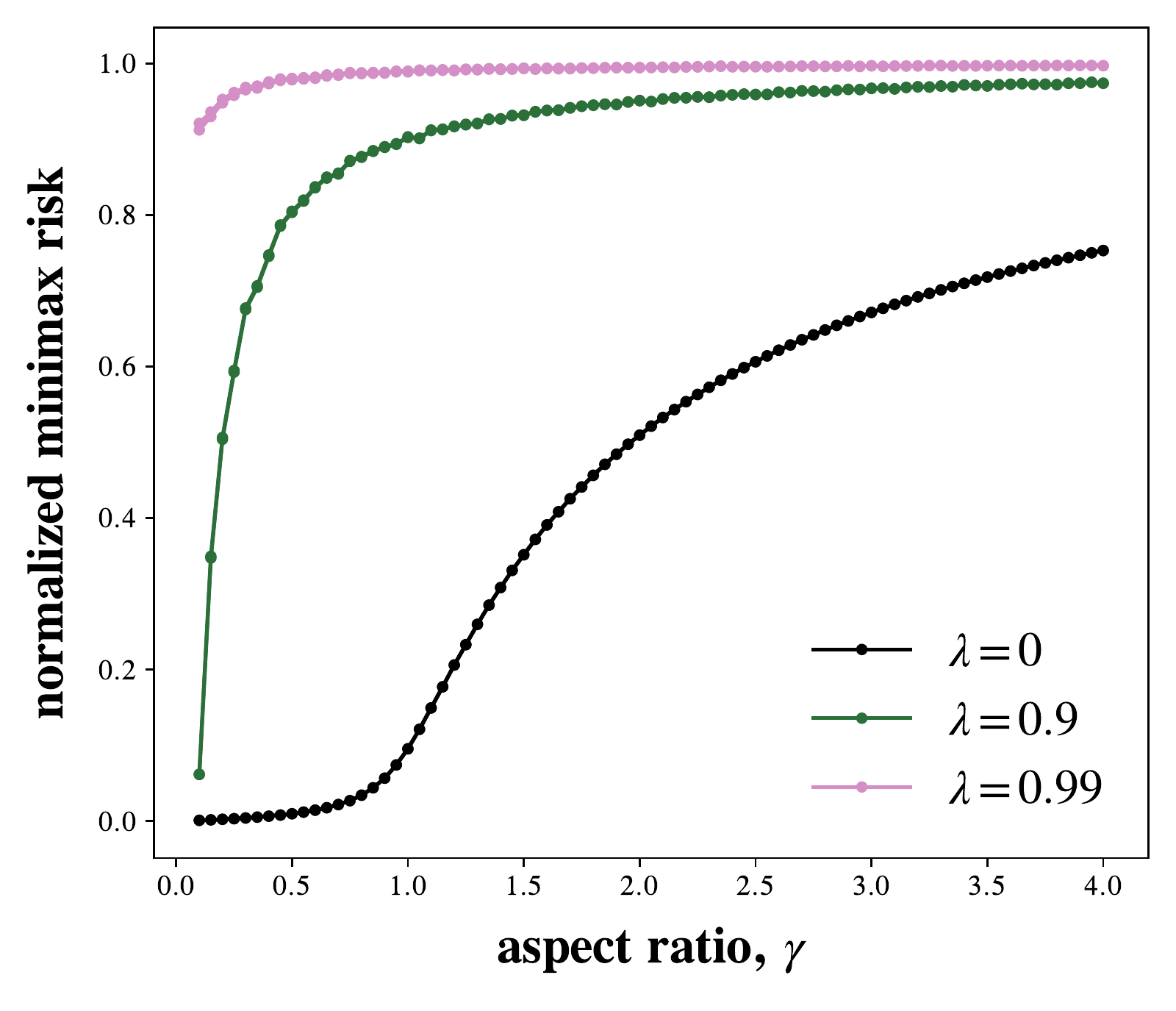}
\caption{$n = 512,~\tau = 10$}
\label{fig:n-512-SNR-100}
\end{subfigure}
\caption{ Simulations of random design regression for three covariate
  laws, $P_\lambda$ as defined in
  equation~\eqref{eqn:covariate-ensemble} with $\lambda \in \{0, 0.9,
  0.99\}$.  For a given choice of the mixture weight $\lambda$ and
  signal-to-noise ratio (SNR) $\tau$, we plot the lower bound
  $\ell_\numobs(\lambda, \tau, \gamma)$ and upper bound
  $u_\numobs(\lambda, \tau, \gamma)$ as $\gamma$ varies between $0.05$
  and $4$.  The normalized minimax risk $\mathfrak{m}_\numobs$ is then
  guaranteed to lie in the region whose upper and lower envelopes are
  given by $u_\numobs$ and $\ell_\numobs$, respectively. To facilitate
  interpretation of these figures, we have shaded this region to
  highlight where we can guarantee the minimax risk
  $\mathfrak{m}_\numobs$ must lie.  The quantities $u_\numobs,
  \ell_\numobs, \mathfrak{m}_\numobs$ are all defined in
  display~\eqref{eqn:rescaled-bounds}.  In
  panels~\eqref{fig:n-128-SNR-1} and~\eqref{fig:n-128-SNR-100}, we set
  the sample size $\numobs = 128$, and set the SNR as $\tau = 1, 10$,
  respectively.  In panels~\eqref{fig:n-512-SNR-1}
  and~\eqref{fig:n-512-SNR-100}, we set the sample size $\numobs =
  512$, and set the SNR as $\tau = 1, 10$, respectively.  The plots
  above demonstrate that as $\lambda$ increases, the minimax risks are
  much worse.  Numerically, in the setting where $n = 512$ and $\tau =
  10$---as depicted in panel~\eqref{fig:n-512-SNR-100}---our upper and
  lower bounds guarantee that the minimax risk for the isotropic
  ensemble (depicted with $\lambda = 0$ above) can be over 806 times
  larger than the minimax risk for the ensemble with $\lambda = 0.99$.
  It should be noted that in this comparison the first and second
  moments of the ensemble are held fixed (see
  equation~\eqref{eqn:sim-moment-restrictions}), and hence the
  differences between the lines plotted in any given panel can only be
  explained by differences in higher-order moments within the ensemble
  $\{P_\lambda\}$.  The figures also demonstrate that the gap between
  our upper and lower bounds is fairly small, particularly whenever $d
  > 5$.}
\label{fig:random-reg-sim}
\end{figure}

In our simulations, we take $\lambda \in \{0, 0.9, 0.99\}$ and vary
$\gamma \in [0.05, 4]$.  The results of these simulations are
presented in Figure~\ref{fig:random-reg-sim}; see the caption for a
detailed description and commentary.  The general pattern should be
clear: the covariate law can have a dramatic impact on the overall
rate of estimation, even when restricting some moments such as we have
with the relations~\eqref{eqn:sim-moment-restrictions}.

\subsubsection{Mixing time effects in Markovian linear regression}
\label{sec:markov-simulation}

Covariates need not be drawn in an \IID{} manner, and any dependencies
can be expected to affect the minimax risk.  Here we illustrate this
general phenomena via some simulations for the Markov regression
example as outlined in Section~\ref{sec:markov-regression}.  We seek
to study a wide range of possible mixing conditions for the Markovian
covariate model.  In order to do so, we consider covariates generated
from the Markovian model~\eqref{eqn:markov-covariates} with
\[
r_t = \frac{\psi(t-1)}{\psi(t)},
\]
where $\psi \colon \N \cup \{0\} \to \R_+$ is a nondecreasing function
satisfying $\psi(0) = 1$ and $\lim_{t \to \infty} \psi(t) = \infty$.
With this choice, it is easily checked that, marginally
\[
x_t \sim \Normal{0}{1 - \frac{1}{\psi(t)}}.
\]
Therefore, $x_t\to \Normal{0}{1}$ in distribution as $t \to \infty$, and the
rate of convergence is of order $1/\psi(t)$. 

We now illustrate how the minimax rate, as determined in
Corollary~\ref{cor:markov}, for this problem behaves for different
choices of the function $\psi$ and the signal-to-noise ratio (SNR). As
in Section~\ref{sec:random-design-simulation}, we normalize the
minimax risk by the squared radius so that it only depends on $\tau =
\tfrac{\radius}{\sigma}$. The quantity we then plot is
\[
\Phi_T(\tau) \defn \frac{\Phi_T(\tau, 1)}{\tau^2},
\]
where $\Phi_T(\radius, \sigma)$ is the functional appearing in Corollary~\ref{cor:markov}.

In the simulation, we consider the following choices of scaling function $\psi$, 
\[
5^t, \quad t + 1, \quad 1 + \log(t + 1), \quad \mbox{and} \quad 
1 + \log\big(1 + \log(t + 1)\big). 
\]
With the choice $\psi(t) = 5^t$, the underlying Markov chain converges geometrically to the
standard Normal law. On the other hand, the choice $\psi(t) = \log(1 + \log(1 + t)) + 1$ exhibits
much slower convergence---the variational distance between the law of $x_t$ and 
$\Normal{0}{1}$ is of order $O(1/(\log \log t))$.

\begin{figure}
\centering 
\begin{subfigure}{0.48\textwidth}
\includegraphics[width=\textwidth]{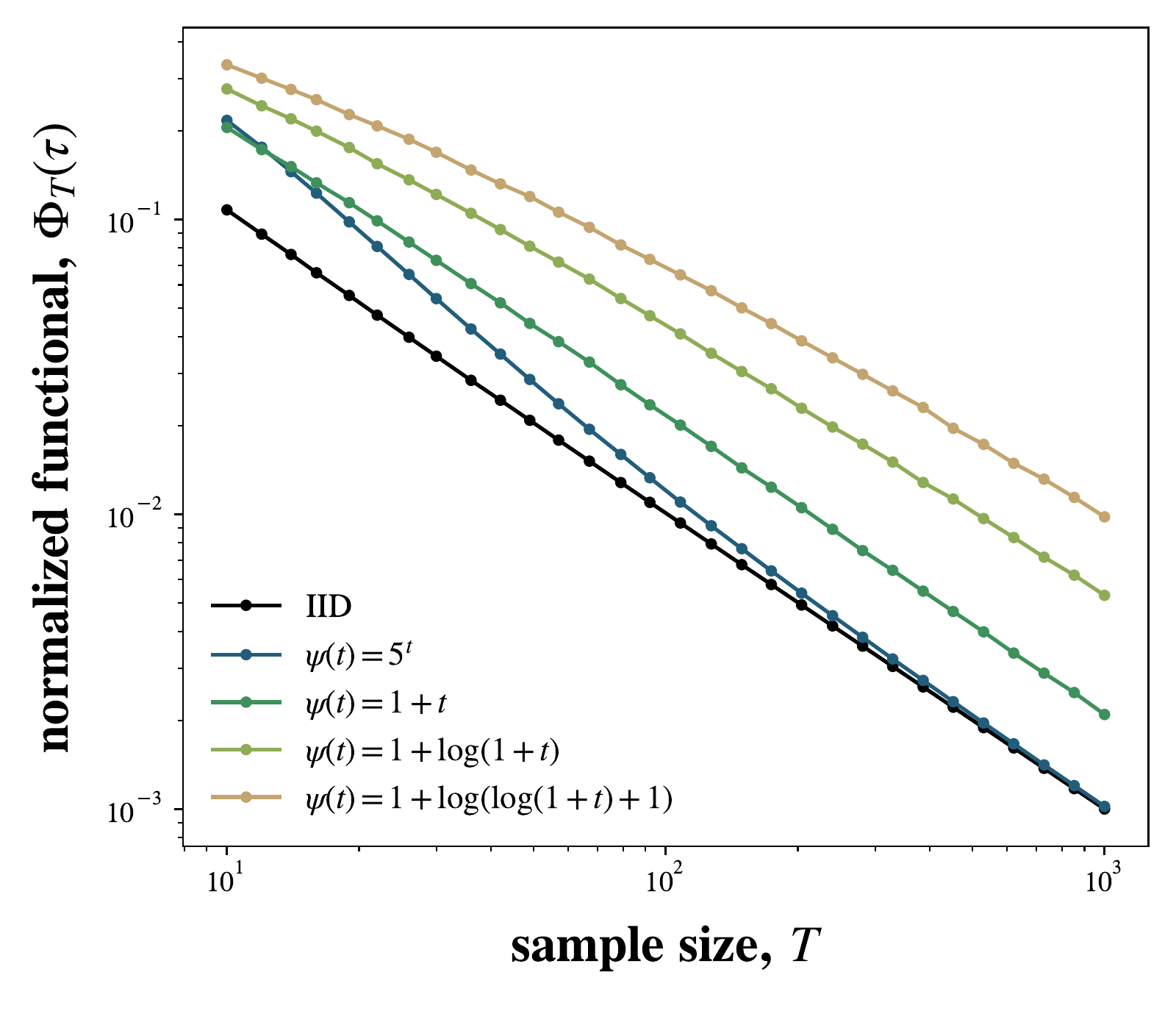}
\caption{$\tau = 1$}
\vspace*{1em}
\label{fig:markov-SNR-1}
\end{subfigure}
\hfill
\begin{subfigure}{0.48\textwidth}
\includegraphics[width=\textwidth]{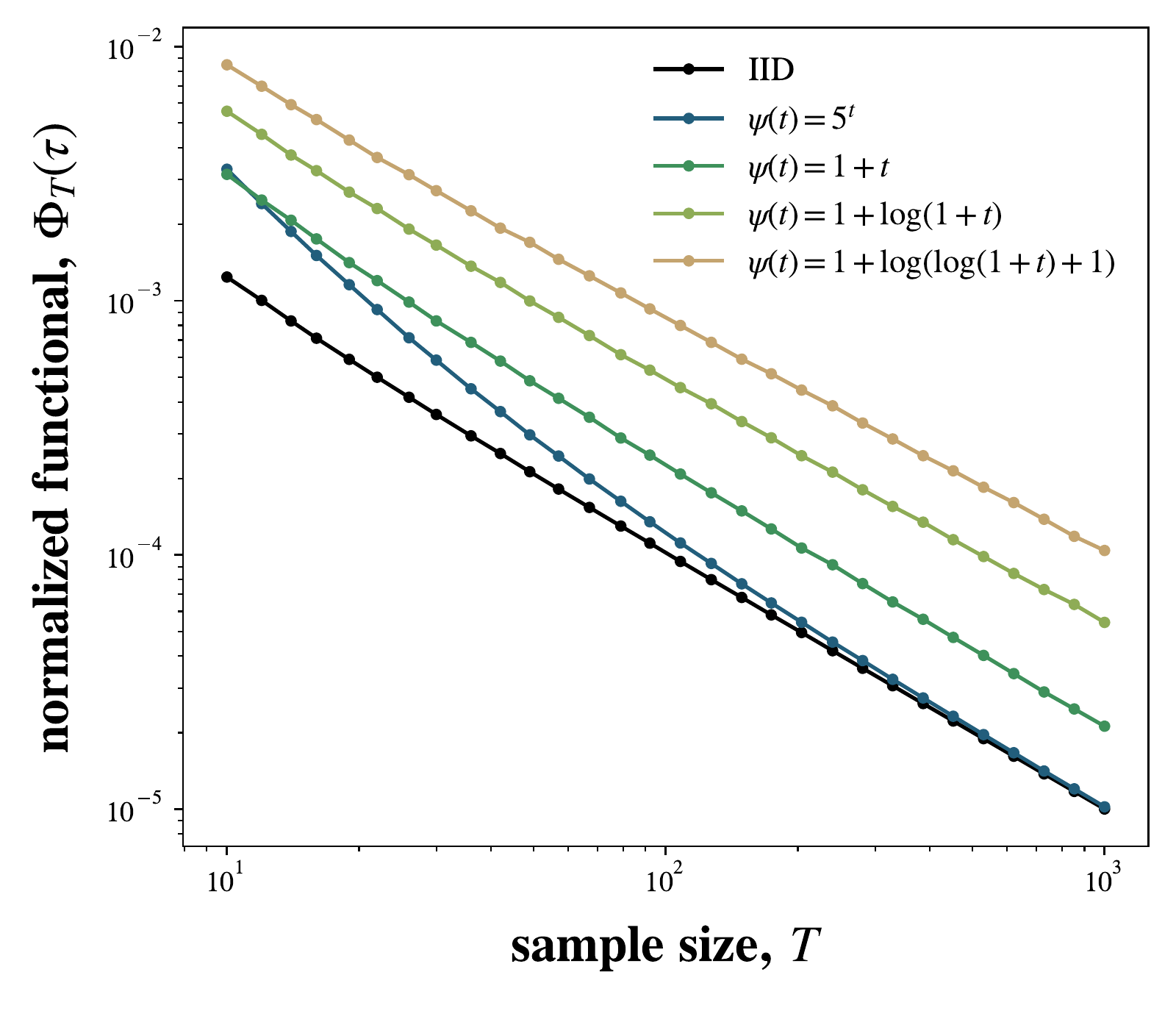}
\caption{$\tau= 10$}
\vspace*{1em}
\label{fig:markov-SNR-10}
\end{subfigure}
\caption{
  Simulations for five distributions of Markovian covariates.
  In panel~\eqref{fig:markov-SNR-1}, we set the SNR parameter as $\tau = 1$, and
  in panel~\eqref{fig:markov-SNR-10}, we set the SNR parameter as $\tau = 10$. 
  As the scaling function $\psi$ grows more slowly, the chain converges to its stationary distribution more
  slowly, and the minimax rate decays more slowly, as indicated by the displayed behavior of our
  functional $T \mapsto \Phi_T(\tau)$.}
\label{fig:markov-sim}
\end{figure}

We simulate each of these chains, computing the normalized functional
$\Phi_T(\tau)$ over the course of 5000 Monte Carlo trials. The
sample size $T$ is varied between 10 and $3162$.  In the simulation we
also include the choice $r_t \equiv 0$, which corresponds to \iid{}
covariates.  The results of the simulation are presented in
Figure~\ref{fig:markov-sim}; see the caption for more details and
commentary.

\subsection*{Acknowledgements}
We thank Jaouad Mourtada for a helpful conversation and useful email
exchanges; we also thank Peter Bickel for a helpful discussion
regarding his prior work on the Gaussian sequence model.  RP was
partially supported by a UC Berkeley Chancellor's Fellowship via the
ARCS Foundation.  MJW and RP were partially funded by ONR grant
N00014-21-1-2842 and National Science Foundation grant NSF-DMS grant
2015454. MJW and RP gratefully acknowledge funding support from
Meta via the UC Berkeley AI Research (BAIR) Commons initiative. 

\appendix

\section{Proofs from Section~\ref{sec:main-result}}

\subsection{Proof of Proposition~\ref{prop:concavity-functional}}\label{app:proof-concavity-functional}

The constraint set is evidently convex, as it is formed by the intersection of
of two convex sets: the $\dimension \times \dimension$ real, symmetric positive definite matrices with the 
hyperplane $\{\Omega : \trace(\ConstraintMat^{-1} \Omega) \leq \radius^2\}$. 

We claim that the objective function $f$ is concave over the set of symmetric positive definite matrices. 
It can be expressed as
\begin{align*}
f(\Omega) =  \E_\xi[g(\LinOp{\xi}^\T \NoiseCovariance^{-1} \LinOp{\xi}, \Omega)], \quad
\mbox{where} \quad g(X, \Omega) \defn
\trace(\EstimationMat^{\!1/2} (X + \Omega^{-1})^{-1}
\EstimationMat^{\!1/2}).
  \end{align*}
Evidently to establish that $f$ is concave, it is enough to show that $g(X, \cdot)$ is concave
for every symmetric positive semidefinite $X$. 
In order to establish this claim, let us fix some $\eps > 0$, and define
\mbox{$X(\eps) \defn X + \eps
  I_\dimension$.}  By the joint concavity of the harmonic mean
of positive operators~\cite[Corollary~37.2]{Sim19}, it follows that
for any pair of positive definite matrices $\Omega, \Omega'$, we have
\begin{align*}
\Big(X(\eps) + \Big(\frac{\Omega + \Omega'}{2}\Big)^{-1}\Big)^{-1}
\succeq \frac{1}{2} \Big(X(\eps) + \Omega^{-1}\Big)^{-1} +
\frac{1}{2} \Big(X(\eps) + (\Omega')^{-1}\Big)^{-1}.
\end{align*}
Passing to the limit as $\eps \to 0$ yields
\begin{align*}
\Big(X + \Big(\frac{\Omega + \Omega'}{2}\Big)^{-1}\Big)^{-1}
\succeq \frac{1}{2} \Big(X + \Omega^{-1}\Big)^{-1} +
\frac{1}{2} \Big(X + (\Omega')^{-1}\Big)^{-1}.
\end{align*}
Since the trace is a monotone mapping on positive definite matrices,
and $g$ is continuous in its second argument, 
we obtain the claimed concavity of $g$.

\subsection{Proof of Proposition~\ref{prop:to-population}}\label{app:proof-to-population}

To establish the upper bound, it suffices to show that for each positive
definite $\Omega \succ 0$ with $\trace(\ConstraintMat^{-1/2} \Omega \ConstraintMat^{-1/2}) \leq
\frac{\numobs \radius^2}{\sigma^2}$ that the following inequality holds
\begin{equation}
\label{EqnAuxiliary}
\trace\Big( \E\big[ \big(\EmpCov + \Omega^{-1}
  \big)^{-1} \Sigma_P \big]\Big) \leq \Big(1 +
\frac{\radius^2 \kappa^2}{\sigma^2}\Big) \trace\big(
\big(\Sigma_P + \Omega^{-1} \big)^{-1}
\Sigma_P\Big).
\end{equation}

Our proof of the auxiliary claim~\eqref{EqnAuxiliary} is based on
exchangeability and operator convexity, and is similar to previous work on
the analysis of ridge regression estimators~\cite{MouRos22}.  
Let $x_{\numobs + 1}$ be a fresh sample drawn
independently from $\{x_i\}_{i = 1}^\numobs$ with the same distribution $P$.  
Letting $\E$ denote the expectation over the full sequence $\{x_i\}_{i=1}^{\numobs + 1}$, 
we have
\begin{equation}
\label{eqn:relation-one-more}
\E\big[\big(\EmpCov + \Omega^{-1}
  \big)^{-1} \Sigma_P \big] = \numobs
\E\big[ \big(\numobs \EmpCov + \numobs \Omega^{-1}
  \big)^{-1} (\psi(x_{\numobs + 1}) \otimes \psi(x_{\numobs +
    1}))\big].
  \end{equation}
Define $\hat \Sigma_{\numobs + 1} \defn (\numobs + 1)^{-1}
\sum_{i=1}^{\numobs + 1} \psi(x_i) \otimes \psi(x_i)$.  Then, by the
Sherman--Morrison lemma~\cite[Section 0.7.4]{HorJoh13}, it follows
that
\begin{align*}
(\numobs \EmpCov + \numobs \Omega^{-1})^{-1} \psi(x_{\numobs + 1}) = (1 +
  \langle (\numobs \EmpCov + \numobs \Omega^{-1})^{-1} \psi(x_{\numobs +
      1}), \psi(x_{\numobs + 1})\rangle \big((\numobs + 1) \hat \Sigma_{\numobs +
    1} + \numobs \Omega^{-1}\big)^{-1} \psi(x_{\numobs + 1}).
  \end{align*}
Additionally, by the Cauchy--Schwarz inequality, we have
\begin{align*}
\langle (\numobs \EmpCov + \numobs \Omega^{-1})^{-1} \psi(x_{\numobs + 1}), \psi(x_{\numobs
  + 1}) \rangle \leq \frac{1}{\numobs} \opnorm{\ConstraintMat^{\!-1/2} \Omega} \|\ConstraintMat^{1/2} \psi(x_{\numobs + 1})\|_2^2
\leq \frac{\radius^2 \kappa^2}{\sigma^2},
\end{align*}
where the last inequality holds $P$-almost surely.
Applying the previous two displays in
relation~\eqref{eqn:relation-one-more}, it follows that
  \begin{align}
\trace \E\big[ \big(\EmpCov + \Omega^{-1} \big)^{-1}
  \Sigma_P \big] &\leq 
\Big(1 + \frac{\radius^2 \kappa^2}{\sigma^2}\Big)
\trace\E\big[ \big(\hat \Sigma_{\numobs + 1} +
  \Omega^{-1}\big)^{-1} \psi(x_{\numobs + 1}) \otimes \psi(x_{\numobs + 1})
\big] \nonumber  \\
& = \Big(1 + \frac{\radius^2
  \kappa^2}{\sigma^2}\Big) \frac{1}{\numobs + 1} \sum_{i=1}^{\numobs +
  1}\trace\E\big[ \big(\hat \Sigma_{\numobs + 1} +
  \Omega^{-1}\big)^{-1} \psi(x_{i}) \otimes \psi(x_{i})  \big] \\ &= \Big(1
+ \frac{\radius^2 \kappa^2}{\sigma^2}\Big) \trace\E\big[
  \big(\hat \Sigma_{\numobs + 1} + \Omega^{-1}\big)^{-1} \hat
  \Sigma_{\numobs + 1} \big] \label{eqn:pop-step-1}\\
& \leq \Big(1 + \frac{\radius^2
  \kappa^2}{\sigma^2}\Big) \trace\Big(
  \big(\Sigma_P + \Omega^{-1}\big)^{-1}\Sigma_P \Big). \label{eqn:pop-step-2}
  \end{align}
Above step~\eqref{eqn:pop-step-1} follows by the exchangeability of
$\{\psi(x_{i})\}_{i=1}^{\numobs + 1}$, and step~\eqref{eqn:pop-step-2} follows by the cyclicity and
linearity of the trace, as well as the the fact that for any fixed
symmetric positive definite matrix $B$, the mapping $A \mapsto (A +
B)^{-1}A = I_\dimension - (A+B)^{-1}B$ is concave over the set
over symmetric positive semidefinite matrices (see Bhatia~\cite[page~19]{Bha07}).

\subsection{Proof of Corollary~\ref{cor:limit-relation}}\label{app:proof-limit-relation}

Combining Theorems~\ref{thm:main-lower} and~\ref{thm:main-lower}, we find that 
\begin{equation}\label{eqn:sandwich-minimax}
\Phi(T, \bP, \NoiseCovariance, \tfrac{\radius}{2}, \EstimationMat, \ConstraintMat)  
\leq 
\mathfrak{M}(T, \bP, \NoiseCovariance, \radius, \EstimationMat, \ConstraintMat) 
\leq
\Phi(T, \bP, \NoiseCovariance, \radius, \EstimationMat, \ConstraintMat).  
\end{equation}
Evidently, by definition of the functional $\Phi$ (see definition~\eqref{eqn:defn-functional-general}), the map $\radius \to \Phi(T, \bP, \NoiseCovariance, \radius, \EstimationMat, \ConstraintMat)$ is nondecreasing. Moreover since 
$\LinOp{\xi}^\T \NoiseCovariance^{-1} \LinOp{\xi}$ is invertible with probability 1, 
it is a bounded function. Therefore, 
\begin{align*}
\lim_{\radius \to \infty} \frac{\Phi(T, \bP, \NoiseCovariance, \radius, \EstimationMat, \ConstraintMat)}{\Phi(T, \bP, \NoiseCovariance, \radius/2, \EstimationMat, \ConstraintMat)}
=1,
\end{align*}
which in view of the sandwich relation~\eqref{eqn:sandwich-minimax}, furnishes the claim.

\section{Proofs and calculations from Section~\ref{sec:examples}}
\label{app:examples}

\subsection{Proof and calculations from Section~\ref{sec:examples-parametric}}

\subsubsection{Proof of equation~\eqref{eqn:gaussian-regression-functional}}
\label{app:proof-gauss-calc}

From the definition of the functional~\eqref{eqn:d-functional-iid}, we
have
\begin{align*}
d_\numobs(\Normal{0}{I_\dimension}, \radius, \sigma^2, I_\dimension, I_\dimension) 
= \sup \Big\{\,
\E[\trace((\EmpCov + \tfrac{\sigma^2
    \dimension}{\numobs \radius^2} M^{-1} )^{-1})] : M \succ
0,~~\trace(M) = d\,\Big\}.
\end{align*}
In this section, all expectations are over $x_i \simiid
\Normal{0}{I_\dimension}$.  
We claim that the supremum above is achieved at $M = I_\dimension$.

\begin{lemma}
\label{lem:attainment}
For any positive definite matrix $M \succ 0$ such that $\trace(M) =
d$, we have
\begin{align*}
\E[\trace((\EmpCov + \tfrac{\sigma^2 \dimension}{\numobs \radius^2}
  M^{-1} )^{-1})] \leq \E[\trace((X^\T X + \tfrac{\dimension
    \sigma^2}{\radius^2} I_\dimension )^{-1})]
\end{align*}
\end{lemma}
\noindent Assuming Lemma~\ref{lem:attainment}, we then have 
\begin{align*}
d_\numobs(\Normal{0}{I_\dimension}, \radius, \sigma^2, I_\dimension, I_\dimension) 
= 
\E[\trace((\EmpCov + \tfrac{\sigma^2
    \dimension}{\numobs \radius^2} I_\dimension )^{-1})]
= 
d_{\rm Dicker}(\numobs, \dimension, \radius, \sigma),    
\end{align*}
which establishes~\eqref{eqn:gaussian-regression-functional}, as needed. 

\myparagraph{Proof of Lemma~\ref{lem:attainment}} 
Define the function $\phi \colon (\Sigma, M) 
\mapsto (\Sigma + \tfrac{\dimension \sigma^2}{\numobs
  \radius^2} M^{-1})^{-1}$, 
  where $\Sigma, M$ are assumed symmetric positive semidefinite and 
  $M$ is nonsingular.  
For each $\Sigma \succeq 0$, it is well known that
$\phi(\Sigma, \cdot)$ is operator concave~\cite[Corollary~37.2]{Sim19}---for any 
collection $\{M_i\}_{i=1}^d$ of
symmetric positive definite matrices, one has
\begin{equation}\label{eqn:concavity}
\frac{1}{\dimension} \sum_{i=1}^\dimension  \phi(\Sigma, M_i)
\preceq \phi(\Sigma, \tfrac{1}{\dimension} \sum_{i=1}^\dimension M_i),
\qquad\mbox{for any}~\Sigma \in \bS_+^\dimension.
\end{equation}
Now let $M \succ 0$ satisfying $\trace(M) = \dimension$ be given.  
Diagonalize $M$ so that $M = U \Lambda U^\T$, where $\Lambda =
\diag(\lambda) \succ 0$, and $U$ is orthogonal. 
Consider the cyclic permutations of $\Lambda$, given by
\begin{align*}
\Lambda^{(j)} = \diag(\lambda^{(j)}), \quad \mbox{where} \quad
\lambda^{(j)}_i = \lambda_{i + j}.
\end{align*}
Above, the arithmetic $i + j$ occurs modulo $d$.  
By rotational invariance of the Gaussian and the fact that $x_i$ has iid
coordinates, we have
\begin{align*}
  \E\trace((\EmpCov + \tfrac{d \sigma^2}{\numobs \radius^2}
  M^{-1})^{-1}) & = \E \trace((\EmpCov + \tfrac{d \sigma^2}{\numobs
    \radius^2} \Lambda^{-1})^{-1}) \\
& = \E \Big[\frac{1}{d} \sum_{j=1}^d \trace(( \EmpCov + \tfrac{d
      \sigma^2}{\numobs \radius^2} (\Lambda^{(j)})^{-1})^{-1})\Big]
  \\
& = \trace \Big\{\E \Big[\frac{1}{d} \sum_{j=1}^d \phi(\EmpCov,
    \Lambda^{(j)}) \Big]\Big\} \\ &\leq \trace \Big\{\E \Big[
    \phi(\EmpCov, \overline{\Lambda}) \Big]\Big\} \qquad \mbox{where}
  \quad \overline{\Lambda} \defn \tfrac{1}{d}
  \sum_{j=1}^d\Lambda^{(j)},
\end{align*}
The final inequality above uses the concavity
inequality~\eqref{eqn:concavity}, where we have taken $M_i = \Lambda^{(i)}$.  Now
note that
\begin{align*}
\overline{\Lambda} = \frac{\trace(\Lambda)}{\dimension}
I_\dimension = \frac{\trace(M)}{\dimension} I_\dimension =
I_\dimension.
\end{align*}
Combining the preceding displays furnishes the claim.

\subsubsection{Proof of the lower bound in equation~\eqref{eqn:our-dicker-result}}
\label{app:lower-bound-tail-bound-calc}

We apply our our sharp lower bound in Theorem~\ref{thm:main-lower} with 
$\Omega = \tfrac{\radius^2}{\dimension} I_\dimension$ and 
$\tau^2 = 1 - \tfrac{1}{2\dimension-1}$. Let us define 
$u = (1 - \tfrac{1}{2\dimension - 1})(1 - \P\{Z > 2d^2-d\})$, where $Z$ is a $\chi^2$-random variable 
with $d$-degrees of freedom. 
Note that $d(d-1) \geq \sqrt{dt} + t$ for $t = \tfrac{d^{3/2}}{4}$ for all $d \geq 2$. 
Therefore by standard tail bounds for $\chi^2$-variates~\cite[pp.~1325]{LauMas00}, we have $u \leq \exp(-d^{3/2}/4)$. 
Applying the sharp lower bound~\eqref{eqn:gen-sharp-lower} in 
Theorem~\ref{thm:main-lower} then yields the claim.

\subsubsection{Proof of equation~\eqref{eqn:mourtada-limit-relation}}
\label{app:mourtada-limit-proof}

Using the semidefinite inequality 
\begin{align*}
\big(\EmpCov + \Omega^{-1}\big)^{-1} \preceq \EmpCov^{-1}, 
\end{align*}
and the choice $\Omega = \tfrac{\numobs}{\sigma^2} \tfrac{\radius^2}{\dimension} I_\dimension$, 
we have the sandwich relation 
\begin{equation}\label{eqn:mourtada-limit-sandwich}
\trace \E_{P^\numobs}\big[ \Sigma_P^{1/2} (\EmpCov + \tfrac{\sigma^2}{\numobs} \tfrac{\dimension}{\radius^2} 
I_\dimension)^{-1}  \Sigma_P^{1/2} \big]
\leq 
d_\numobs(P, \radius, \sigma^2, I_\dimension, \Sigma_P) 
\leq 
\trace \E_{P^\numobs}\big[ \Sigma_P^{1/2} \EmpCov^{-1}  \Sigma_P^{1/2} \big],
\end{equation}
for all $\radius > 0$. Since $\radius \mapsto 
d_\numobs(P, \radius, \sigma^2, I_\dimension, \Sigma_P)$ is nondecreasing, 
the display above also demonstrates that this map has a limit. 
Now, note that by continuity, $P^\numobs$-almost surely we have 
\begin{align*}
\lim_{\radius \to \infty} \trace(\Sigma_P^{1/2} (\EmpCov + \tfrac{\sigma^2}{\numobs} \tfrac{\dimension}{\radius^2} 
I_\dimension)^{-1}  \Sigma_P^{1/2}) = 
\trace(\Sigma_P^{1/2} \EmpCov^{-1}  \Sigma_P^{1/2}).
\end{align*}
Thus, using the sandwich relation~\eqref{eqn:mourtada-limit-relation} and Fatou's lemma, 
we have 
\begin{multline*}
\trace \E_{P^\numobs}\big[ \Sigma_P^{1/2} \EmpCov^{-1}  \Sigma_P^{1/2} \big]
\leq 
\liminf_{\radius \to \infty} \trace \E_{P^\numobs}
\big[ \Sigma_P^{1/2} (\EmpCov + \tfrac{\sigma^2}{\numobs} \tfrac{\dimension}{\radius^2} 
I_\dimension)^{-1}  \Sigma_P^{1/2} \big] 
\\ \leq
\lim_{\radius \to \infty} d_\numobs(P, \radius, \sigma^2, I_\dimension, \Sigma_P) 
\leq 
\trace \E_{P^\numobs}\big[ \Sigma_P^{1/2} \EmpCov^{-1}  \Sigma_P^{1/2} \big],
\end{multline*}
which establishes relation~\eqref{eqn:mourtada-limit-relation}, as required. 

\subsubsection{Proof of minimax relation~\eqref{eqn:minimax-relation-markovian}}
\label{app:proof-minimax-relation-markovian}

Let us state the claim corresponding to relation~\eqref{eqn:minimax-relation-markovian}
somewhat more precisely. 
We define the functional 
\begin{align*}
\Phi_T(\radius, \sigma) \defn 
\E\bigg[\Big(\frac{1}{\radius^2} + \frac{z^\T M z}{\sigma^2}\Big)^{-1}\bigg]
\end{align*}
Then the following lemma corresponds to the claim underlying relation~\eqref{eqn:minimax-relation-markovian}. 
\begin{lemma}\label{lem:minimax-relation-markov}
The minimax risk under the Markovian observation model 
defined by the displays~\eqref{eqn:markov-covariates} and~\eqref{eqn:markov-responses}
satisfies
\begin{align*}
\frac{1}{4}\, \Phi_T(\radius, \sigma)
\leq
\inf_{\thetahat} \sup_{|\theta^\star|\leq \radius} 
\E \big[(\thetahat - \thetastar)^2 \big]
\leq
\Phi_T(\radius, \sigma).
\end{align*}
\end{lemma}
\noindent The remainder of this section is devoted to the proof of this claim

Note that if we define $\xi = (x_1,\dots, x_T)$, and 
$\LinOp{\xi} = x$, then the observation model~\eqref{eqn:markov-responses} can be written 
\begin{align*}
y = \LinOp{\xi}\thetastar + \NoiseCovariance^{1/2} w, 
\end{align*}
where $w \sim \Normal{0}{I_T}$ and $\NoiseCovariance = \sigma^2 I_T$. 
We have $\ConstraintMat = 1 = \EstimationMat$, since we are considering a univariate 
estimation problem. 
Therefore, since the functional~\eqref{eqn:defn-functional-general} is attained 
at $\Omega = \radius^2$, in order to establish Lemma~\ref{lem:minimax-relation-markov}, 
it is sufficient to show that 
\begin{equation}\label{eqn:sufficient-condition-lemma-markov}
\LinOp{\xi}^\T \NoiseCovariance^{-1} \LinOp{\xi} = \frac{x^\T x}{\sigma^2} = 
\frac{z^\T M z}{\sigma^2}.
\end{equation}
However, from display~\eqref{eqn:markov-covariates}, by induction we can establish 
that 
\begin{align*}
x_t = \sum_{s = 1}^t \sqrt{c_{st}} \, z_s, 
\end{align*}
where the coefficients $\{c_{st}\}$ are defined as in display~\eqref{eqn:gen-markov-matrix}.
Then, it follows that
\begin{align*}
  x^\T x = \sum_{t=1}^T \sum_{s, s' = 1}^t \sqrt{c_{st} c_{s' t}} z_{s} z_{s'}
  = \sum_{s, s' = 1}^T \underbrace{\sum_{t = \twomax{s}{s'}} \sqrt{c_{st} c_{s' t}}}_{= M_{ss'}} z_s z_{s'}.
\end{align*}
Using the display above, we establish the relation~\eqref{eqn:sufficient-condition-lemma-markov}, which in turn establishes Lemma~\ref{lem:minimax-relation-markov}, as needed. 

\subsection{Proof and calculations from Section~\ref{sec:examples-nonparametric}}

\subsubsection{Proof of limit relation~\eqref{eqn:limit-relation-gauss-seq-model}}
\label{app:proof-limit-relation-gauss-seq-model}

To lighten notation in this section, let us define the shorthands
\begin{subequations}
\label{eqn:both-shorthands-trunc-gauss-seq-model}
\begin{align}
\mathfrak{M}_k &\defn \mathfrak{M}_k\Big(\{\eps_j\}_{j=1}^k, \Theta_k(a, C)\Big), 
\quad \mbox{and,} \label{eqn:shorthand-k-trunc-gauss-seq-model}\\
\mathfrak{M} &\defn 
\mathfrak{M}\Big(\{\eps_j\}_{j=1}^\infty, \Theta(a, C)\Big) 
\defn \inf_{\thetahat} \sup_{\thetastar \in \Theta(a, C)} 
\E\Big[ \sum_{j=1}^\infty (\thetahat_j(y) - \thetastar_j)^2 \Big]. 
\label{eqn:shorthand-full-gauss-sequence-model}
\end{align}
\end{subequations}
We begin by stating the following sandwich relation for the minimax risks. 
\begin{lemma}\label{lem:sandwich-relation-gauss-seq-model-truncation}
The sequence of minimax risks $\{\mathfrak{M}_k\}$ and infinite-dimensional risk 
$\mathfrak{M}$ satisfies the sandwich relation 
\begin{equation}\label{eqn:lem-sandwich-relation-trunc-gauss-seq}
\mathfrak{M}_k \leq \mathfrak{M} \leq \mathfrak{M}_k + \frac{C^2}{a_{k+1}^2}, 
\end{equation}
for all $k \geq 1$.
\end{lemma} 
Assuming Lemma~\ref{lem:sandwich-relation-gauss-seq-model-truncation} 
for the moment, note that it implies for any divergent sequence $a_k \to \infty$ that
\begin{align*}
\lim_{k \to \infty} \mathfrak{M}_k = \mathfrak{M}.
\end{align*}
In view of the shorthands~\eqref{eqn:both-shorthands-trunc-gauss-seq-model}, the display above 
establishes our desired limit relation~\eqref{eqn:limit-relation-gauss-seq-model}. 

\myparagraph{Proof of Lemma~\ref{lem:sandwich-relation-gauss-seq-model-truncation}}
We begin by establishing the lower bound. Note that $\Theta_k(a, C) \subset \Theta(a, C)$, 
hence we have 
\begin{align*}
\mathfrak{M} &\geq \inf_{\thetahat} \sup_{\thetastar \in \Theta_k(a, C)} 
\E \Big[\sum_{j=1}^\infty(\thetahat_j((y_i)_{i=1}^\infty) - 
\thetastar_j)^2\Big] \\ 
&\geq \inf_{\thetahat} \sup_{\thetastar \in \Theta_k(a, C)} 
\E \Big[\sum_{j=1}^k(\thetahat_j((y_i)_{i=1}^\infty) - 
\thetastar_j)^2\Big],
\end{align*}
where the last equation arises since $\thetastar_j = 0$ for $j > k$
and thus any minimax optimal estimator over $\Theta_k(a, C)$ satisfies
$\thetahat_j \equiv 0$ for all $j > k$. The righthand side differs
from $\mathfrak{M}_k$ in that $\thetahat$ is a function of the full
sequence $y = (y_i)_{i=1}^\infty$. However, note that due to the
independence of the noise variables $z_i$, for the observation
model~\eqref{eqn:gauss-sequence-model} restricted to $\Theta_k(a, C)$,
the vector $y^{(k)} = (y_i)_{i=1}^k$ is a sufficient statistic.  Hence
we have for each $k \geq 1$,
\begin{align*}
\mathfrak{M} \geq 
\inf_{\thetahat} \sup_{\thetastar \in \Theta_k(a, C)} 
\E \Big[\sum_{j=1}^k(\thetahat_j(y^{(k)}) - 
\thetastar_j)^2\Big] = \mathfrak{M}_k,
\end{align*}
which establishes the lower bound in relation~\eqref{eqn:lem-sandwich-relation-trunc-gauss-seq}.

To establish the upper bound, note that we certainly may restrict 
the infimum in the definition of $\mathfrak{M}$ to those estimators 
taking values in $\R^k$ which only are a function of $y^{(k)}$. 
Indeed, we then find
\begin{align}
\mathfrak{M} &\leq 
\inf_{\thetahat \in \R^k} \sup_{\thetastar \in \Theta(a, C)} 
\E \Big[\sum_{j=1}^k(\thetahat_j(y^{(k)}) - 
\thetastar_j)^2 +  \sum_{j > k} (\thetastar_j)^2 \Big]
\label{ineq:minimax-upper-trunc-gauss-seq-model-1} \\ 
&\leq \mathfrak{M}_k 
+ \sup_{\thetastar \in \Theta(a, C)}
\sum_{j > k} (\thetastar_j)^2. \label{ineq:minimax-upper-trunc-gauss-seq-model}
\end{align}
The inequality~\eqref{ineq:minimax-upper-trunc-gauss-seq-model}
arises by taking the supremum over the two terms of the risk in 
display~\eqref{ineq:minimax-upper-trunc-gauss-seq-model-1}, and noting the 
first term only depends on the first $k$ coordinate of $\thetastar \in \Theta(a, C)$, 
and hence the supremum may be taken over $\Theta_k(a, C)$ in the first term so 
as to obtain $\mathfrak{M}_k$. 

Now observe by Hölder's inequality, and the membership $\thetastar \in \Theta(a, C)$,
\begin{align*}
\sum_{j > k} (\thetastar_j)^2 = 
\sum_{j > k} \frac{1}{a_j^2} (a_j^2 (\thetastar_j)^2)
\leq \Big(\max_{j > k} \frac{1}{a_j^2}\Big) C^2 = \frac{C^2}{a_{k + 1}^2}, 
\end{align*}
with the last equality arising because $j \mapsto a_j^2$ is assumed nondecreasing. 
Combining the display above with inequality~\eqref{ineq:minimax-upper-trunc-gauss-seq-model}
establishes the upper bound in~\eqref{eqn:lem-sandwich-relation-trunc-gauss-seq}, 
and thus establishes Lemma~\ref{lem:sandwich-relation-gauss-seq-model-truncation} 
as needed.

\subsubsection{Proof of relation~\eqref{eqn:sandwich-relation-gauss-sequence-model}}
\label{app:proof-sandwich-relation-gauss-seq-model}

Let us continue to adopt the shorthands $\mathfrak{M}_k$ and 
$\mathfrak{M}$ defined, respectively, in the 
displays~\eqref{eqn:shorthand-k-trunc-gauss-seq-model} 
and~\eqref{eqn:shorthand-full-gauss-sequence-model}. Moreover, 
we also use the shorthands 
\begin{align*}
R_k^\star \defn 
R_k^\star\Big(\eps, a, C\Big), 
\quad \mbox{and} \quad 
R^\star \defn R^\star(\eps, a, C),
\end{align*}
corresponding to the functionals~\eqref{eqn:functional-k-trunc-gauss-seq-model}
and~\eqref{eqn:functional-full-gauss-seq-model}, respectively.

We prove the following lemma. 
\begin{lemma} \label{lem:proof-sandwich-relation-gauss-seq-model}
The functionals $R_k^\star$, $R^\star$ and minimax risks $\mathfrak{M}_k$
satisfy 
\begin{subequations}
\begin{align}
\frac{1}{4} R_k^\star &\leq \mathfrak{M}_k \leq R_k^\star \quad \mbox{for all $k \geq 1$, and,} \label{ineq:sandwich-relation-k-trunc-functional-gauss-seq-model} \\ 
&\lim_{k \to \infty} R_k^\star = R^\star. 
\label{eqn:limit-relation-k-trunc-functional-gauss-seq-model}
\end{align}
\end{subequations}
\end{lemma} 
Assuming Lemma~\ref{lem:proof-sandwich-relation-gauss-seq-model} for the moment, 
note that the two inequalities immediately imply the sandwich relation~\eqref{eqn:sandwich-relation-gauss-sequence-model}, simply by applying the 
sandwich~\eqref{ineq:sandwich-relation-k-trunc-functional-gauss-seq-model} 
to the terms $\mathfrak{M}_k$ and then applying the limit relations~\eqref{eqn:limit-relation-gauss-seq-model} and~\eqref{eqn:limit-relation-k-trunc-functional-gauss-seq-model}. 
Consequently, it suffices to establish Lemma~\ref{lem:proof-sandwich-relation-gauss-seq-model}.

\myparagraph{Proof of Lemma~\ref{lem:proof-sandwich-relation-gauss-seq-model}}
Recall the settings of the parameters 
$T^{(k)}, \NoiseCovariance^{(k)}, \EstimationMat^{(k)}, \radius^{(k)}, 
\ConstraintMat^{(k)}$, corresponding to the $k$ dimensional minimax risk $\mathfrak{M}_k$, 
as given in~\eqref{eqn:parameters-k-trunc-gauss-sequence-model}. 
We claim that
\begin{equation}\label{eqn:equiv-functionals-k-trunc-gauss-seq-model}
\Phi(T^{(k)}, \bP,  \NoiseCovariance^{(k)}, \radius^{(k)}, \EstimationMat^{(k)},  
\ConstraintMat^{(k)}) = 
R_k^\star.
\end{equation}
(Note by our construction of $T^{(k)}$ the choice of $\bP$ is irrelevant.) 
Then the sandwich relation~\eqref{ineq:sandwich-relation-k-trunc-functional-gauss-seq-model}
follows by applying Theorems~\ref{thm:main-upper} and~\ref{thm:main-lower}
to the minimax risk $\mathfrak{M}_k$. 

To see that relation~\eqref{eqn:equiv-functionals-k-trunc-gauss-seq-model}
holds, note that by definition~\ref{eqn:defn-functional-general}, we have 
\begin{align*}
\Phi(T^{(k)}, \bP,  \NoiseCovariance^{(k)}, \radius^{(k)}, \EstimationMat^{(k)}, 
\ConstraintMat^{(k)})
= \sup_{\Omega \succ 0}  \Big\{\, \trace\Big( 
(\Omega^{-1} + (\NoiseCovariance^{(k)})^{-1})^{-1} \Big) 
: \sum_{j=1}^k a_j^2 \Omega_{jj} \leq C^2 \,\Big\}.
\end{align*}
We claim that the supremum above can be reduced to 
diagonal $\Omega$. To see why, first note that for every nonzero $\lambda \in \R$
\begin{align*}
\big(\Omega^{-1} + (\NoiseCovariance^{(k)})^{-1} \big)^{-1} 
\preceq 
\lambda^2 \Omega + (1-\lambda)^2 \NoiseCovariance^{(k)}.  
\end{align*}
This follows from Lemma~\ref{lem:matrix-lower-bound}, with the choices 
\begin{align*}
A = \NoiseCovariance^{(k)}, \quad B = \Omega^{-1}, 
\quad \mbox{and} \quad 
D = \lambda I. 
\end{align*}
Consequently, we have for every nonzero $u \in \R^k$, that 
\begin{align*}
u^\T \big(\Omega^{-1} + (\NoiseCovariance^{(k)})^{-1} \big)^{-1} u 
\leq 
\inf_{\lambda \in \R}
\lambda^2 u^\T \Omega u + 
(1-\lambda)^2 u^\T \NoiseCovariance^{(k)}u 
= \Big(\frac{1}{u^\T \Omega u} + \frac{1}{u^\T \NoiseCovariance^{(k)} u} \Big)^{-1}. 
\end{align*}
Hence taking $u$ to be elements of the standard basis $e_i$, and summing 
over $i = 1, \dots, k$, we obtain, 
\begin{align*}
\trace\Big(\big(\Omega^{-1} + (\NoiseCovariance^{(k)})^{-1} \big)^{-1} \Big) 
\leq \sum_{i=1}^k 
\Big(\frac{1}{\Omega_{ii}} + \frac{1}{\eps_i^2} \Big)^{-1} = 
\sum_{i=1}^k \frac{\Omega_{ii} \eps_i^2}{\Omega_{ii} + \eps_i^2}. 
\end{align*}
Moreover, by taking $\Omega$ to be diagonal, the 
inequality above holds with equality. 
Thus, 
\begin{align*}
\Phi(T^{(k)}, \bP,  \NoiseCovariance^{(k)}, \radius^{(k)}, \EstimationMat^{(k)}, 
\ConstraintMat^{(k)})
&= 
\sup_{\Omega_{jj} > 0 }  \Big\{\, 
\sum_{j=1}^k \frac{\Omega_{jj} \eps_j^2}{\Omega_{jj} + \eps_j^2}
: \sum_{j=1}^k a_j^2 \Omega_{jj} \leq C^2 \,\Big\}  \\ 
&=
\sup_{\tau_j^2 > 0}  \Big\{\, 
\sum_{j=1}^k \frac{\tau_j^2 \eps_j^2}{\tau_j^2 + \eps_j^2}
: \sum_{j=1}^k a_j^2 \tau_j^2 \leq C^2 \,\Big\}   \\ 
&= R^\star_k,
\end{align*}
which establishes the relation~\eqref{eqn:equiv-functionals-k-trunc-gauss-seq-model}.
Note that in the last equality, we have dropped the inequality constraints  
$\tau_j^2 > 0$, due to the continuity of the map $\tau \mapsto 
\sum_{i=1}^k \frac{\tau_j^2 \eps_j^2}{\tau_j^2 + \eps_j^2}$ over $\tau \in \R^k$. 

We now turn to establishing 
the relation~\eqref{eqn:limit-relation-k-trunc-functional-gauss-seq-model}. 
Note that for any $\tau \in \R^\N$ with $\sum_{j=1}^\infty a_j^2 \tau_j^2 \leq C^2$, 
we have 
\begin{align*}
\sum_{j=1}^k \frac{\tau_j^2 \eps_j^2}{
  \tau_j^2 + \eps_j^2
} 
\leq 
\sum_{j=1}^\infty \frac{\tau_j^2 \eps_j^2}{
  \tau_j^2 + \eps_j^2
} 
\leq 
\sum_{j=1}^k \frac{\tau_j^2 \eps_j^2}{
  \tau_j^2 + \eps_j^2
} 
+ \sup_{\tau \in \R^\N : \sum_{j=1}^\infty a_j^2 \tau_j^2 \leq C^2}
\sum_{j > k}^\infty \tau_j^2 
\end{align*}
By Hölder's inequality, the second term is bounded above by $C^2/a_{k+1}^{2}$, hence
in view of definitions~\eqref{eqn:functional-k-trunc-gauss-seq-model} 
and~\eqref{eqn:functional-full-gauss-seq-model}, we have the sandwich relation 
\begin{align*}
R_k^\star \leq R^\star \leq R_k^\star + \frac{C^2}{a_{k+1}^{2}}, 
\end{align*}
which holds for all $k \geq 1$. Since $a_{k} \to \infty$, the 
limit relation~\eqref{eqn:limit-relation-k-trunc-functional-gauss-seq-model} follows.

\subsubsection{Proof of limit relation~\eqref{eqn:limit-relation-rkhs-model}}
\label{app:proof-limit-relation-rkhs-model}

We claim that the following sandwich relation holds for the minimax risks 
in this case. 
\begin{lemma}\label{lem:sandwich-relation-rkhs-model-truncation}
For all $k \geq 1$, we have 
\begin{equation}\label{eqn:lem-sandwich-relation-trunc-rkhs}
\mathfrak{M}^{(k)}_\numobs(\radius, \sigma^2, P) \leq 
\mathfrak{M}_\numobs(\radius, \sigma^2, P)  \leq 
\mathfrak{M}^{(k)}_\numobs(\radius, \sigma^2, P)  
+ \radius^2 \mu_{k+1}. 
\end{equation}
\end{lemma} 
\noindent Assuming Lemma~\ref{lem:sandwich-relation-rkhs-model-truncation}, note that 
since $\mu_k \to 0$ as $k \to \infty$, it immediately implies limit relation~\eqref{eqn:limit-relation-rkhs-model}

\myparagraph{Proof of Lemma~\ref{lem:sandwich-relation-rkhs-model-truncation}}
The proof is quite similar to Lemma~\ref{lem:sandwich-relation-gauss-seq-model-truncation}.
We now prove inequality~\eqref{eqn:lem-sandwich-relation-trunc-rkhs}. We begin by defining the sets 
\begin{align*}
\cB(\radius) = \{ \theta \in \ell^2(\N) : \|\theta\|_2 \leq \radius \}, \quad \mbox{and} \quad 
\cB_k(\radius) = \{\theta \in \cB_k(\radius) : \theta_j = 0,~~\mbox{for all}~j > k\}.
\end{align*}
By Parseval's identity, we may rewrite the minimax risks in the following form 
\begin{subequations}
\begin{align}
\mathfrak{M}_k &\equiv \mathfrak{M}^{(k)}_\numobs(\radius, \sigma^2, P) = 
\inf_{\thetahat} \sup_{\substack{\thetastar \in \cB_k(\radius)\\\nu \in \cP(\sigma^2 I_\numobs)}} 
\E\Big[\sum_{j=1}^k \mu_j (\thetahat_j(y_1, \dots, y_\numobs, 
\Phi_k(x_1), \dots, \Phi_k(x_\numobs)) - \thetastar_j)^2\Big],
\label{eqn:k-trunc-minimax-risk-parametric-form-rkhs}\\ 
\mathfrak{M} &\equiv \mathfrak{M}_\numobs(\radius, \sigma^2, P) = 
\inf_{\thetahat} \sup_{\substack{\thetastar \in \cB(\radius)\\\nu \in \cP(\sigma^2 I_\numobs)}} 
\E\Big[\sum_{j=1}^\infty \mu_j (\thetahat_j
(y_1, \dots, y_\numobs, 
\Phi(x_1), \dots, \Phi(x_\numobs)) - \thetastar_j)^2\Big]. 
\label{eqn:full-minimax-risk-parametric-form-rkhs}
\end{align}
\end{subequations}
Evidently, we have $\mathfrak{M} \geq \mathfrak{M}_k$, 
since $\cB_k(\radius) \subset \cB(\radius)$ and $(y, \Phi_k(x))$ are sufficient in this submodel. 
Similarly, the upper bound follows since by restricting to those estimators 
$\thetahat$ with $\thetahat_j = 0$ 
for all $j > k$ that are functions of $(y, \Phi_k(x))$, we have 
\begin{align*}
\mathfrak{M} \leq \mathfrak{M}_k + 
\sup_{\theta \in \cB(\radius)} \sum_{j > k} \mu_j \theta_j^2 
 = \mathfrak{M}_k + \radius^2 \mu_{k + 1}, 
\end{align*}
which establishes the upper bound. 

\subsubsection{Proof of relation~\eqref{eqn:sandwich-relation-rkhs-model}}
\label{app:proof-sandwich-relation-rkhs-model}

Applying Corollary~\ref{cor:iid-result} to the minimax risk 
$\mathfrak{M}_k(\radius, \sigma^2, P)$, we find that 
\begin{align*}
\frac{1}{4} \, \frac{\sigma^2}{\numobs} d_\numobs^{(k)} 
\leq \mathfrak{M}_k(\radius, \sigma^2, P) 
\leq \frac{\sigma^2}{\numobs} d_\numobs^{(k)}, 
\end{align*}
since the quantity $d_\numobs^{(k)}$ equals the functional 
for this minimax risk (see equation~\eqref{eqn:functional-k-trunc-rkhs-model}).
Therefore passing to the superior limit and applying the 
limit relation~\eqref{eqn:limit-relation-rkhs-model}, we obtain the result.

\subsubsection{Proof of relation~\eqref{eqn:loosened-characterization-rkhs}}
\label{app:proof-loosened-characterization-rkhs-model}

First, we define the population counterparts of the functional 
$d_\numobs^{(k)}$, as defined in equation~\eqref{eqn:functional-k-trunc-rkhs-model}.
 Note that under $P$, we have $\E \Sigma_{\numobs}^{(k)} = \diag(\mu_1, \dots, \mu_k)$. 
 We denote this matrix by $M_k$. Hence, the population 
 counterpart to $d_\numobs^{(k)}$ is 
\begin{align}
\overline{d}_{\numobs}^{(k)} &\defn 
\sup_{\Omega \succ 0} \,
\Big\{\trace\Big(M_k( M_k  +
\Omega^{-1})^{-1}\Big) :  \trace(\Omega) \leq \frac{\numobs
  \radius^2}{\sigma^2}\Big\} \\ 
&= 
\sup_{\Omega \succ 0} \,
\Big\{\trace\Big(( I_k  +
\Omega^{-1})^{-1}\Big) :  \trace(M_k^{-1} \Omega) \leq \frac{\numobs
  \radius^2}{\sigma^2}\Big\}.
\end{align}
Using Proposition~\ref{prop:to-population} and the sandwich relation~\eqref{eqn:lem-sandwich-relation-trunc-rkhs}, we find 
\begin{align*}
\frac{1}{4}\, \frac{\sigma^2 }{\numobs} \overline{d}_{\numobs}^{(k)} 
\leq \minimaxrisk_{\numobs}^{(k)}(\radius, \sigma^2, P) \leq 
\Big(1 + \frac{\kappa^2 \radius^2}{\sigma^2}\Big)  \frac{\sigma^2 }{\numobs} \overline{d}_{\numobs}^{(k)}  + \mu_{k+1} \radius^2.
\end{align*}
Since $\mu_k \to 0$ as $j \to \infty$, it suffices to show that 
\begin{equation}
\label{eqn:limit-relation-trunc}
\lim_{j \to \infty} \overline{d}_{n}^{(k)} = \overline{d}_{n}^\star.
\end{equation}

\myparagraph{Proof of relation~\eqref{eqn:limit-relation-trunc}}

This limit relation can be established via an argument based on
Lagrange multipliers.  First, by an eigendecomposition of the variable
$\Omega \succ 0$, we have
\begin{align}
\overline{d}_{n}^{(k)} & = \sup\Big\{\sum_{j=1}^k \frac{\tau_j}{1 +
  \tau_j} : \tau_j > 0, \sum_{j = 1}^k \frac{\tau_j}{\mu_j} \leq
\frac{\numobs \radius^2}{\sigma^2} \Big\} \nonumber \\
\label{eqn:d-reformulation}
& = \sup\Big\{\sum_{j=1}^k \frac{\mu_j
  \gamma_j}{\tfrac{\sigma^2}{\numobs \radius^2} + \mu_j\gamma_j} :
\gamma_j \geq 0, \sum_{j = 1}^k \gamma_j \leq 1 \Big\}.
\end{align} 
The final equality arises by a rescaling and continuity argument. Note
that we may drop the nonnegativity constraint, since the sequence
$\{\mu_j\}$ is nonnegative.  We can
compute~\eqref{eqn:d-reformulation} by introducing dual variables.  In
particular, we have
\begin{align}
\label{eqn:d-variational-reformulation}   
\overline{d}_{n}^{(k)} & = \sup_{\gamma} \inf_{\lambda} \sum_{j=1}^k
\frac{\mu_j \gamma_j }{\tfrac{\sigma^2}{\numobs \radius^2} +
  \mu_j\gamma_j} - \frac{\tfrac{\numobs
    \radius^2}{\sigma^2}}{\lambda^2} \Big(\sum_{j=1}^k \gamma_j -
1\Big).
\end{align}
By simple calculus, we see that the saddle point $(\gamma^\star,
\lambda^\star)$ satisfies
\begin{align*}
\quad \sum_{j=1}^k \tau_j^\star = 1 \quad \mbox{and} \quad
\Big(\frac{\sigma^2}{\numobs \radius^2}\Big)^2
\frac{\mu_j}{(\tfrac{\sigma^2}{\numobs \radius^2} + \mu_j
  \gamma_j^\star)^2} = \frac{1}{(\lambda^\star)^2}, ~~\mbox{for}~j =
1, \ldots, k.
\end{align*}
Using the fact that $\gamma_j^\star \geq 0$, we obtain $\gamma_j^\star
= \tfrac{\sigma^2}{\numobs \radius^2}
\tfrac{1}{\sqrt{\mu_j}}(\lambda^\star - \tfrac{1}{\sqrt{\mu_j}})_+$,
where $\lambda^\star$ is chosen such that
\begin{align*}
\frac{\sigma^2}{\numobs \radius^2}\sum_{j=1}^k \frac{1}{\sqrt{\mu_j}}
\Big(\lambda^\star - \frac{1}{\sqrt{\mu_j}} \Big)_+ =
\frac{\sigma^2}{\numobs \radius^2} \sum_{j=1}^k \gamma^\star_j = 1 .
\end{align*}
Using equation~\eqref{eqn:d-reformulation}, it follows that
\begin{align*}
\overline{d}_{\numobs}^{(k)} = \sum_{j=1}^k
\frac{1}{\lambda^\star}\Big(\lambda^\star -
\tfrac{1}{\sqrt{\mu_j}}\Big)_+.
\end{align*}
The result then follows by appealing to the following numerical
result, with $a_j = 1/\sqrt{\mu_j}$.

\begin{lemma}
\label{lem:numerical-seq}
Let $\{a_j\}_{j \geq 1}$ denote a nonnegative, divergent
sequence,\footnote{Formally, $\{a_j\} \subset \R_+$ and $\lim_{j \to
  \infty} a_j = +\infty$.}  and define $a_\star \defn \inf_{j \geq 1}
a_j$.  Consider the functions $f_n, f \colon [a_\star, +\infty) \to
  \R_+$ given by
\begin{align*}
f_n(t) \defn \sum_{k=1}^n a_k(t - a_k)_+ \quad \mbox{and} \quad f(t)
\defn \sum_{k=1}^\infty a_k(t - a_k)_+,
\end{align*}
and define $\tau_n$ and $\tau$ via the relations $f_n(\tau_n) =
f_n(\tau) = 1$.  Then:
\begin{enumerate}
\item [(i)] The function $f$ and values $\tau_n, \tau$ are
  well-defined; and
  \item [(ii)] We have $\tau_n = \tau$ for $n$ sufficiently large.
\end{enumerate}
\end{lemma}

\myparagraph{Proof of Lemma~\ref{lem:numerical-seq}}

Since the sequence $a_k$ diverges to infinity, we may assume without
loss of generality that $a_k > 0$ for all $k$.  For the first claim,
note that $f$ is well-defined.  Indeed fix $t \geq a_\star$. Then,
there exists $n$ sufficiently large such that $t < a_k$ for all $k
\geq n$.  Consequently, $f(t) = f_n(t)$.  Similarly, note that $f_n,
f$ are strictly increasing, continuous functions with $f(a_\star) =
f_n(a_\star) = 0$, and $f(x), f_n(x) \to \infty$ in the limit as $x
\to \infty$.  Therefore, $\tau_n, \tau$ exist and are uniquely defined
by the equations $f_n(\tau_n) = 1$ and $f(\tau) = 1$, respectively.
By the argument given previously, $f_n(\tau) = f(\tau)$ for all $n$
large enough, and therefore, by uniqueness $\tau = \tau_n$ for $n$
large enough.

\subsubsection{Proof of Sobolev rate calculation}
\label{app:proof-sobolev-rate-calculation}

Let $\{\mu_j\}_{j \geq 1}$ denote the eigenvalue sequence associated
to the integral operator for the order $\beta$ Sobolev space on $[0,
  1]^\dimension$.  We then define the auxiliary functions
\begin{align*}
f(\lambda) \defn \sum_{k=1}^\infty \frac{1}{\sqrt{\mu_k}} \Big(
\lambda - \frac{1}{\sqrt{\mu_k}}\Big)_+ \quad \mbox{and} \quad
d_\numobs(\lambda) \defn \sum_{k=1}^\infty
\frac{1}{\lambda}\Big(\lambda - \frac{1}{\sqrt{\mu_k}}\Big)_+.
\end{align*}
In view of relation~\eqref{eqn:loosened-characterization-rkhs}, 
it follows that the minimax risk over the ball of radius $\radius > 0$ 
within the order-$\beta$ Sobolev space in $[0, 1]^\dimension$ 
is equal (up to constant pre-factors) to
\begin{equation}
\label{eqn:from-corollary}
\frac{\sigma^2}{\numobs} d_\numobs(\lambda_\numobs^\star) \quad
\mbox{where} \quad f(\lambda_\numobs^\star) = \frac{\numobs
  \radius^2}{\sigma^2},
\end{equation}
whenever $\radius \lesssim \sigma$.\footnote{In this subsection,
we allow the relations $\asymp, \lesssim, \gtrsim$ to hide 
constants which depend on $\beta, d$ but not on $\numobs, \radius,
\sigma$.}  In order to simplify the description of the rate above, we
claim that
\begin{equation}
\label{eqn:sob-scaling}
d_\numobs(\lambda_\numobs^\star) \asymp \Big(\frac{\sigma^2}{\numobs
  \radius^2}\Big)^{-\frac{d}{2\beta + d}}.
\end{equation}
Assuming equation~\eqref{eqn:sob-scaling} for the moment, combination
with display~\eqref{eqn:from-corollary} yields the minimax risk, which
is $\radius^2 (\tfrac{\sigma^2}{\numobs
  \radius^2})^{\tfrac{2\beta}{2\beta + d}}$, up to constant factors.
This is the claimed result.

\myparagraph{Proof of relation~\eqref{eqn:sob-scaling}}

We begin by determining $\lambda_\numobs^\star$, apart from constants.
For $\beta > d/2$, the eigenvalues $\mu_j$ satisfy $\mu_j \asymp
j^{-2\alpha}$, where $\alpha \defn \beta / d$. Therefore, it follows
that
\begin{align*}
f(\lambda) \asymp g(\lambda) \defn \sum_{k=1}^\infty k^{\alpha}
(\lambda - k^\alpha)_+.
\end{align*}
Note that both $f$ and $g$ are increasing functions. If $g(\lambda)
\asymp g(\lambda')$, it follows that $\lambda \asymp \lambda'$, since
$g$ is piecewise affine, and thus locally Lipschitz.  It follows that
$\lambda^\star_\numobs \asymp \tilde \lambda_\numobs^\star$, where
$g(\tilde \lambda_\numobs^\star) \asymp \frac{\numobs
  \radius^2}{\sigma^2}$.  A similar argument shows that
\begin{align*}
d_\numobs(\lambda) \asymp \tilde d_\numobs(\lambda) \defn
\sum_{k=1}^\infty \frac{(\lambda - k^\alpha)_+}{\lambda}
\end{align*}

Our argument is based on establishing the following relations,
\begin{equation}
\label{eqn:scaling-relations-Sobolev}
g(\lambda) \stackrel{\rm(i)}{\asymp} \lambda^{2 + 1/\alpha} \quad
\mbox{and} \quad \tilde d_\numobs(\lambda) \stackrel{\rm (ii)}{\asymp}
\lambda^{1/\alpha}.
\end{equation}
Assuming these bounds for a moment, we explain how the claimed result
on the minimax risk follows.  First, note that since
$f(\lambda_\numobs^\star) = \tfrac{\numobs \radius^2}{\sigma^2}$, the
argument above implies that $\lambda_\numobs^\star \asymp \tilde
\lambda_\numobs^\star$ where $\tilde \lambda_\numobs^\star$ satisfies
$g(\lambda) \asymp \tfrac{\numobs \radius^2}{\sigma^2}$.  Therefore,
from equation~\eqref{eqn:scaling-relations-Sobolev}(i), it follows
that $\tilde \lambda_\numobs^\star \asymp (\tfrac{\sigma^2}{\numobs
  \radius^2})^{-\frac{\alpha}{2\alpha + 1}}$.  Then, using
equation~\eqref{eqn:scaling-relations-Sobolev}(ii), it follows that
$\tilde d_\numobs(\lambda_\numobs^\star)\asymp \tilde d_\numobs(\tilde
\lambda_\numobs^\star) \asymp (\tfrac{\sigma^2}{\numobs
  \radius^2})^{-\frac{1}{2\alpha + 1}}$, which establishes the claimed
inequality~\eqref{eqn:sob-scaling}, after recalling $\alpha =
\beta/d$, and clearing the denominator of the exponent.

We now demonstrate scaling
relation~\eqref{eqn:scaling-relations-Sobolev}(i), so that we show
that $g(\lambda) \asymp \lambda^{2 + 1/\alpha}$, for all $\lambda \geq
1$.  In order to establish this claim, choose the integer $k$ such
that $\lambda \in (k^\alpha, (k+1)^{\alpha}]$. Then
\begin{align*}
g(\lambda) \leq \lambda \sum_{j=1}^k j^\alpha \leq \lambda
\frac{(k+1)^{\alpha + 1}}{\alpha + 1} \lesssim \lambda k^{\alpha + 1}
\lesssim \lambda^{2 + 1/\alpha}.
\end{align*}
Above, we used an integral approximation for the summation.  On the
other hand, when $\lambda \in (k^\alpha, (k+1)^{\alpha}]$, we have
\begin{align*}
g(\lambda) \geq g(k^\alpha) \geq (k^\alpha - \ceil{k/2}^\alpha)
\sum_{j=1}^{\ceil{k/2}} j^\alpha \gtrsim k^{2\alpha + 1}.
\end{align*}
To simplify, the last equality (up to constants) is obtained by an
integration argument.  Therefore, we have
\begin{align*}
\inf_{k \geq 1} \inf_{\lambda \in (k^\alpha, (k+1)^\alpha]} \frac{g(\lambda)}{\lambda^{2 + 1/\alpha}} 
\geq \inf_{k \geq 1} \frac{g(k^\alpha)}{(k+1)^{2\alpha + 1}} 
\gtrsim 1. 
\end{align*}
Thus, we have the bound $g(\lambda) \gtrsim \lambda^{2 + 1/\alpha}$
for all $\lambda$, as needed.

We now demonstrate the scaling
relation~\eqref{eqn:scaling-relations-Sobolev}(ii), so that we show
$\tilde d_\numobs(\lambda) \asymp \lambda^{1/\alpha}$.  To see this,
note first that for $\lambda \in (k^\alpha, (k+1)^\alpha]$, we have
  the trivial bound
\begin{align*}
\tilde d_\numobs(\lambda) = \sum_{j=1}^k (1 - \lambda^{-1} j^\alpha)_+
\leq k \lesssim \lambda^{1/\alpha}.
\end{align*}
On the other hand, we have the lower bound
\begin{align*}
\tilde d_\numobs(\lambda) \geq \sum_{j=1}^{\ceil{k/2}} \frac{(k^\alpha
  - j^\alpha)}{(k+1)^\alpha} \geq \frac{k+1}{2} \cdot
\Big(\frac{k}{k+1} \frac{(k^\alpha -
  \ceil{k/2}^\alpha)}{(k+1)^\alpha}\Big) \gtrsim k+1 \gtrsim
\lambda^{1/\alpha}.
\end{align*}

\subsubsection{Proof of relation~\eqref{eqn:very-loose-covshift-lower}}
\label{app:proof-of-covshift-relation}

Note that the kernel regularity condition is not necessary for our lower bound. 
Indeed, note that we first have 
\begin{align*}
\inf_{\delta > 0} \Big\{\delta^2 + \frac{\sigma^2 B}{\numobs \radius^2}
d(\delta)\Big\} &= 
\inf_{d \geq 1 } 
\Big\{\, 
\mu_{d} + 
\frac{\sigma^2 B d}{\numobs \radius^2} \, \Big\}
\end{align*}
Let $d_\numobs^\star$ be the largest integer $d$ such that 
$\mu_d \geq \tfrac{\sigma^2 B d}{\numobs \radius^2}$; this must 
exist since $\mu_d \to 0$. As the two sequences are nonincreasing and 
strictly increasing, respectively, the display above is bounded above by 
\begin{align*}
4 \Big( \twomin{ \mu_{d_\numobs^\star}}{  \frac{\sigma^2 B d_\numobs^\star}{\numobs \radius^2}}\Big)
\leq 4  \frac{\sigma^2 B d_\numobs^\star}{\numobs \radius^2}. 
\end{align*}
Hence, it suffices to establish that the lower bound
$\frac{\sigma^2 B d_\numobs^\star}{\numobs \radius^2}$ can be obtained 
from our result~\eqref{eqn:loose-covshift-lower}. 

Note that if $\mu_d \geq \frac{\sigma^2 B d}{\numobs \radius^2}$
then the choice of $\lambda$ in the lower bound~\eqref{eqn:loose-covshift-lower}, 
given by 
\begin{align*}
\lambda_j = \frac{\sigma^2 B}{\numobs \radius^2} \frac{1}{\mu_j} \1\{j \leq d\}, 
\quad \mbox{for}~j = 1, 2, 3, \dots,
\end{align*}
satisfies $\sum_j \lambda_j \leq 1$. Evaluating the corresponding lower bound, 
with the maximal choice $d = d_\numobs^\star$ yields the lower bound
$\tfrac{\sigma^2 B d}{\numobs \radius^2}$, as needed.

\section{Proofs and calculations from Section~\ref{sec:proof-main}}
\label{app:proofs}
\subsection{Deferred proofs from Section~\ref{sec:proof-upper}}\label{app:proofs-upper}

In this section, we collect proofs of the results underlying the argument establishing our 
upper bound in Section~\ref{sec:proof-upper} of the paper. 

\subsubsection{Proof of Lemma~\ref{lem:suprema-equality}}\label{app:proof-suprema-equality}

Clearly the lefthand side is less than the right hand side as for
$\theta \in \Theta(\radius, \ConstraintMat)$ we have $\theta \otimes
\theta \succeq 0$, and $\trace(\ConstraintMat^{\!-1/2} \theta \otimes
\theta \ConstraintMat^{\!-1/2}) = \vecnorm{\theta}_{\ConstraintMat^{-1}}^2\leq
\radius^2$.  

For the reverse inequality, fix $\Omega \in \cK(\radius,
\ConstraintMat)$.  We diagonalize the positive semidefinite matrix
$\ConstraintMat^{\!-1/2} \Omega \ConstraintMat^{\!-1/2} = UDU^\T$, and define
$\theta(\eps) = \ConstraintMat^{\!1/2} UD^{1/2} \eps$, where $\eps \in
\{\pm 1\}^\dimension$.  Evidently,
\begin{align*}
  \vecnorm{\theta(\eps)}_{\ConstraintMat^{-1}}^2 = \twonorm{U D^{1/2} \eps}^2 = \trace(D) = \trace(\ConstraintMat^{\!-1/2} \Omega \ConstraintMat^{\!-1/2}) \leq \radius^2.
\end{align*}
Thus, for all $\eps \in \{\pm 1\}^\dimension$, the vector $\theta(\eps)$ lies in the set $\Theta(\radius, \ConstraintMat)$. 
Consequently, we have 
\begin{align}
\sup_{\theta \in \Theta(\radius, \ConstraintMat)} r(\hat
\theta_C, \theta) 
&\geq \max_{\eps \in \{\pm 1\}^\dimension}
r(\thetahat_C, \theta(\eps)) \nonumber\\
& \geq \E_{\eps} r(\thetahat_C,
\theta(\eps)) \label{eqn:replace-max-by-expectation} \\
& =  r(\thetahat_C, \Omega). \label{eqn:lin-of-exp}
\end{align}
Note that $\Omega \in \cK(\radius, \ConstraintMat)$ was arbitrary in this argument, and hence 
passing to supremum over $\Omega$ gives us the desired reverse inequality. 
Above, display~\eqref{eqn:replace-max-by-expectation} follows by lower bounding the maximum over $\eps \in
\{\pm 1\}^\dimension$ by the expectation over $\eps$ where $\eps_i$ are
\iid{} Rademacher variables.  
The relation~\eqref{eqn:lin-of-exp} follows by noting that $r(\thetahat_C,
\theta(\eps)) = r(\thetahat_C,
\theta(\eps) \otimes \theta(\eps))$, and moreover this latter quantity is linear in the rank-one matrix 
$\theta(\eps) \otimes \theta(\eps)$, as justified by Lemma~\ref{lem:curvature-of-linear-risk}. 
By linearity of expectation we can bring the expectation inside, and use the fact that
\begin{align*}
  \E_\eps[\theta(\eps) \otimes \theta(\eps)] = \ConstraintMat^{\!1/2} UDU^\T \ConstraintMat^{\!1/2} = \Omega. 
\end{align*}

\subsubsection{Proof of Lemma~\ref{lem:curvature-of-linear-risk}}\label{app:proof-curvature-risk}

Inspecting the definition of $r$ (see equation~\eqref{eqn:definition-linear-risk-matrix-version}), 
we see that it is affine in $\Omega$. To verify that it is convex in $C$, note that $r$ can be equivalently 
expressed as 
\begin{multline*}
r(\thetahat_C, \Omega) = 
\E_\xi \Big[\fronorm{\EstimationMat^{\!1/2} (C(\LinOp{\xi}) \LinOp{\xi}^\T \NoiseCovariance^{-1} \LinOp{\xi} - I_\dimension) \Omega^{1/2}}^2
+ \fronorm{\EstimationMat^{\!1/2} (C(\LinOp{\xi}) \LinOp{\xi}^\T \NoiseCovariance^{-1/2}}^2\Big].
\end{multline*}
Evidently, the display above is convex in $C$. 

\subsubsection{Proof of Proposition~\ref{prop:minimizer-of-linear-risk}}\label{app:proof-minimizer-linear-risk}

In order to prove Proposition~\ref{prop:minimizer-of-linear-risk}, we
need two results regarding the harmonic mean of positive
(semi)definite matrices.  For our results, it is important to allow
once of these matrices to be (possibly) singular, and so we study
(twice) the harmonic mean of $A$ and the Moore-Penrose pseudoinverse
$B^\dagger$---that is, the quantity $(A^{-1} + B)^{-1}$, where $B
\succeq 0$ and $A \succ 0$. Note that since $(B^\dagger)^\dagger = B$,
these results also imply bounds for the mean $(A^{-1} +
B^\dagger)^{-1}$.  See the reference~\cite[chap.~4]{Bha07} for
additional details about the harmonic mean of positive matrices.

\begin{lemma}
\label{lem:semidefinite-result}
Suppose that $A, B$ are two symmetric positive semidefinite matrices, and that $A$ is nonsingular.
For any $x \in \R^\dimension$ and any $y$ in the range of $B$, we have
\begin{align*}
(x - y)^\T A (x-y) + y^\T B^\dagger y \geq x^\T (A^{-1} + B)^{-1} x,
\end{align*}
where $B^\dagger$ denotes the Moore-Penrose pseudoinverse associated with $B$. 
\end{lemma}
\begin{proof}
Using $BB^\dagger B = B$, the claim is equivalent to showing that
$\inf_{x, u} g(x, u) \geq 0$ where
\begin{align*}
g(x, u) \defn (x - Bu)^\T A (x- Bu) + u^\T B u- x^\T(A^{-1} + B)^{-1}
x.
\end{align*}
Define $f(u) = \inf_{x} g(x, u)$. A calculation demonstrates that
\begin{align}
\label{eqn:f-calculation}  
  f(u) & = u^\T \Big[B + BAB - BA(A - (A^{-1} + B)^{-1})^\dagger A B
    \Big] u \nonumber \\ &= u^\T BA^{1/2} \Big[K^\dagger + I - (I - (I
    + K)^{-1})^\dagger\Big] A^{1/2} B u.
\end{align}
Above, $K \defn A^{1/2} B A^{1/2}$. Diagonalizing $K$, we may write $K
= UDU^\T$ and therefore $K^\dagger = UD^\dagger U^T$. Applying the
similarity transformation under $U$, we have
\begin{equation}
\label{eqn:matrix-ineq-sim}
U^\T(K^\dagger + I - (I - (I + K)^{-1})^\dagger)U = D^\dagger + I - (I
- (I + D)^{-1})^\dagger = I - D^\dagger D \succeq 0.
\end{equation}
Therefore, combining displays~\eqref{eqn:f-calculation}
with~\eqref{eqn:matrix-ineq-sim}, we obtain
\begin{align*}
\inf_{x, u} g(x, u) = \inf_{u} f(u) \geq 0,
\end{align*}
which establishes the desired claim.
\end{proof}
\begin{lemma}
\label{lem:matrix-lower-bound}
Suppose that $A, B$ are two symmetric positive semidefinite matrices, and that $A$ is nonsingular.
If $D^\T \in \R^{\dimension \times \dimension}$ has range included in the
range of $B$, then
\begin{align*}
(I - D) A (I - D)^\T + D B^\dagger D^\T \succeq (A^{-1} + B)^{-1}.
\end{align*}
Moreover equality holds with the choice $D = (A^{-1} + B)^{-1} B$.
\end{lemma}
\begin{proof}
Let $x \in \R^\dimension$ and note that if $y \defn D^\T x$, then
\begin{align*}
x^\T \Big[(I - D) A (I - D)^\T + D B^\dagger D^\T\Big] x &= (x-y)^\T A
(x- y) + y^\T B^\dagger y \\ &\geq x^\T (A^{-1} + B)^{-1} x,
\end{align*}
where the final inequality follows
from Lemma~\ref{lem:semidefinite-result}, since $y$ lies in the range of
$B$.  As the inequality holds for arbitrary $x \in \R^\dimension$, we
have established the desired matrix inequality.  To see the attainment
at $D = (A^{-1} + B)^{-1}B$, first note that $D^\T = B(A^{-1} + B)^{-1}$.  Therefore the
range of $D^\T$ is exactly the range of $B$. Additionally, since $I -
D = (A^{-1} + B)^{-1}A^{-1}$, we have
\begin{align*}
(I - D) A(I-D)^\T + D B^\dagger D^\T = (A^{-1} + B)^{-1} (A^{-1} + B
  B^\dagger B) (A^{-1} + B)^{-1} = (A^{-1} + B)^{-1},
\end{align*}
as required.
\end{proof}

We are now in a situation to prove Proposition~\ref{prop:minimizer-of-linear-risk}. 

\myparagraph{Proof of Proposition~\ref{prop:minimizer-of-linear-risk}}

From display~\eqref{eqn:definition-linear-risk-matrix-version}, to establish the claim, it suffices
to lower bound the following matrix in the semidefinite ordering,
\begin{multline}\label{eqn:matrix-to-lower-bound}
(C(\LinOp{\xi}) \LinOp{\xi}^\T \NoiseCovariance^{-1} \LinOp{\xi} - I_\dimension) \Omega (C(\LinOp{\xi}) \LinOp{\xi}^\T \NoiseCovariance^{-1} \LinOp{\xi} - I_\dimension)^\T 
 \\+ C(\LinOp{\xi}) \LinOp{\xi}^\T \NoiseCovariance^{-1} \LinOp{\xi} C(\LinOp{\xi})^\T. 
 \end{multline}
This matrix can be written as $(I - D) \Omega(I- D)^\T + D B^\dagger D^\T$ where we defined 
\begin{align*}
B \defn \LinOp{\xi}^\T \NoiseCovariance^{-1} \LinOp{\xi}, \quad \mbox{and,} \quad 
D \defn C(\LinOp{\xi}) \LinOp{\xi}^\T \NoiseCovariance^{-1} \LinOp{\xi}.
\end{align*}
Evidently, the range of $D^\T$ is included in the range of $B$, and so it follows from Lemma~\ref{lem:matrix-lower-bound}
that the matrix in equation~\eqref{eqn:matrix-to-lower-bound} is lower bounded in the semidefinite ordering by 
\begin{equation}\label{eqn:final-matrix-lower}
(\Omega^{-1} +  \LinOp{\xi}^\T \NoiseCovariance^{-1} \LinOp{\xi})^{-1}. 
\end{equation}
Moreover, Lemma~\ref{lem:matrix-lower-bound} also demonstrates this is established by taking 
\begin{align*}
D = (\Omega^{-1} +  \LinOp{\xi}^\T \NoiseCovariance^{-1} \LinOp{\xi})^{-1}\LinOp{\xi}^\T \NoiseCovariance^{-1} \LinOp{\xi},
\end{align*}
which arises from taking $C(\LinOp{\xi}) = (\Omega^{-1} +  \LinOp{\xi}^\T \NoiseCovariance^{-1} \LinOp{\xi})^{-1}$, as claimed.  
Evaluating this lower bound matrix~\eqref{eqn:final-matrix-lower} 
in~\eqref{eqn:definition-linear-risk-matrix-version} establishes equality~\eqref{eqn:minimizer-linear-risk}.

\subsubsection{Proof of equation~\eqref{eqn:upper-bound-conclusion}}
\label{app:proof-upper-conclusion}

Let us formally state our claim, equivalent to
equation~\eqref{eqn:upper-bound-conclusion}, as a lemma.
\begin{lemma} 
Let $\cK_+(\radius, \ConstraintMat)$ denote the subset of nonsingular
matrices in $\cK(\radius, \ConstraintMat)$---that is, the set
$\{\Omega \succ 0 : \Omega \in \cK(\radius, \ConstraintMat)\}$.  Then,
we have
\begin{align*}
\sup_{\Omega \in \cK(\radius, \ConstraintMat)} 
\inf_C r(\thetahat_C, \Omega) = 
\sup_{\Omega \in \cK_{+}(\radius, \ConstraintMat)} 
\inf_C r(\thetahat_C,\Omega). 
\end{align*}
\end{lemma} 

We prove this claim now. Evidently, since $\cK_+(\radius,
\ConstraintMat) \subset \cK(\radius, \ConstraintMat)$ it suffices to
show that the lefthand side is less than or equal to the righthand
side.  To begin, we note that for each $\lambda > 0$, we have
\begin{multline*}
\sup_{\Omega \in \cK(\radius, \ConstraintMat)} \inf_C r(\thetahat_C,
\Omega) \stackrel{\rm (a)}{\leq} \sup_{\Omega \in \cK(\radius,
  \ConstraintMat)} \inf_C r(\thetahat_C, \Omega + \tfrac{(\radius +
  \lambda)^2 - \radius^2}{\dimension} \ConstraintMat) \stackrel{\rm
  (b)}{\leq} \sup_{\Omega \in \cK_+(\radius + \lambda,
  \ConstraintMat)} \inf_C r(\thetahat_C, \Omega) \eqcolon f(\lambda).
\end{multline*}
Inequality (a) above follows since $r(\thetahat_C, \Omega) \leq
r(\thetahat_C, \Omega')$ for any $\Omega \preceq \Omega'$---this
follows immediately from
display~\eqref{eqn:definition-linear-risk-matrix-version}.  Here we
have taken $\Omega' \defn \Omega + \tfrac{(\radius + \lambda)^2 -
  \radius^2}{\dimension} \ConstraintMat \succeq \Omega$.  Inequality
(b) then follows by noting that $\Omega'$ is symmetric positive
(strictly) definite, and $\trace(\ConstraintMat^{\!-1/2}
\Omega'\ConstraintMat^{\!-1/2}) \leq (\radius + \lambda)^2$, since
$\Omega \in \cK(\radius, \ConstraintMat)$.  Since the displayed
relation above holds for any $\lambda > 0$, it suffices to show that
\begin{equation}
\label{eqn:suff-condition-eqn}
\inf_{\lambda > 0} f(\lambda) = f(0).
\end{equation}
By Proposition~\ref{prop:minimizer-of-linear-risk}, we have
\begin{align*}
f(\lambda) &= \sup_{\Omega} \Big\{\, \E \trace\Big(
\EstimationMat^{\!\!1/2} (\Omega^{-1} + \LinOp{\xi}^\T
\NoiseCovariance^{-1} \LinOp{\xi})^{-1} \EstimationMat^{\!\!1/2} \Big)
: \nonumber \\ &\qquad\qquad\qquad\qquad\qquad\qquad\qquad \Omega
\succ 0,\trace(\ConstraintMat^{\!-1/2} \Omega
\ConstraintMat^{\!-1/2})\leq (\radius + \lambda)^2 \,\Big\} \\ &=
\sup_{\Omega} \Big\{\, \E \trace\Big( \EstimationMat^{\!\!1/2} (
(\tfrac{\radius + \lambda}{\radius})^{-2} \Omega^{-1} + \LinOp{\xi}^\T
\NoiseCovariance^{-1} \LinOp{\xi})^{-1} \EstimationMat^{\!\!1/2} \Big)
: \nonumber \\ &\qquad\qquad\qquad\qquad\qquad\qquad\qquad \Omega
\succ 0,\trace(\ConstraintMat^{\!-1/2} \Omega
\ConstraintMat^{\!-1/2})\leq \radius^2 \,\Big\} \\ &\leq
\Big(\frac{\radius + \lambda}{\radius}\Big)^2 \; \sup_{\Omega}
\Big\{\, \E \trace\Big( \EstimationMat^{\!\!1/2} ( ( \Omega^{-1} +
\LinOp{\xi}^\T \NoiseCovariance^{-1} \LinOp{\xi})^{-1}
\EstimationMat^{\!\!1/2} \Big) : \nonumber
\\ &\qquad\qquad\qquad\qquad\qquad\qquad\qquad \Omega \succ
0,\trace(\ConstraintMat^{\!-1/2} \Omega \ConstraintMat^{\!-1/2})\leq
\radius^2 \,\Big\} \\ &= \Big(\frac{\radius + \lambda}{\radius}\Big)^2
\, f(0).
\end{align*}
Hence we have established the sandwich relation
\begin{align*}
f(0) \leq f(\lambda) \leq \Big(\frac{\radius + \lambda}{\radius}\Big)^2 \, f(0), \qquad \mbox{for all}~\lambda > 0.
\end{align*}
Note that $f(0) \leq f(\lambda') \leq f(\lambda)$ whenever $0 < \lambda' \leq \lambda$. 
Thus, $\inf_{\lambda > 0} = \lim_{\lambda \to 0^+} f(\lambda) = f(0)$, which establishes~\eqref{eqn:suff-condition-eqn}, 
completing the proof of the claim.

\subsection{Deferred proofs from Section~\ref{sec:proof-lower}}\label{app:proofs-lower}

In this section, we collect proofs of the results underlying the argument 
establishing our lower bound in Section~\ref{sec:proof-lower} of the paper. 

\subsubsection{Proof of Lemma~\ref{lem:minimax-reduction}} \label{app:proof-minimax-reduction}
By parameterizing $\theta^\star = \EstimationMat^{\!-1/2} \eta^\star$, we have  
\begin{align}
\minimaxrisk^{\rm G}(&T, \bP, \NoiseCovariance, \radius, \ConstraintMat, \EstimationMat)  \nonumber \\
&= \inf_{\etahat} \sup_{\eta^\star \in \Theta( \radius^2 \EstimationMat^{\!1/2} \ConstraintMat \EstimationMat^{\!1/2} )} 
\E_{\xi, w\sim\Normal{0}{I_\numobs}} 
\Big[\vecnorm[\big]{\etahat(\LinOp{\xi} \EstimationMat^{\!-1/2}, \LinOp{\xi} \EstimationMat^{\!-1/2} \eta^\star + \NoiseCovariance^{1/2} w) - \eta^\star}_2^2\Big] \nonumber \\ 
&=\inf_{\etahat} \sup_{\eta^\star \in \Theta( \radius^2 \EstimationMat^{\!1/2} \ConstraintMat \EstimationMat^{\!1/2} )} 
\E_{\xi, z \sim\Normal{0}{I_{r(\xi)}}} 
\Big[\vecnorm[\big]{\etahat(Q_\xi, Q_\xi \eta^\star + V_{\xi} \Lambda_\xi^{\!1/2} z) - \eta^\star}_2^2\Big] \label{eqn:sufficient-statistic}\\ 
&=\inf_{\etahat} \sup_{\eta^\star \in \Theta( \radius^2 \EstimationMat^{\!1/2} \ConstraintMat \EstimationMat^{\!1/2} )}  \E_{\omega \sim \tilde \bP, z \sim\Normal{0}{I_{r(\xi)}}} 
\Big[\vecnorm[\big]{\etahat(\omega, V_{\xi} V_{\xi}^\T \eta^\star + V_{\xi} \Lambda_\xi^{\!-1/2} z) - \eta^\star}_2^2\Big] \label{eqn:sufficient-statistic-2}\\ 
&= \minimaxrisk^{\rm G}_{\rm red}(\tilde \bP, \radius^2 \EstimationMat^{\!1/2} \ConstraintMat \EstimationMat^{\!1/2}). \nonumber 
\end{align}

We justify some of the relations in the display above. 
Since the density of $v = \LinOp{\xi} \EstimationMat^{\!-1/2} \eta^\star + \NoiseCovariance^{1/2} w$ is, up to constants independent of $\eta^\star$, proportional to
\begin{align*}
\exp\Big(-\frac{1}{2}\big\{ \langle \eta^\star, \EstimationMat^{\!-1/2} \LinOp{\xi}^\T \NoiseCovariance^{-1} \LinOp{\xi} \EstimationMat^{\!-1/2} \eta^\star \rangle 
- 2 \langle v, \NoiseCovariance^{-1} \LinOp{\xi} \EstimationMat^{\!-1/2} \eta^\star\rangle \big\}\Big),
\end{align*}
factorization arguments imply $Q_\xi \defn \EstimationMat^{\!-1/2}
\LinOp{\xi}^\T \NoiseCovariance^{-1}
\LinOp{\xi}\EstimationMat^{\!-1/2}$ and $v' \defn
\EstimationMat^{\!-1/2} \LinOp{\xi}^\T \NoiseCovariance^{-1} v$ are
sufficient statistics for $\eta^\star$. Note that $v'$ is distributed
$\Normal{Q_\xi \eta^\star}{Q_\xi}$. Thus, as consequence of the
Rao-Blackwell theorem, any minimax optimal estimator is a function of
$(Q_\xi, v')$, and hence display~\eqref{eqn:sufficient-statistic}
follows. Similarly, any optimal estimator function is a function of
any bijective function of $(Q_\xi, v')$. Evidently one can construct
$Q_\xi$ from $\omega \defn (r(\xi), V_{\xi}, \Lambda_\xi)$, and vice
versa.  On the other hand, $v'$ lies in the range of $G(\xi) \defn
\EstimationMat^{\!-1/2} \LinOp{\xi}^\T \NoiseCovariance^{-1/2}$, which
is the same as the range of $G(\xi) G(\xi)^\T = Q_\xi$; consequently
one may replace $v'$ with $Q_\xi^\dagger v' \equiv V_{\xi}
(\Lambda_\xi)^{-1} V_{\xi}^\T v'$, which is distributed $\Normal{V_{\xi}
  V_{\xi}^\T \eta^\star}{V_{\xi} (\Lambda_\xi)^{-1} V_{\xi}^\T}$, and so
that display~\eqref{eqn:sufficient-statistic-2} follows.

\subsubsection{Proof of~\Cref{lem:lambda-lower}}
\label{app:proof-lambda-lower}

In this argument, we use the notation $B(\etahat, \pi \mid \omega)$
to denote the Bayes risk of estimator $\etahat$, conditional on
$\omega$, for the original observation $\Upsilon$. Formally, it is the
expectation $\E[\vecnorm{\etahat(\Upsilon) - \eta}_2^2]$, where the
expectation is over $\Upsilon \sim \Normal{VV^\T}{V \Lambda^{-1}
  V^\T}$.

The main observation is that if we consider the projection of
$\Upsilon_\lambda$ onto the range of $V$, we will recover a random
variable with the same distribution as $\Upsilon$, and therefore the
risks are the same.  Formally, let $\etahat$ be any estimator which
is constant over the fibers of the operator $V V^\T$. Equivalently, it
can be written
\begin{align*}
\etahat(y) = \etahat_0(VV^\T y), \qquad \mbox{for some
  measurable}~\hat\eta_0.
\end{align*}
Let this class of estimators be denoted by $\cE_V$. Then we evidently have 
\begin{equation}\label{eqn:upper-lambda}
B_\lambda(\pi \mid \omega) \leq \inf_{\etahat\in \cE_V} B_\lambda(\etahat, \pi \mid \omega). 
\end{equation}
To complete the proof of the claim, we claim that
\begin{equation}\label{eqn:main-equality}
B_\lambda(\etahat, \pi \mid \omega) = B(\etahat, \pi \mid \omega), \qquad \mbox{for any}~\etahat\in \cE_V.
\end{equation}
This follows immediately from the fact that $VV^\T \Upsilon_\lambda = \Upsilon$ with probability 1. 
We note that combination with~\eqref{eqn:upper-lambda} furnishes the claim, since it 
implies that 
\begin{align*}
B_\lambda(\pi \mid \omega) \leq \inf_{\etahat\in \cE_V} B(\etahat, \pi \mid \omega) = B(\pi \mid \omega).
\end{align*}
The final equality occurs since for any measurable estimator $\etahat\not \in \cE_V$, we can define 
$\etahat_{V}(y) = \etahat(VV^\T y)$, and since $\Upsilon = VV^\T \Upsilon$ with probability $1$, 
and therefore $B(\etahat_V, \pi \mid \omega) = B(\etahat, \pi \mid \omega)$, which establishes this claim. 

\subsubsection{Proof of Lemma~\ref{lem:lambda-bayes-risk}}
\label{app:proof-lambda-bayes-risk}
Let $\etahat_\pi$ denote the posterior mean $y \mapsto \E[\eta \mid
  \Upsilon_\lambda = y]$.  Then, as the posterior mean $\etahat_\pi$
minimizes the Bayes risk $\etahat \mapsto B_\lambda(\etahat, \pi \mid
\omega)$ over all measurable estimators $\etahat$, it suffices to
compute the risk of $\etahat_\pi$.  Note that, by definition of
conditional expectation, we have
\begin{align*}
\etahat_\pi(y) = \frac{1}{p(y)}\int \eta \, p(y \mid \eta) \, \pi(\ud \eta).
\end{align*}
We now compute the derivative of $p(y)$. Exchanging integration and differentiation,\footnote{
  This is valid since $y \mapsto p(y\mid \eta)$ is differentiable for each $\eta$, 
  and for each $y$, we have  $\eta \mapsto p(y \mid \eta)$ and $\eta \mapsto \nabla_y p(y \mid \eta) = \Sigma_\lambda^{-1}(X_\lambda \eta - y)$ 
  are $\pi$-integrable (since $0 \leq p(y \mid \eta) \leq 1$, and the gradient is an affine function of $\eta$). 
  }
\begin{align*}
 \Sigma_\lambda \nabla p(y) = \int (X_\lambda \eta - y) \, p(y \mid \eta) \, \pi(\ud \eta).
\end{align*}
Therefore, we conclude that 
\begin{align*}
\etahat_\pi(y) = X_\lambda^{-1} \Big(y + \Sigma_\lambda \nabla \log p(y)\Big).
\end{align*}
Finally, to compute risk of the posterior mean $\etahat_\pi(\Upsilon_\lambda) \defn \E[\eta \mid \Upsilon_\lambda]$, 
we add and subtract the observation $X_\lambda^{-1} \Upsilon_\lambda$, and find that
\begin{align*}
\E_{(\eta, \Upsilon_\lambda)}\Big[(\eta -
  \etahat_\pi(\Upsilon_\lambda)) \otimes (\eta -
  \etahat_\pi(\Upsilon_\lambda))\Big] =X_\lambda^{-1} \Sigma_\lambda
X_{\lambda}^{-1} - X_\lambda^{-1} \Sigma_\lambda \E[\nabla \log
  p(\Upsilon_\lambda) \otimes \nabla \log p(\Upsilon_\lambda)]
\Sigma_\lambda X_{\lambda}^{-1}.
\end{align*}
Identifying the Fisher information in the display above, factoring the expression, and taking the trace yields the desired result. 

\subsubsection{Proof of Lemma~\ref{lem:fisher-info-lower}}\label{app:proof-fisher-info-lower}
Note that $\pi_{\tau, \Pi}$ is evidently absolutely continuous with respect to Lebesgue measure. 
In particular, on the interior of $\Theta(K)$, $\pi_{\tau, \Pi}$ and $\pi_{\tau, \Pi}^{\rm G}$ have the 
same Lebesgue density up to rescaling by $\pi_{\tau, \Pi}^{\rm G}(\Theta(K))$. Denote this density by $f_{\tau, \Pi}$. 
Therefore, we have 
\begin{align*}
\Information{\pi_{\tau, \Pi}^{\rm G}} &= \E_{\eta \sim \pi_{\tau, \Pi}^{\rm G}} \1_{\Theta(K)}(\eta) \nabla \log 
f_{\tau, \Pi}(\eta) \otimes \nabla \log f_{\tau, \Pi}(\eta) + 
\E_{\eta \sim \pi_{\tau, \Pi}^{\rm G}} \1_{\Theta(K)^c}(\eta) \nabla \log 
f_{\tau, \Pi}(\eta) \otimes \nabla \log f_{\tau, \Pi}(\eta) \\ 
&\succeq \E_{\eta \sim \pi_{\tau, \Pi}^{\rm G}} \1_{\Theta(K)}(\eta) \nabla \log 
f_{\tau, \Pi}(\eta) \otimes \nabla \log f_{\tau, \Pi}(\eta)  \\
&= \pi_{\tau, \Pi}^{\rm G}(\Theta(K)) \Information{\pi_{\tau, \Pi}}.
\end{align*}
The final equality arises since the boundary of $\Theta(K)$ has Lebesgue measure zero. 
Using the well known relation $\Information{\pi_{\tau, \Pi}^{\rm G}} = (\tau^2 \Pi)^{-1}$~\cite[Example 6.3]{LehCas98},
the above display implies that 
\begin{align*}
\Information{\pi_{\tau, \Pi}}^{-1} \succeq \pi_{\tau, \Pi}^{\rm G}(\Theta(K)) \tau^2 \Pi 
= \tau^2 (1 - \pi_{\tau, \Pi}^{\rm G}(\Theta(K)^c)) \Pi.
\end{align*}
To ensure that $\eta \sim \pi_{\tau, \Pi}^{\rm G}$ lies in $\Theta(K)$ with decent probability, we take 
$\Pi$ to satisfy the relation $\trace(K^{-1} \Pi) \leq 1$. Then defining 
\begin{align*}
c(\tau, \Pi) \defn \tau^2 (1 - \pi_{\tau, \Pi}^{\rm G}(\Theta(K)^c)), 
\end{align*}
completes the proof of the claim. 

\subsubsection{Proof of Lemma~\ref{lem:constant-lower}}\label{app:proof-constant-lower}
Fix $\Pi \succ 0$ such that $\trace(\Pi^{1/2} K^{-1} \Pi^{1/2}) \leq 1$. 
Let $\lambda = (\lambda_1,\dots, \lambda_d)$ denote the eigenvalues of $\Pi^{1/2} K^{-1} \Pi^{1/2}$.
The vector satisfies the inequalities $\lambda \succ 0, \lambda^\T \1 \leq 1$. 
Moreover, by the rotational invariance of the Gaussian, we have for $g \sim \Normal{0}{I_\dimension}$, that
\begin{align*}
\pi^{\rm G}_{\tau, \Pi}(\Theta(K)^c) = \P\Big\{ \tau^2 g^\T \Pi^{1/2} K^{-1} \Pi^{1/2} g > 1\Big\} 
= \P\Big\{\tau^2 \sum_{i=1}^\dimension \lambda_i g_i^2 > 1 \Big\}.
\end{align*}
Let us make the choice $\tau^2 = 1/2$. 
Then, note for any $\lambda \succ 0, \lambda^\T \1 \leq 1$, by Markov's inequality,
\[
\P\Big\{\sum_{i=1}^\dimension \lambda_i g_i^2 > 2 \Big\}
\leq \frac{\sum_{i=1}^\dimension
  \lambda_i \E[g_i^2]}{2} = \frac{1}{2}.
\]
Hence, using this bound in the definition of $c(\tau, \Pi)$, we find
\begin{align*}
c_\ell(K) \geq \inf_{\lambda \succ 0, \lambda^\T \1 \leq 1} 
c(1/2, \diag(\lambda)) \geq \frac{1}{4},
\end{align*}
which completes the proof of the claim.

\bibliographystyle{abbrv}
\bibliography{references}
\end{document}